\documentclass{article}
\usepackage{amssymb,latexsym,pstricks
}

\newtheorem{theorem}{Theorem}[section]

\newtheorem{corollary}{Corollary}[section]
\newtheorem{definition}{Definition}[section]

\newtheorem{remark}{Remark}[section]
\newtheorem{lemma}{Lemma}[section]

\newtheorem{example}{Example}[section]

\def\R{\mathbb{R}}

\def\CP{\mathbb{CP}}
\def\BZ{\mathrm{BZ}}
\def\({\left(}
\def\){\right)}
\def\[{\left[}
\def\]{\right]}
\def\e{{\bf e}}
\def\v{{\bf v}}
\def\w{{\bf w}}
\def\z{{\bf z}}
\def\s{\mbox{s}}
\def\res{\mbox{Res}}
\def\<{\langle}
\def\>{\rangle}

\newcommand{\Gl}{\mbox{GL}}

\def\A{{\bf A}}
\def\B{{\bf B}}
\def\C{{\bf C}}

\def\tA{{\bf {\tilde A}}}
\def\tB{{\bf {\tilde B}}}
\def\tC{{\bf {\tilde C}}}
\def\ta{{\tilde a}}
\def\tb{{\tilde b}}
\def\tc{{\tilde c}}
\def\Id{\mbox{{\bf Id}}}
\def\trace{\mbox{tr}}
\def\proof{\smallskip\noindent {\it Proof --- \ }}
\def\endproof{\hfill$\square$\medskip}

\title{Some Explicit Solutions of the Additive Deligne-Simpson Problem
and Their Applications (Preprint)}

\author{%
Oleg Gleizer\\[.05in]
{\normalsize University of California, Los
Angeles~~~~~~~~~~~~~~~~~~~~~~~~~~~~~~~~~~~~~~~~~~~
oleg@math.ucla.edu}\\[.2in]}

\date{{\normalsize April 22, 2002}}

\begin{document}
\maketitle

\begin{abstract}

    In this paper we construct three infinite series and two extra
triples ($E_8$ and ${\hat E}_8$) of complex matrices $\B$, $\C$,
and $\A=\B+\C$ of special spectral types associated to
C.~Simpson's classification in \cite{S} and P.~Magyar, J.~Weyman,
A.~Zelevinsky classification in \cite{MWZ}. This enables us to
construct Fuchsian systems of differential equations which
generalize the hypergeometric equation of Gauss-Riemann. In a
sense, they are the closest relatives of the famous equation,
because their triples of spectral flags have finitely many orbits
for the diagonal action of the general linear group in the space
of solutions. In all the cases except for $E_8$, we also
explicitly construct scalar products such that $\A$, $\B$, and
$\C$ are self-adjoint with respect to them. In the context of
Fuchsian systems, these scalar products become monodromy invariant
complex symmetric bilinear forms in the spaces
of solutions. \\

    When the eigenvalues of $\A$, $\B$, and $\C$ are real, the
matrices and the scalar products become real as well. We find
inequalities on the eigenvalues of $\A$, $\B$, and $\C$ which make
the scalar products positive-definite. \\

    As proved by A.~Klyachko, spectra of three hermitian (or real
symmetric) matrices $\B$, $\C$, and $\A=\B+\C$  form a polyhedral
convex cone in the space of triple spectra. He also gave a
recursive algorithm to generate inequalities describing the cone.
The inequalities we obtain describe non-recursively some faces of
the Klyachko cone.

\end{abstract}

\section{Introduction}
\label{sec:intro}

 Let $V$ be a vector space over complex numbers such that $\dim V=n$ where $1<n<\infty$.
Let $\B$, $\C$ be linear operators in $V$ and let $\A=\B+\C$. We
call the pair $\B,\C$ {\it irreducible} if the operators do not
preserve simultaneously any proper subspace of $V$. \\

 Let $O_{\A}$ be the adjoint orbit of $\A$ in $\mbox{End}\,V$  under the $\Gl(V)$
action. We call the triple $\A$, $\B$, $\C$ {\it rigid}, if any
other triple $\B'$, $\C'$, $\A'=\B'+\C'$ such that $\A'\in
O_{\A}$, $\B'\in O_{\B}$, and
$\C'\in O_{\C}$ is conjugate to the triple $\A$, $\B$, $\C$. \\

 For a linear operator ${\bf A}\in \mbox{End} \, V$, we call the multiset
of its eigenvalues the {\it spectrum of ${\bf A}$}. This means
that each eigenvalue $\lambda_i$ is taken with its multiplicity
$m_i$. Any ordering $\lambda_1$, $\lambda_2, \cdots ,\lambda_k$ of
distinct eigenvalues of $\A$ allows us to represent the spectrum
of $\A$ by a vector $\s(\A) = ( \underbrace{\lambda_1 \cdots
\lambda_1}_ {m_1~\mbox{{\it times}}},\underbrace{\lambda_2 \cdots
\lambda_2}_{m_2~\mbox{{\it times}}}, \cdots ,
\underbrace{\lambda_k \cdots \lambda_k}_{m_k~\mbox{{\it times}}} )
\in \mathbb{C}^n$. For a diagonalizable operator $\A$, we call the
partition $(m_1,m_2,\cdots ,m_k)$ of $n$ the {\it spectral type}
of $\A$. With slight abuse of terminology, we also call the
spectral type of $\A$ any composition obtained by some ordering of
$\lambda_1,\cdots ,\lambda_k$. We say that a vector $(x_1,\cdots
x_n,y_1,\cdots y_n,z_1,\cdots z_n)\in (\mathbb{C}^n)^3$ satisfies
the {\it trace condition} if $\sum_{i=1}^n x_i=\sum_{i=1}^n
(y_i+z_i)$. Then $(\s(\A),\s(\B),\s(\C))$ belongs to the
hyperplane in $(\mathbb{C}^n)^3$ defined by the trace condition.
We call this hyperplane the {\it space of triple spectra}. Let
$\alpha=(m_1,m_2,\cdots ,m_p)$, $\beta=(n_1,n_2,\cdots ,n_q)$, and
$\gamma=(k_1,k_2,\cdots ,k_r)$ (compositions of $n$) be the
spectral types of $\A$, $\B$, and $\C$. Then
$\(\s(\A),\s(\B),\s(\C)\)$ lies in the part
$S(\alpha,\beta,\gamma) \subset \mathbb{C}^{3n}$ defined as
follows. A vector $(x,y,z)\in \mathbb{C}^{3n}$ is in
$S(\alpha,\beta,\gamma)$ if $x_1=x_2=\cdots =x_{m_1} \ne
x_{m_1+1}=\cdots =x_{m_1+m_2}\ne \cdots$ and the same for
$y$ and $z$. \\

 Consider the following table of triples of spectral types.

\begin{equation}
\label{equation:S_list}
\begin{tabular}{|l|l|l|}
\hline
 & & \\
hypergeometric family & $(1,m-1),(1^m),(1^m)$ & $m\ge 2$ \\
 & & \\
\hline
 & & \\
even family & $(m,m),(1,m-1,m),(1^{2m})$ & $m\ge 2$ \\
 & & \\
\hline
 & & \\
odd family & $(m+1,m),(1,m,m),(1^{2m+1})$ & $m\ge 2$ \\
 & & \\
\hline
 & & \\
extra case & $(4,2),(2,2,2),(1^6)$ & \\
 & & \\
\hline
\end{tabular}
\end{equation}
\medskip

\noindent Here and later $(1^n)$ is a shorthand for
$(\underbrace{1,1,\cdots ,1}_
{n~times})$. \\

\begin{theorem}[Simpson, Kostov]
\label{thm:Slist}

Let $(\alpha,\beta,\gamma)$ be a triple of spectral types such
that at least one of them is $(1^n)$. The following conditions are
equivalent:
\begin{enumerate}
\item for a generic point $(x,y,z)\in S(\alpha,\beta,\gamma)$ there exists a
rigid irreducible triple \\
$(\A=\B+\C,\B,\C)$ of diagonalizable operators such that
$(\s(\A),\s(\B),\s(\C))=(x,y,z)$;
\item $(\alpha,\beta,\gamma)$ is one of the triples in (\ref{equation:S_list}).
\end{enumerate}
\end{theorem}

\begin{remark}
\label{remark:Kostov} {\rm This theorem is an additive version of
Theorem 4 from \cite{S}. This version easily follows from
Simpson's results. A more elementary proof of Theorem 4 from
\cite{S} and a proof of Theorem \ref{thm:Slist} were given by
V.~Kostov in \cite{Ko1}. }
\end{remark}

 The first main result of this paper is that for each triple $(\alpha,\beta,\gamma)$
of spectral types from (\ref{equation:S_list}) and a generic
vector from $S(\alpha,\beta,\gamma)$ we explicitly construct the
corresponding triple $(\A,\B,\C)$. \\

    Recently there appeared algorithms to produce all rigid
irreducible $r$-tuples of matrices $M_1,\cdots ,M_r$ such that
$M_1+\cdots +M_r=0$, see \cite{CB1} and \cite{DR}. We use a
different (less general, but more powerful for our particular
purposes) tool: P.~Magyar, J.~Weyman, and A.~Zelevinsky in
\cite{MWZ} constructed all indecomposable triple partial flag
varieties with finitely many orbits for the diagonal action of the
general linear group. Their list (\ref{equation:MWZ_list}) is
strikingly similar to the list (\ref{equation:S_list}) of Simpson.
It has just one more family: the $E_8$-family. A triple of {\it
spectral flags} (for the definition, see Section
\ref{sec:triple_flag_var}) of the matrices $\A$, $\B$, and $\C$
constructed from Simpson's list (\ref{equation:S_list}) gives a
representative of the open orbit of the corresponding triple flag
variety from (\ref{equation:MWZ_list}). \\

    Our $\A$, $\B$, and $\C$ have the following common features.
$\B$ is block upper-triangular, $\C$ is block lower-triangular.
The block sizes of $\B$ and $\C$ are given by the compositions
$\beta$ and $\gamma$ respectively. Each entry of $\A$, $\B$, and
$\C$ is a ratio of products of linear forms in the eigenvalues of
$\A$, $\B$, and $\C$. The coordinates of all eigenvectors of $\A$,
$\B$, and $\C$ are also ratios of products of linear forms in the
eigenvalues. The linear forms are remarkably simple: all the
coefficients are equal to either $1$ or $-1$. As a corollary of
our construction, we obtain the following.

\begin{theorem}
\label{thm:form} For every composition $(\alpha, \beta, \gamma)$
from Simpson's list (\ref{equation:S_list}), there exist open
subsets $S''(\alpha, \beta, \gamma)\subset S'(\alpha, \beta,
\gamma) \subset S(\alpha, \beta, \gamma)$ with the following
properties.

\begin{enumerate}
\item Each of $S'(\alpha, \beta, \gamma)$ and $S''(\alpha, \beta, \gamma)$ is
obtained from $S(\alpha, \beta, \gamma)$ by removing finitely many
hyperplanes.

\item If $\(\s(\A), \s(\B), \s(\C)\) \in S'(\alpha, \beta, \gamma)$,
then there exists a non-zero symmetric bilinear form on $V$ such
that $\A$, $\B$, and $\C$ are self-adjoint with respect to it.

\item If $\(\s(\A), \s(\B), \s(\C)\) \in S''(\alpha, \beta, \gamma)$,
then there exists a non-degenerate symmetric bilinear form on $V$
such that $\A$, $\B$, and $\C$ are self-adjoint with respect to
it.
\end{enumerate}
\end{theorem}

\noindent This theorem is proved case by case in Theorems
\ref{thm:hg:form}, \ref{thm:even:form}, \ref{thm:odd:form},
\ref{thm:extra:form} for the bilinear forms given by the formulas
(\ref{equation:hg:form}), (\ref{equation:even:form}),
(\ref{equation:odd:form}), (\ref{equation:extra:form})
correspondingly.

\begin{remark}
\label{remark:thm:form} The main virtue of this theorem is not the
proof of existence of the objects, but an explicit construction of
all of them.
\end{remark}

    In Simpson's list (\ref{equation:S_list}), the last composition
$\gamma$ is always $(1^n)$. Thus, the matrix $\C$ has all
eigenvalues distinct. Let $\v_i$ be the eigenvector of $\C$
corresponding to the eigenvalue $c_i$. If $\C$ is self-adjoint
with respect to a scalar product $\<*,*\>$ on $V$, then
$\<\v_i,\v_j\>=l_i\, \delta_{ij}$. If we manage to find $l_i$ such
that the matrix $\B$ becomes self-adjoint with respect to
$\<*,*\>$ as well, then $\A$ is also self-adjoint with respect to
$\<*,*\>$ as $\A=\B+\C$. We find the $l_i$ and it turns out that
they are also ratios of products of linear forms in the
eigenvalues of $\A$, $\B$, and $\C$. And again all the
coefficients of the forms are equal to either $1$ or $-1$. The set
of linear forms that appear in the $l_i$ includes the set of
linear forms that appear in the matrix elements of $\A$, $\B$, and
$\C$ and in the coordinates of their eigenvectors. The hyperplanes
one has to remove from $S(\alpha, \beta, \gamma)$ to obtain
$S''(\alpha, \beta, \gamma)$ of Theorem \ref{thm:form} are exactly
the zero levels of the linear forms that appear in the $l_i$. The
explicit formulas we find for the $l_i$ give explicit description of
these hyperplanes. \\

    Probably the most important applications of our explicit construction
is to {\it Fuchsian systems} (see Section
\ref{sec:applications_to_Fuchsian_systems} for the definition).
Let $z_1$, $z_2$, $z_3$ be distinct points of $\CP^1$. Consider
the following system of differential equations

\begin{equation}
\label{equation:fs}
\frac{df}{dz}=\[\frac{\B}{z-z_2}+\frac{\C}{z-z_3}-\frac{\A}{z-z_1}\]f(z)
\end{equation}

\noindent where $\A=\B+\C$, $z\in\CP^1\setminus \{z_1,z_2,z_3\}$
and $f$ takes values in $V$. The matrices $A$, $B$, and $C$ are
called the {\it residue matrices} of (\ref{equation:fs}). Their
eigenvalues are called {\it local exponents}. Real parts of the
local exponents are the rates of growth of solutions of
(\ref{equation:fs}) expanded analytically towards the
corresponding singularities (and restricted to sectors centered at the singularities).
Thus, at each singularity the space
of solutions stratifies into a flag. Local basis changes near each
singularity turn this flag into a flag variety. The triple of flag
varieties of the Gauss-Riemann equation has finitely many orbits
for the action of the general linear group in the space of
solutions. The Fuchsian systems constructed by means of our
matrices exhaust the list of Fuchsian systems (with more than two
singularities) having this property. In this sense, they are the
simplest Fuchsian systems possible and we expect their solutions
to be interesting functions. \\

    It was known to F.~Klein that if the hypergeometric equation of
Gauss-Riemann had real local exponents, then there existed a
monodromy invariant hermitian form in the space of solutions. If
the local exponents were generic, then the form was non-degenerate
and unique up to a real constant multiple. We prove the same for
all the Fuchsian systems constructed from (\ref{equation:S_list}).
Indeed, when all the eigenvalues of $\A$, $\B$, and $\C$ are real,
the form $\<*,*\>$ becomes real as well. So do the matrices $\A$,
$\B$, and $\C$ themselves. Thus, $\A$, $\B$, and $\C$ become
matrices of real operators acting on the real space $V_{\R}$ and
self-adjoint with respect to the real symmetric bilinear form
$\<*,*\>_{\R}$. Let us extend the form $\<*,*\>_{\R}$ to the
hermitian form $(*,*)$ on $V$. This form gives rise to the
monodromy invariant hermitian form in the space of solutions of
(\ref{equation:fs}). Once again, the forms are constructed
explicitly. For the hypergeometric family, this result is not new.
The Fuchsian systems from the hypergeometric family are equivalent
to the generalized hypergeometric equations studied by F.~Beukers
and G.~Heckman in \cite{BH}. Among other things, they construct
the hermitian form. Also the generalized hypergeometric equations
were studied in what became later known as the {\it Okubo normal
form} by K.~Okubo in \cite{O}. For the generalized hypergeometric
equations in the Okubo normal form, the monodromy invariant
hermitian form was constructed by Y.~Haraoka in \cite{H2}. \\

    As proved by A.~Klyachko in \cite{Kl}, if a hermitian form is
positive definite, then the spectra of hermitian matrices $\B$,
$\C$ and $\A=\B+\C$ form a polyhedral convex cone in the space of
triple spectra. His proof contains a recursive algorithm to
compute the inequalities describing the cone. We call this cone
the {\it Klyachko cone} and we call the inequalities the {\it
Klyachko inequalities}. Beukers and Heckman in \cite{BH} give
explicitly the inequalities on the local exponents of the
generalized hypergeometric equation which make the monodromy
invariant hermitian form in the space of solutions positive
definite. Thus, they describe non-recursively a non-trivial face
of the Klyachko cone. We do the same for all the Fuchsian systems
constructed from (\ref{equation:S_list}). Hence the second
important application of our results is an explicit description of
some interesting faces of the Klyachko cone. Beukers and Heckman
in \cite{BH} use their criterion of positivity of the hermitian
form to see when solutions of the generalized hypergeometric
equations are algebraic functions. It is also needed to know the
signature of the form for applications to number theory, see
\cite{BH} and \cite{DR}. Our construction provides tools to answer
similar questions about the solutions of our Fuchsian systems. \\

    Let $\lambda$, $\mu$, and $\nu$ be highest weights of $\Gl(V)$.
Let $V_{\lambda}$, $V_{\mu}$, and $V_{\nu}$ be the corresponding
rational irreducible representations of $\Gl(V)$. Let $V_{\lambda}
\otimes V_{\mu}=\sum_{\nu} c_{\lambda \mu}^{\nu} V_{\nu}$ be the
decomposition of the tensor product of $V_{\lambda}$ and $V_{\mu}$
into the sum of irreducible representations. It follows from the
results of A.~Klyachko \cite{Kl} combined with a refinement by
A.~Knutson and T.~Tao \cite{KT} that the lattice points of the
Klyachko cone are precisely the triples of weights $\lambda$,
$\mu$, $\nu$ with non-zero Littlewood-Richardson coefficient
$c_{\lambda \mu}^{\nu}$. Thus, our techniques allow us to
explicitly describe some triples of highest weights with non-zero
Littlewood-Richardson coefficients. In fact, for all the cases we
consider, the Littlewood-Richardson coefficients are equal to $1$. \\

    The paper is organized as follows. In Section
\ref{sec:main_results}, we formulate main results of the paper.
Namely, we list the triples $(\A,\B,\C)$, the scalar products, and
the Klyachko inequalities for all the partitions from Simpson's
list (\ref{equation:S_list}). In Section \ref{sec:proofs},
theorems of Section \ref{sec:main_results} are proved and
elaborated. For example, in subsection
\ref{subsection:HGM_from_EM} we construct the residue matrices of
the $m$-hypergeometric system as submatrices of the residue
matrices of the $m$-even system. \\

    Although the actual construction of the matrices relied heavily on
the explicit description of representatives of the open orbits in
\cite{MWZ}, it turned out that once the answers were known, it was
much simpler to prove them by inspection. We start using the
results of \cite{MWZ} directly only in Section
\ref{sec:triple_flag_var}. In the Section, we construct the
matrices $\A$, $\B$, and $\C$ which give rise to the $E_8$-family
of Magyar, Weyman, and Zelevinsky. We also prove the following
theorem.

\begin{theorem}
\label{thm:irr} Let $(\alpha,\beta,\gamma)$ be a triple of
spectral types from Theorem \ref{thm:Slist}. If $\(\s(\A), \s(\B),
\s(\C)\)$ is a point of $S''(\alpha, \beta, \gamma)$ then the
triple $(\A,\B,\C)$ is irreducible.
\end{theorem}

\noindent The proof of the Theorem for the hypergeometric family
is on page \pageref{proof:hg:irr}, for all other families -- on
page \pageref{proof:irr}.

\begin{remark}
\label{remark:genericity} The genericity condition of Theorem
\ref{thm:Slist} was somewhat a mysterious one. In general, it still is.
A theorem of Katz (see \cite{NKatz}) excludes the coexistence of
irreducible and reducible triples in rigid cases.
If one deals with a reducible triple, then except for the ``big''
trace condition one also has a ``small'' trace condition coming from
the reduced submatrices. Thus, people call generic spectra that stay
away from all ``small trace condition'' hyperplanes possible (see Kostov's
papers). In our cases however, Theorem \ref{thm:irr} gives an explicit meaning
to the genericity condition: ``generic'' means ``not in $S''$''.
\end{remark}

\begin{remark}
If we take a
different ordering of the eigenvalues of $\A$, $\B$, and $\C$, then the scalar
products $\<\v_i,\v_i\>$ change. So does the set $S'$, but the set $S''$ does not.
\end{remark}

    Let $\A$, $\B$, and $\C$ be self-adjoint with respect to a
non-zero symmetric bilinear form on $V$. Let
$(\alpha,\beta,\gamma)$ be a triple of spectral types from Theorem
\ref{thm:Slist}. Then the following corollary of Theorem
\ref{thm:irr} strengthens the third statement of Theorem
\ref{thm:form}.

\begin{corollary}
\label{corollary:uniqueness_of_form} If $(\s(\A),\s(\B),\s(\C))$
is a point of $S''(\alpha,\beta,\gamma)$, then the form is unique
up to a constant multiple.
\end{corollary}

\proof If $(\s(\A),\s(\B),\s(\C))$ is a point of
$S''(\alpha,\beta,\gamma)$, then it follows from Theorem
\ref{thm:irr} that the triple $(\A,\B,\C)$ is irreducible. If a
triple $(\A,\B,\C)$ is irreducible, then uniqueness of the form
follows from Schur's lemma. \endproof \\

    In Section \ref{sec:LR}, we introduce the {\it
Berenstein-Zelevinsky triangles} which provide a geometric version
of the celebrated Littlewood-Richardson rule. For the $E_8$
family, we do not have formulas for the hermitian form as nice as
we have for other families. However, the Berenstein-Zelevinsky
triangles enable us to compute the Klyachko inequalities for the
$E_8$-family as well. \\

    Section \ref{sec:applications_to_Fuchsian_systems} contains no new
results. In the section, we provide (very) basic facts about
Fuchsian systems and raise questions we plan to answer in
subsequent publications. In particular, we quote some results from
\cite{Y1}, \cite{H1}, and \cite{H2} which are very similar to (but
different from) ours. \\

    Most of the proofs of the paper boil down to proofs of certain
rational identities. These identities are collected in Section
\ref{sec:appendix} (the Appendix).

\section*{Acknowledgments}

    I would like to thank G.~Heckman who acquainted me with the
problem when I was a student at the Master Class Program of the
University of Utrecht. I would like to thank W.~Crawley-Boevey,
H.~Derksen, Y.~Haraoka, A.~Khovansky, A.~Klyachko, A.~Knutson,
V.~Kostov, P.~Magyar, A.~Postnikov, and C.~Simpson for answering
a lot of questions. \\

    I would like to thank the referee for very helpful suggestions
which made this paper more clear and easier to read. \\

    This paper was a basis for my Ph.D. thesis at Northeastern University
of Boston. I would especially like to thank my adviser
A.~Zelevinsky for his tremendous support which made both the paper
and the thesis possible. This paper was partially supported by NSF
through research grants of A.~Zelevinsky.

\section{Main Results}
\label{sec:main_results}

\subsection{Hypergeometric Family}
\label{subsec:results:hg}

 Let us pick a vector $(a_1,\underbrace{a_2,\cdots ,a_2}_{m-1~times},b_1,\cdots ,b_m,
c_1,\cdots ,c_m)$ from $S((1,m-1),(1^m),(1^m))$. Recall that this
means $a_1\ne a_2$, all $b_i$ are distinct, all $c_i$ are
distinct, and the trace condition holds:
$a_1+(m-1)a_2=\sum_{i=1}^m (b_i+c_i)$. Define the matrix elements
of $\B$ and $\C$ as follows:

\begin{equation}
\label{equation:hg:B,C}
\begin{array}{ll}
B_{ij}=\left\{
\begin{array}{lll}
0                 & ,if & i<j \\
b_i               & ,if & i=j \\
b_i+c_{m+1-i}-a_2 & ,if & i>j
\end{array},
\right. &

C_{ij}=\left\{
\begin{array}{lll}
b_i+c_{m+1-i}-a_2 & ,if & i<j \\
c_{m+1-i}         & ,if & i=j \\
0                 & ,if & i>j
\end{array}.
\right.
\end{array}
\end{equation}

\noindent Here is an example with $m=5$.

\begin{example}
\label{example:hg:B,C}

$$
\begin{array}{l}

\B=\left[
\begin{array}{lllll}
b_1 & b_1+c_5-a_2 & b_1+c_5-a_2 & b_1+c_5-a_2 & b_1+c_5-a_2 \\
0   & b_2         & b_2+c_4-a_2 & b_2+c_4-a_2 & b_2+c_4-a_2 \\
0   & 0           & b_3         & b_3+c_3-a_2 & b_3+c_3-a_2 \\
0   & 0           & 0           & b_4         & b_4+c_2-a_2 \\
0   & 0           & 0           & 0           & b_5
\end{array}
\right] \\
 \\
\C=\left[
\begin{array}{lllll}
c_5         & 0           & 0           & 0           & 0 \\
b_2+c_4-a_2 & c_4         & 0           & 0           & 0 \\
b_3+c_3-a_2 & b_3+c_3-a_2 & c_3         & 0           & 0 \\
b_4+c_2-a_2 & b_4+c_2-a_2 & b_4+c_2-a_2 & c_2         & 0 \\
b_5+c_1-a_2 & b_5+c_1-a_2 & b_5+c_1-a_2 & b_5+c_1-a_2 & c_1
\end{array}
\right]

\end{array}
$$
\end{example}

\noindent It is clear that $\s(\B)=\{b_1,\cdots ,b_m\}$ and
$\s(\C)=\{c_1,\cdots ,c_m\}$. Since all the $b_i$ and all the
$c_i$ are distinct, $\B$ and $\C$ are diagonalizable.

\begin{theorem}
\label{thm:hg:A}
If $\B$ and $\C$ are given by (\ref{equation:hg:B,C}), then $\A=\B+\C$ is diagonalizable and \\
$\s(\A)=\{a_1,\underbrace{a_2,a_2,\cdots
,a_2}_{m-1~\mbox{times}}\}$.
\end{theorem}

 For $i=1,\cdots ,m$, let $\v_i=(v_i^1,\cdots ,v_i^{i-1},1,0,\cdots ,0)$ be the
eigenvector of $\B$ with the eigenvalue $b_i$.

\begin{lemma}
\label{lemma:hg:eigenbasis_of_B} For $1\le j<i\le m$, we have
$$
v_i^j =
\displaystyle{\frac{b_j+c_{m+1-j}-a_2}{b_i-b_j}}\prod\limits_{k=1}^{i-j-1}
\displaystyle{\frac{b_i+c_{m+1-j-k}-a_2}{b_i-b_{j+k}}}
$$
Here and in the sequel, all empty products are understood to be
equal to $1$.
\end{lemma}

\noindent Here is an example with $m=5$.

\begin{example}
\label{example:hg:m=5:eigenvect_of_B}
$$
\begin{array}{l}
\v_1=\( 1,0,0,0,0 \) \\
 \\
\v_2=\( \frac{b_1+c_5-a_2}{b_2-b_1},1,0,0,0 \) \\
  \\
\v_3=\( \frac{(b_1+c_5-a_2)(b_3+c_4-a_2)}{(b_3-b_1)(b_3-b_2)},
           \frac{b_2+c_4-a_2}{b_3-b_2},1,0,0 \) \\
 \\
\v_4=\(
\frac{(b_1+c_5-a_2)(b_4+c_3-a_2)(b_4+c_4-a_2)}{(b_4-b_1)(b_4-b_2)(b_4-b_3)},
           \frac{(b_2+c_4-a_2)(b_4+c_3-a_2)}{(b_4-b_2)(b_4-b_3)},
           \frac{b_3+c_3-a_2}{b_4-c_3},1,0 \) \\
 \\
\v_5=\(
\frac{(b_1+c_5-a_2)(b_5+c_2-a_2)(b_5+c_3-a_2)(b_5+c_4-a_2)}
               {(b_5-b_1)(b_5-b_2)(b_5-b_3)(b_5-b_4)},
          \frac{(b_2+c_4-a_2)(b_5+c_2-a_2)(b_5+c_3-a_2)}{(b_5-b_2)(b_5-b_3)(b_5-b_4)},
          \frac{(b_3+c_3-a_2)(b_5+c_2-a_2)}{(b_5-b_3)(b_5-b_4)}, \right. \\
\phantom{\v_5=\( \right.} \left. \frac{b_4+c_2-a_2}{b_5-b_4},1 \)
\end{array}
$$
\end{example}

 We define a scalar product on $V$ by setting

\begin{equation}
\label{equation:hg:form} \<\v_i,\v_j\>=\delta_{ij}\,
\displaystyle{\frac{\prod\limits_{k=i+1}^m (b_i-b_k)}
{\prod\limits_{k=1}^{i-1} (b_i-b_k)}}\,
\displaystyle{\frac{\prod\limits_{k=m+2-i}^m (b_i+c_k-a_2)}
{\prod\limits_{k=1}^{m+1-i} (b_i+c_k-a_2)}}.
\end{equation}

\noindent Here is an example of the form with $m=5$.

\begin{example}
\label{example:hg:m=5:form}
$$
\begin{array}{ll}
\<\v_1,\v_1\> = & \frac{(b_1-b_2)(b_1-b_3)(b_1-b_4)(b_1-b_5)}
{(b_1+c_1-a_2)(b_1+c_2-a_2)(b_1+c_3-a_2)(b_1+c_4-a_2)(b_1+c_5-a_2)} \\
 & \\
\<\v_2,\v_2\> = & \frac{(b_2-b_3)(b_2-b_4)(b_2-b_5)(b_2+c_5-a_2)}
{(b_2-b_1)(b_2+c_1-a_2)(b_2+c_2-a_2)(b_2+c_3-a_2)(b_2+c_4-a_2)} \\
 & \\
\<\v_3,\v_3\> = &
\frac{(b_3-b_4)(b_3-b_5)(b_3+c_4-a_2)(b_3+c_5-a_2)}
{(b_3-b_1)(b_3-b_2)(b_3+c_1-a_2)(b_3+c_2-a_2)(b_3+c_3-a_2)} \\
 & \\
\<\v_4,\v_4\> = &
\frac{(b_4-b_5)(b_4+c_3-a_2)(b_4+c_4-a_2)(b_4+c_5-a_2)}
{(b_4-b_1)(b_4-b_2)(b_4-b_3)(b_4+c_1-a_2)(b_4+c_2-a_2)} \\
 & \\
\<\v_5,\v_5\> = &
\frac{(b_5+c_2-a_2)(b_5+c_3-a_2)(b_5+c_4-a_2)(b_5+c_5-a_2)}
{(b_5-b_1)(b_5-b_2)(b_5-b_3)(b_5-b_4)(b_5+c_1-a_2)}
\end{array}
$$
\end{example}

 Let $S'((1,m-1),(1^m),(1^m))$ be obtained from $S((1,m-1),(1^m),(1^m))$
by removing the hyperplanes which are zero levels of the linear
forms in the denominators of $\<\v_i,\v_i\>$. Let
$S''((1,m-1),(1^m),(1^m))$ be obtained from
$S'((1,m-1),(1^m),(1^m))$ by removing the hyperplanes which are
zero levels of the linear forms in the numerators of
$\<\v_i,\v_i\>$. It is clear that if $\(\s(\A),\s(\B),\s(\C)\)$
lies in $S'((1,m-1),(1^m),(1^m))$, then the form
(\ref{equation:hg:form}) is well-defined. It is clear that if
$\(\s(\A),\s(\B),\s(\C)\)$ lies in $S''((1,m-1),(1^m),(1^m))$,
then the form (\ref{equation:hg:form}) is non-degenerate.

\begin{theorem}
\label{thm:hg:form} The operators $\A$, $\B$, and $\C$ are
self-adjoint with respect to the scalar product
(\ref{equation:hg:form}).
\end{theorem}

 Now suppose that the eigenvalues of the matrices $\A$, $\B$, and $\C$
are real numbers.
Let $b_1>b_2>\cdots >b_m$ and $c_1>c_2>\cdots >c_m$. We call a real
symmetric bilinear (or hermitian) form {\it sign-definite} if it
is either positive-definite or negative-definite.

\begin{theorem}
\label{thm:hg:sign_of_form} Let the form $\<*,*\>$ be defined by
(\ref{equation:hg:form}). Then it is sign-definite precisely in
the following two situations:

\begin{equation}
\label{equation:hg:form:pos}
\begin{array}{|c|c|}
\hline
 & \\
\begin{array}{lcccl}
b_m+c_1     & > & a_2    & > & b_m+c_2     \\
b_{m-1}+c_2 & > & a_2    & > & b_{m-1}+c_3 \\
\vdots      &   & \vdots &   & \vdots      \\
b_2+c_{m-1} & > & a_2    & > & b_2+c_m     \\
b_1+c_m     & > & a_2
\end{array} &

\begin{array}{lcccl}
b_1+c_{m-1} & > & a_2    & > & b_1+c_m     \\
b_2+c_{m-2} & > & a_2    & > & b_2+c_{m-1} \\
\vdots      &   & \vdots &   & \vdots      \\
b_{m-1}+c_1 & > & a_2    & > & b_{m-1}+c_2 \\
\phantom{b_{m-1}+c_1} &   & a_2 & > & b_m+c_1
\end{array} \\
 & \\
\hline
\end{array}.
\end{equation}
\medskip

\noindent If the form $\epsilon\, \<*,*\>$ is positive-definite
for $\epsilon=\pm 1$, then $\epsilon=sign(a_1-a_2)$. If the
inequalities of the first column hold, then $a_1>a_2$. If the
inequalities of the second column hold, then $a_1<a_2$.
\end{theorem}

\subsection{Even Family}  
\label{subsec:results:even}

 Let us pick a vector $(\underbrace{a_1\cdots ,a_1}_{m~times},
\underbrace{a_2,\cdots ,a_2}_{m~times},b_1,\underbrace{b_2,\cdots
,b_2}_{m-1~times},
\underbrace{b_3,\cdots ,b_3}_{m~times},c_1,\cdots ,c_{2m})$ from \\
$S((m,m),(1,m-1,m),(1^{2m}))$. Recall that this means $a_1\ne
a_2$, all $b_i$ are distinct, all $c_j$ are distinct, and the trace
condition holds: $ma_1+ma_2=b_1+(m-1)b_2+mb_3+ \sum_{i=1}^{2m}
c_i$. Let us set up the following notation:

\begin{equation}
\label{eq:p,q}
\begin{array}{ll}
p_i^{jk}=c_i+b_j-a_k, & q_{ij}=c_i+c_j+b_2+b_3-a_1-a_2.
\end{array}
\end{equation}
\bigskip

 We now define the matrices $\B$ and $\C$ by setting

\begin{center}
\pspicture[](0,0)(9,7) \label{picture:even:B}
\psset{unit=1cm,linewidth=1pt,linecolor=blue} \rput(0,-0.5){
\pspolygon(1,1)(7,1)(7,7)(1,7) \psline(2,1)(2,7) \psline(4,1)(4,7)
\psline(1,4)(7,4) \psline(1,6)(7,6) \rput(0,4.1){$\B=$}
\rput(7.8,4.1){,where} \rput(1.5,7.3){1} \rput(3.1,7.3){m-1}
\rput(5.4,7.2){m} \rput(0.7,2.7){m} \rput(0.6,5){m-1}
\rput(0.8,6.6){1} \rput(1.5,6.5){$b_1$}
\rput(3.1,6.5){$B_{1,1+j}$} \rput(5.4,6.5){$B_{1,m+j}$}
\rput(1.5,5){0} \rput(3.1,5){$b_2\,Id_{m-1}$}
\rput(5.4,5){$B_{1+i,m+j}$} \rput(1.5,2.7){0} \rput(3,2.7){0}
\rput(5.4,2.7){$b_3\, Id_m$} }
\endpspicture
\end{center}

\begin{equation}
\label{equation:even:B:II} B_{1,1+j}=(-1)^{m+1-j}\,
\displaystyle{\frac{\prod\limits_{k=j+1}^m q_{k,2m-j}}
{\prod\limits_{k=m+1}^{2m-1-j} (c_k-c_{2m-j})}}~~~~(1\le j\le
m-1);
\end{equation}

\begin{equation}
\label{equation:even:B:III} B_{1,m+j}=(-1)^{m-j}\,
p_{m+1-j}^{31}\displaystyle{\frac{\prod\limits_{k=m+j}^{2m-1}
q_{m+1-j,k}}{\prod\limits_{k=1}^{m-j} (c_k-c_{m+1-j})}}~~~~(1\le
j\le m);
\end{equation}

\begin{equation}
\label{equation:even:B:VI}
\begin{array}{c}
B_{1+i,m+j}=(-1)^{m-j}\, p_{m+1-j}^{31}\displaystyle{\frac{\prod
\limits_{k=1\atop k\ne m+1-j}^{i} q_{k,2m-i}
\prod\limits_{k=m+j\atop k\ne 2m-i}^{2m-1}
q_{m+1-j,k}}{\prod\limits_{k=2m+1-i}^{2m-1} (c_{2m-i}-c_k)
\prod\limits_{k=1}^{m-j} (c_k-c_{m+1-j})}}\\
 \\
(1\le i\le m-1,~~1\le j\le m);
\end{array}
\end{equation}
\bigskip

\begin{center}
\pspicture[](0,0)(9,7) \label{picture:even:C}
\psset{unit=1cm,linewidth=1pt,linecolor=blue} \rput(0,-0.5){
\pspolygon(0.8,1)(7,1)(7,7)(0.8,7) \psline(2,1)(2,7)
\psline(4,1)(4,7) \psline(0.8,4)(7,4) \psline(0.8,6)(7,6)
\rput(-0.2,4.1){$\C=$} \rput(7.8,4.1){,where} \rput(1.4,7.3){1}
\rput(3.1,7.3){m-1} \rput(5.4,7.2){m} \rput(0.5,2.7){m}
\rput(0.4,5){m-1} \rput(0.6,6.6){1} \rput(1.4,6.5){$c_{2m}$}
\rput(3.1,6.5){0} \rput(5.4,6.5){0} \rput(1.45,5){$C_{1+i,1}$}
\rput(2.6,5.8){$c_{2m-1}$} \rput(3.6,4.2){$c_{m+1}$}
\rput(2.7,5.3){.} \rput(3,5){.} \rput(3.3,4.7){.}
\rput(3.7,5.6){0} \rput(2.4,4.5){0} \rput(5.4,5){0}
\rput(4.3,3.8){$c_m$} \rput(6.8,1.2){$c_1$} \rput(5.1,2.9){.}
\rput(5.5,2.5){.} \rput(5.9,2.1){.} \rput(6.5,3.4){0}
\rput(4.6,1.7){0} \rput(1.45,2.7){$C_{m+i,1}$}
\rput(3,2.7){$C_{m+i,1+j}$}
}
\endpspicture
\end{center}

\begin{equation}
\label{equation:even:C:IV} C_{1+i,1}=-\frac{\prod\limits_{k=1}^i
q_{k,2m-i}} {\prod\limits_{k=2m+1-i}^{2m-1}
(c_{2m-i}-c_k)}~~~~(1\le i\le m-1);
\end{equation}

\begin{equation}
\label{equation:even:C:VII} C_{m+i,1}=-p_{m+1-i}^{32}\,
\frac{\prod\limits_{k=m+1}^{m-1+i} q_{m+1-i,k}}
{\prod\limits_{k=m+2-i}^m(c_{m+1-i}-c_k)}~~~~(1\le i\le m);
\end{equation}

\begin{equation}
\label{equation:even:C:VIII}
\begin{array}{c}
C_{m+i,1+j}=(-1)^{m+1-j}\, p_{m+1-i}^{32}\,
\frac{\prod\limits_{k=m+1\atop k\ne 2m-j}^{m-1+i} q_{m+1-i,k}
\prod\limits_{k=j+1\atop k\ne m+1-i}^m q_{k,2m-j}}
{\prod\limits_{k=m+2-i}^m (c_{m+1-i}-c_k) \prod\limits_{k=m+1}^{2m-1-j} (c_k-c_{2m-j})} \\
 \\
(1\le i\le m~~1\le j\le m-1).
\end{array}
\end{equation}

\noindent Here is an example with $m=3$.

\begin{example}
\label{example:even:m=3:B,C}

$$
\B=\left[
\begin{array}{c|cc|ccc}
b_1 & -\frac{q_{25}\, q_{35}}{c_4-c_5} & q_{34} & \frac{p^{31}_3\,
q_{34}\, q_{35}}{(c_1-c_3)(c_2-c_3)} &
-\frac{p^{31}_2\, q_{25}}{c_1-c_2} & p^{31}_1 \\
 & & & & & \\
\hline
 & & & & & \\
0 & b_2 & 0 & \frac{p^{31}_3\, q_{15}\, q_{34}}
{(c_1-c_3)(c_2-c_3)} & -\frac{p^{31}_2\, q_{15}}{c_1-c_2} & p^{31}_1 \\
 & & & & & \\
0 & 0 & b_2 & \frac{p^{31}_3\, q_{14}\, q_{24}\,
q_{35}}{(c_1-c_3)(c_2-c_3)(c_4-c_5)} & -\frac{p^{31}_2\, q_{14}\,
q_{25}}{(c_1-c_2)(c_4-c_5)} &
\frac{p^{31}_1\, q_{24}}{c_4-c_5} \\
 & & & & & \\
\hline
 & & & & & \\
0 & 0 & 0 & b_3 & 0 & 0 \\
 & & & & & \\
0 & 0 & 0 & 0 & b_3 & 0 \\
 & & & & & \\
0 & 0 & 0 & 0 & 0 & b_3
\end{array}
\right]
$$
\bigskip
$$
\C=\left[
\begin{array}{c|cc|ccc}
c_6 & 0  & 0 & 0 & 0 & 0 \\
 & & & & & \\
\hline
 & & & & & \\
-q_{15} & c_5 & 0 & 0 & 0 & 0 \\
 & & & & & \\
-\frac{q_{14}\, q_{24}}{c_4-c_5} & 0 & c_4 & 0 & 0 & 0 \\
 & & & & & \\
\hline
 & & & & & \\
-p^{32}_3 & -\frac{p^{32}_3\, q_{25}}{c_4-c_5} & p^{32}_3 & c_3 & 0 & 0 \\
 & & & & & \\
-\frac{p^{32}_2\, q_{24}}{c_2-c_3} & -\frac{p^{32}_2\, q_{24}\,
q_{35}}{(c_2-c_3)(c_4-c_5)} &
\frac{p^{32}_2\, q_{34}}{c_2-c_3} & 0 & c_2 & 0 \\
 & & & & & \\
-\frac{p^{32}_1\, q_{14}\, q_{15}}{(c_1-c_2)(c_1-c_3)} &
-\frac{p^{32}_1\, q_{14}\, q_{25}\, q_{35}}
{(c_1-c_2)(c_1-c_3)(c_4-c_5)} & \frac{p^{32}_1\, q_{15}\,
q_{34}}{(c_1-c_2)(c_1-c_3)} & 0 & 0 & c_1
\end{array}
\right]
$$
\end{example}
\bigskip

\noindent It is clear that $\B$ and $\C$ are diagonalizable and
that
their spectra are \\
$\s(\B)=\{b_1,\underbrace{b_2,b_2,\cdots ,b_2,}_{m-1~times},
\underbrace{b_3,b_3,\cdots ,b_3,b_3}_{m~times}\}$,
$\s(\C)=\{c_1,\cdots ,c_{2m}\}$.

\begin{theorem}
\label{thm:even:A}
If $\B$ and $\C$ are as above, then $\A=\B+\C$ is diagonalizable and \\
$\s(\A)=\{\underbrace{a_1,\cdots ,a_1,}_{m~times}
\underbrace{a_2,\cdots ,a_2}_{m~times}\}$.
\end{theorem}

 For $i=1,\cdots ,2m$, let $\v_i=(0,\cdots ,0,1,v^{2m+2-i}_i,\cdots ,v^{2m}_i)$
be the eigenvector of the matrix $\C$ with the eigenvalue $c_i$.

\begin{lemma}
\label{lemma:even:eigenbasis_of_C}
\begin{enumerate}

\item
For $1\le i\le m$, we have $v_i^j=0$ for all $2m+1-i<j$.

\item
For $1\le i\le m-1$, we have $v_{m+i}^j=0$ when $m+1-i<j\le m$ and \\
$v_{m+i}^{m+j}=(-1)^i\, p^{32}_{m+1-j}\,
\displaystyle{\frac{\prod\limits_ {k=m+1\atop k\ne m+i}^{m-1+j}
q_{m+1-j,k}\,\prod\limits_{k=m+1-i\atop k\ne m+1-j}^m
q_{k,m+i}}{(c_{m+1-j}-c_{m+i})\prod\limits_{k=m+2-j}^m
(c_{m+1-j}-c_k)\,\prod\limits
_{k=m+1}^{m-1+i} (c_k-c_{m+i})}}$ \\
for $1\le j\le m$.

\item
For $1\le j\le m-1$ we have $v_{2m}^{1+j}=
\displaystyle{\frac{\prod\limits_{k=1}^j q_{k,2m-j}}
{\prod\limits_{k=2m+1-j}^{2m} (c_{2m-j}-c_k)}}$ and for $1\le j\le
m$ we have

\begin{equation}
\label{equation:even:v_{2m}^{m+j}} v_{2m}^{m+j} =(-1)^{m+1}\,
\displaystyle{\frac{p^{32}_{m+1-j}}{c_{m+1-j}-c_{2m}}}
\displaystyle{\frac{\prod\limits_{k=1\atop k\ne m+1-j}^m q_{k,2m}
\prod\limits_{k=m+1}^{m-1+j} q_{m+1-j,k}}
{\prod\limits_{k=m+1}^{2m-1}(c_k-c_{2m})
\prod\limits_{k=m+2-j}^m(c_{m+1-j}-c_k)}}.
\end{equation}

\end{enumerate}
\end{lemma}

\noindent Here is an example with $m=3$ ($\e_1,\cdots ,\e_{2m}$ is
the standard basis of $V$).

\begin{example}
\label{example:even:m=3:evecC}
$$
\begin{array}{l}
\v_1 =  \e_6,~~ \v_2 = \e_5,~~ \v_3 = \e_4, \\
 \\
\v_4 = \( 0,0,1,-\frac{p^{32}_3}{c_3-c_4},-\frac{p^{32}_2\,
q_{34}}{(c_2-c_3)(c_2-c_4)},
-\frac{p^{32}_1\, q_{15}\, q_{34}}{(c_1-c_2)(c_1-c_3)(c_1-c_4)} \), \\
 \\
\v_5 = \( 0,1,0,\frac{p^{32}_3\,
q_{25}}{(c_3-c_5)(c_4-c_5)},\frac{p^{32}_2\, q_{24}\, q_{35}}
{(c_2-c_3)(c_2-c_5)(c_4-c_5)},\frac{p^{32}_1\, q_{14}\, q_{25}\,
q_{35}}{(c_1-c_2)(c_1-c_3)
(c_1-c_5)(c_4-c_5)} \), \\
 \\
\v_6 = \( 1,\frac{q_{15}}{c_5-c_6},\frac{q_{14}\,
q_{24}}{(c_4-c_5)(c_4-c_6)}, \frac{p^{32}_3\, q_{16}\,
q_{26}}{(c_3-c_6)(c_4-c_6)(c_5-c_6)}, \frac{p^{32}_2\, q_{24}\,
q_{16}\, q_{36}}
{(c_2-c_3)(c_2-c_6)(c_4-c_6)(c_5-c_6)}, \right.\\
 \\
\phantom{\v_6 = \( \right.} \left. \frac{p^{32}_1\, q_{14}\,
q_{15}\, q_{26}\, q_{36}}
           {(c_1-c_2)(c_1-c_3)(c_1-c_6)(c_4-c_6)(c_5-c_6)} \right).
\end{array}
$$
\end{example}

 We define a scalar product on $V$ by setting

\begin{equation}
\label{equation:even:form} \<\v_i,\v_j\>=\delta_{ij}\,
\displaystyle{\frac{\prod\limits_{k=i+1}^{2m}(c_i-c_k)}
{\prod\limits_{k=1}^{i-1}(c_i-c_k)}}\,\displaystyle{\frac{\prod\limits_
{k=2m+1-i\atop k\ne i}^{2m} q_{ik}}{\prod\limits_{k=1\atop k\ne
i}^{2m-i}q_{ik}}} \,\times\,\left\{
\begin{array}{lll}
\displaystyle{\frac{p^{31}_i}{p^{32}_i}} & \mbox{,if} & i\le m; \\
                                         &            &         \\
p^{31}_i\, p^{32}_i                      & \mbox{,if} & i> m.   \\
\end{array}
\right.
\end{equation}

\noindent Here is an example with $m=3$.

\begin{example}
\label{example:even:m=3:norm_of_v_i}
$$
\begin{array}{lccc}
\<\v_1,\v_1\>= &
(c_1-c_2)(c_1-c_3)(c_1-c_4)(c_1-c_5)(c_1-c_6)\times &
\frac{p^{31}_1}{p^{32}_1}\times & \frac{q_{16}}{q_{12}\, q_{13}\, q_{14}\, q_{15}}, \\
 & & & \\
\<\v_2,\v_2\>= &
\frac{(c_2-c_3)(c_2-c_4)(c_2-c_5)(c_2-c_6)}{c_2-c_1}\times &
\frac{p^{31}_2}{p^{32}_2}\times & \frac{q_{25}\, q_{26}}{q_{12}\, q_{23}\, q_{24}}, \\
 & & & \\
\<\v_3,\v_3\>= &
\frac{(c_3-c_4)(c_3-c_5)(c_3-c_6)}{(c_3-c_1)(c_3-c_2)}\times &
\frac{p^{31}_3}{p^{32}_3}\times & \frac{q_{34}\, q_{35}\, q_{36}}{q_{13}\, q_{23}}, \\
 & & & \\
\<\v_4,\v_4\>= &
\frac{(c_4-c_5)(c_4-c_6)}{(c_4-c_1)(c_4-c_2)(c_4-c_3)}\times &
p^{31}_4\, p^{32}_4\times & \frac{q_{34}\, q_{45}\, q_{46}}{q_{14}\, q_{24}}, \\
 & & & \\
\<\v_5,\v_5\>= &
\frac{c_5-c_6}{(c_5-c_1)(c_5-c_2)(c_5-c_3)(c_5-c_4)}\times &
p^{31}_5\, p^{32}_5\times & \frac{q_{25}\, q_{35}\, q_{45}\, q_{56}}{q_{15}}, \\
 & & & \\
\<\v_6,\v_6\> = &
\frac{1}{(c_6-c_1)(c_6-c_2)(c_6-c_3)(c_6-c_4)(c_6-c_5)}\times &
p^{31}_6\, p^{32}_6\times & q_{16}\, q_{26}\, q_{36}\, q_{46}\,
q_{56}.
\end{array}
$$
\end{example}

 The sets $S'((m,m),(1,m-1,m),(1^{2m}))$ and $S''((m,m),(1,m-1,m),(1^{2m}))$ are constructed
from (\ref{equation:even:form}) similarly to the hypergeometric
case (see page \pageref{equation:hg:form}
) and have the same properties.

\begin{theorem}
\label{thm:even:form} The operators $\A$, $\B$, and $\C$ are
self-adjoint with respect to the scalar product
(\ref{equation:even:form}).
\end{theorem}

 Now suppose that the eigenvalues of the matrices $\A$, $\B$, and $\C$
are real numbers.
Let $a_1>a_2$ and $c_1>c_2>\cdots >c_{2m}$.

\begin{theorem}
\label{thm:even:sign_of_form} The form $\<*,*\>$ defined by
(\ref{equation:even:form}) is sign-definite precisely in the
following six situations:

$$
\begin{array}{|l|l|l|}
\hline
 & & \\
\begin{array}{l}
b_1>b_3>b_2                  \\
                             \\
p^{31}_{m-1}>0>p^{31}_m      \\
p^{32}_{2m-1}>0>p^{32}_{2m}  \\
                             \\
q_{1,2m-2}>0>q_{1,2m-1}      \\
q_{2,2m-3}>0>q_{2,2m-2}      \\
q_{3,2m-4}>0>q_{3,2m-3}      \\
\vdots                       \\
q_{m-1,m}>0>q_{m-1,m+1}
\end{array}
 &
\begin{array}{l}
b_1>b_2>b_3                       \\
                                  \\
\phantom{p^{32}_m>}0>p^{31}_1     \\
p^{32}_m>0>p^{32}_{m+1}           \\
                                  \\
q_{1,2m-1}>0>q_{1,2m}             \\
q_{2,2m-2}>0>q_{2,2m-1}           \\
\vdots                            \\
q_{m-1,m+1}>0>q_{m-1,m+2}         \\
\phantom{q_{m-1,m+1}>}0>q_{m,m+1}
\end{array}
 &
\begin{array}{l}
b_2>b_1>b_3                   \\
                              \\
\phantom{p^{32}_m>}0>p^{31}_1 \\
p^{32}_m>0>p^{32}_{m+1}       \\
                              \\
q_{1,2m}>0>q_{2,2m}           \\
q_{2,2m-1}>0>q_{3,2m-1}       \\
\vdots                        \\
q_{m-1,m+2}>0>q_{m,m+2}       \\
q_{m,m+1}>0
\end{array}
 \\
 & & \\
\hline
 & & \\
\begin{array}{l}
b_2>b_3>b_1                 \\
                            \\
p^{31}_1>0>p^{31}_2         \\
p^{32}_{m+1}>0>p^{32}_{m+2} \\
                            \\
q_{2,2m}>0>q_{3,2m}         \\
q_{3,2m-1}>0>q_{4,2m-1}     \\
q_{4,2m-2}>0>q_{5,2m-2}     \\
\vdots                      \\
q_{m,m+2}>0>q_{m+1,m+2}
\end{array}
 &
\begin{array}{l}
b_3>b_2>b_1             \\
                        \\
p^{31}_m>0>p^{31}_{m+1} \\
p^{32}_{2m}>0           \\
                        \\
q_{1,2m}>0>q_{2,2m}     \\
q_{2,2m-1}>0>q_{3,2m-1} \\
\vdots                  \\
q_{m-1,m+2}>0>q_{m,m+2} \\
q_{m,m+1}>0
\end{array}
 &
\begin{array}{l}
b_3>b_1>b_2                      \\
                                 \\
p^{31}_m>0>p^{31}_{m+1}          \\
p^{32}_{2m}>0                    \\
                                 \\
q_{1,2m-1}>0>q_{1,2m}            \\
q_{2,2m-2}>0>q_{2,2m-1}          \\
\vdots                           \\
q_{m-1,m+1}>0>q_{m-1,m+2}        \\
\phantom{q_{m-1,m+1}>}0>q_{m,m+1}
\end{array}
 \\
 & & \\
\hline
\end{array}.
$$
\medskip

\noindent If the form $\epsilon\, \<*,*\>$ is positive-definite
for $\epsilon=\pm 1$, then $\epsilon=sign((b_1-b_2)(b_1-b_3))$. In
each case, the inequalities between $b_1$, $b_2$, and $b_3$ are
implied by other inequalities.
\end{theorem}

\subsection{Odd Family}  
\label{subsection:results:odd}

 Let us pick a vector $(\underbrace{a_1,\cdots ,a_1}_{m+1~times},
\underbrace{a_2,\cdots ,a_2}_{m~times},b_1,\underbrace{b_2,\cdots
,b_2}_{m~times},
\underbrace{b_3,\cdots ,b_3}_{m~times},c_1,\cdots ,c_{2m+1})$ from \\
$S((m+1,m),(1,m,m),(1^{2m+1}))$. Recall that this means $a_1\ne
a_2$, all $b_i$ are distinct, all $c_j$ are distinct, and the trace
condition holds: $(m+1)a_1+ma_2=b_1+mb_2+mb_3+
\sum_{i=1}^{2m+1} c_i$. \\

 We now define the matrices $\B$ and $\C$ by setting

\begin{center}
\pspicture[](0,0)(10,8) \label{picture:odd:B}
\psset{unit=1cm,linewidth=1pt,linecolor=blue} \rput(0,-0.5){
\pspolygon(1,1)(8,1)(8,8)(1,8) \psline(2,1)(2,8) \psline(5,1)(5,8)
\psline(1,4)(8,4) \psline(1,7)(8,7) \rput(1.5,8.3){1}
\rput(0,4.1){$\B=$} \rput(8.8,4.1){,where} \rput(3.4,8.3){m}
\rput(6.4,8.3){m} \rput(0.7,2.5){m} \rput(0.7,5.5){m}
\rput(0.8,7.6){1} \rput(1.5,7.5){$b_1$}
\rput(3.5,7.5){$B_{1,1+j}$} \rput(6.4,7.5){$B_{1,m+1+j}$}
\rput(1.5,5.5){0} \rput(3.5,5.5){$b_2\, Id_m$}
\rput(6.4,5.5){$B_{1+i,m+1+j}$} \rput(1.5,2.5){0}
\rput(3.5,2.5){0} \rput(6.4,2.5){$b_3\, Id_m$} }
\endpspicture
\end{center}

\begin{equation}
\label{equation:odd:B:II} B_{1,1+j}=(-1)^{m-j}\, p^{21}_{2m+1-j}\,
\displaystyle{\frac{\prod\limits_{k=j+1}^m q_{k,2m+1-j}}
{\prod\limits_{k=m+1}^{2m-j} (c_k-c_{2m+1-j})}}~~~~(1\le j\le m);
\end{equation}

\begin{equation}
\label{equation:odd:B:III} B_{1,m+1+j}=(-1)^{m-j}\,
p^{31}_{m+1-j}\, \displaystyle{\frac{\prod\limits_{k=m+1+j}^{2m}
q_{m+1-j,k}} {\prod\limits_{k=1}^{m-j} (c_k-c_{m+1-j})}}~~~~(1\le
j\le m);
\end{equation}

\begin{equation}
\label{equation:odd:B:VI}
\begin{array}{c}
B_{1+i,m+1+j}=(-1)^{m-j}\, p^{31}_{m+1-j}\,
\displaystyle{\frac{\prod\limits_{k=1\atop k\ne m+1-j}^i
q_{k,2m+1-i} \prod\limits_{k=m+1+j\atop k\ne 2m+1-i}^{2m}
q_{m+1-j,k}} {\prod\limits_{k=2m+2-i}^{2m} (c_{2m+1-i}-c_k)
\prod\limits_{k=1}^{m-j} (c_k-c_{m+1-j})}} \\
 \\
(1\le i\le m,~~~~1\le j\le m);
\end{array}
\end{equation}
\bigskip

\begin{center}
\pspicture[](0,0)(10,8) \label{picture:odd:C}
\psset{unit=1cm,linewidth=1pt,linecolor=blue} \rput(0,-0.5){
\pspolygon(0.5,1)(8,1)(8,8)(0.5,8) \psline(2,1)(2,8)
\psline(5,1)(5,8) \psline(0.5,4)(8,4) \psline(0.5,7)(8,7)
\rput(1.2,8.3){1} \rput(-0.3,4.1){$\C=$} \rput(8.8,4.1){,where}
\rput(3.4,8.3){m} \rput(6.4,8.3){m} \rput(0.2,2.5){m}
\rput(0.2,5.5){m} \rput(0.3,7.6){1} \rput(1.2,7.5){$c_{2m+1}$}
\rput(3.5,7.5){0} \rput(6.5,7.5){0} \rput(1.2,5.5){$C_{1+i,1}$}
\rput(2.4,6.8){$c_{2m}$} \rput(4.55,4.25){$c_{m+1}$}
\rput(3.1,5.9){.} \rput(3.5,5.5){.} \rput(3.9,5.1){.}
\rput(4.5,6.4){0} \rput(2.6,4.7){0} \rput(6.5,5.5){0}
\rput(1.27,2.5){$C_{m+1+i,1}$} \rput(3.5,2.5){$C_{m+1+i,1+j}$}
\rput(5.3,3.8){$c_m$} \rput(7.8,1.25){$c_1$} \rput(6.1,2.9){.}
\rput(6.5,2.5){.} \rput(6.9,2.1){.} \rput(7.5,3.4){0}
\rput(5.6,1.7){0} }
\endpspicture
\end{center}

\begin{equation}
\label{equation:odd:C:IV} C_{1+i,1}=-\frac{\prod\limits_{k=1}^i
q_{k,2m+1-i}} {\prod\limits_{k=2m+2-i}^{2m}
(c_{2m+1-i}-c_k)}~~~~(1\le i\le m);
\end{equation}

\begin{equation}
\label{equation:odd:C:VII}
C_{m+1+i,1}=-\frac{\prod\limits_{k=m+1}^{m+i} q_{m+1-i,k}}
{\prod\limits_{k=m+2-i}^m (c_{m+1-i}-c_k)}~~~~(1\le i\le m);
\end{equation}

\begin{equation}
\label{equation:odd:C:VIII}
\begin{array}{c}
C_{m+1+i,1+j}=(-1)^{m-j}\, p^{21}_{2m+1-j}\,
\frac{\prod\limits_{k=m+1\atop k\ne 2m+1-j}^{m+i} q_{m+1-i,k}
\prod\limits_{k=j+1\atop k\ne m+1-i}^m q_{k,2m+1-j}}
{\prod\limits_{k=m+2-i}^m (c_{m+1-i}-c_k) \prod\limits_{k=m+1}^{2m-j} (c_k-c_{2m+1-j})} \\
 \\
(1\le i\le m,~~~~1\le j\le m).
\end{array}
\end{equation}

\noindent Here is an example with $m=3$.

\begin{example}
\label{example:odd:m=3:B,C}

$$
\B=\left[
\begin{array}{c|ccc|ccc}
b_1 & \frac{p^{21}_6\, q_{26}\, q_{36}}{(c_4-c_6)(c_5-c_6)} &
-\frac{p^{21}_5\, q_{35}}{c_4-c_5} & p^{21}_4 & \frac{p^{31}_3\,
q_{35}\, q_{36}}{(c_1-c_3)(c_2-c_3)} & -\frac{p^{31}_2\,
q_{26}}{c_1-c_2} &
p^{31}_1  \\
 & & & & & & \\
\hline
 & & & & & & \\
0 & b_2 & 0 & 0 & \frac{p^{31}_3\, q_{16}\,
q_{35}}{(c_1-c_3)(c_2-c_3)} &
-\frac{p^{31}_2\, q_{16}}{c_1-c_2} & p^{31}_1 \\
 & & & & & & \\
0 & 0 & b_2 & 0 & \frac{p^{31}_3\, q_{15}\, q_{25}\,
q_{36}}{(c_1-c_3)(c_2-c_3)(c_5-c_6)} &
-\frac{p^{31}_2\, q_{15}\, q_{26}}{(c_1-c_2)(c_5-c_6)} & \frac{p^{31}_1\, q_{25}}{c_5-c_6} \\
 & & & & & & \\
0 & 0 & 0 & b_2 & \frac{p^{31}_3\, q_{14}\, q_{24}\, q_{35}\,
q_{36}}{(c_1-c_3)(c_2-c_3) (c_4-c_5)(c_4-c_6)} & -\frac{p^{31}_2\,
q_{14}\, q_{34}\, q_{26}}{(c_1-c_2)(c_4-c_5)(c_4-c_6)} &
\frac{p^{31}_1\, q_{24}\, q_{34}}{(c_4-c_5)(c_4-c_6)} \\
 & & & & & & \\
\hline
 & & & & & & \\
0 & 0 & 0 & 0 & b_3 & 0 & 0 \\
 & & & & & & \\
0 & 0 & 0 & 0 & 0 & b_3 & 0 \\
 & & & & & & \\
0 & 0 & 0 & 0 & 0 & 0 & b_3
\end{array}
\right]
$$
\bigskip
$$
\C=\left[
\begin{array}{c|ccc|ccc}
c_7 & 0  & 0 & 0 & 0 & 0 & 0 \\
 & & & & & & \\
\hline
 & & & & & & \\
-q_{16} & c_6 & 0 & 0 & 0 & 0 & 0 \\
 & & & & & & \\
-\frac{q_{15}\, q_{25}}{c_5-c_6} & 0 & c_5 & 0 & 0 & 0 & 0 \\
 & & & & & & \\
-\frac{q_{14}\, q_{24}\, q_{34}}{(c_4-c_5)(c_4-c_6)} & 0 & 0 & c_4 & 0 & 0 & 0 \\
 & & & & & & \\
\hline
  & & & & & & \\
-q_{34} & \frac{p^{21}_6\, q_{34}\, q_{26}}{(c_4-c_6)(c_5-c_6)} &
-\frac{p^{21}_5\, q_{34}}{c_4-c_5} & p^{21}_4 & c_3 & 0 & 0 \\
  & & & & & & \\
-\frac{q_{24}\, q_{25}}{c_2-c_3} & \frac{p^{21}_6\, q_{24}\,
q_{25}\, q_{36}} {(c_2-c_3)(c_4-c_6)(c_5-c_6)} & -\frac{p^{21}_5\,
q_{24}\, q_{35}}{(c_2-c_3)(c_4-c_5)} &
\frac{p^{21}_4\, q_{25}}{c_2-c_3} & 0 & c_2 & 0 \\
  & & & & & & \\
-\frac{q_{14}\, q_{15}\, q_{16}}{(c_1-c_2)(c_1-c_3)} &
\frac{p^{21}_6\, q_{14}\, q_{15}\, q_{26}\,
q_{36}}{(c_1-c_2)(c_1-c_3)(c_4-c_6)(c_5-c_6)} & -\frac{p^{21}_5\,
q_{14}\, q_{16}\, q_{35}} {(c_1-c_2)(c_1-c_3)(c_4-c_5)} &
\frac{p^{21}_4\, q_{15}\, q_{16}}{(c_1-c_2)(c_1-c_3)} & 0 & 0 &
c_1
\end{array}
\right].
$$
\end{example}
\bigskip

\noindent It is clear that $\B$ and $\C$ are diagonalizable and
that their spectra are $\s(\B)=\{b_1,\underbrace{b_2,b_2,\cdots
,b_2,}_{m~\mbox{{\it times}}} \underbrace{b_3,b_3,\cdots
,b_3}_{m~\mbox{{\it times}}}\}$, $\s(\C)=\{c_1,c_2,\cdots
,c_{2m+1}\}$.

\begin{theorem}
\label{thm:odd:A} Let $\B$ and $\C$ be as above and let
$\A=\B+\C$. Then $\A$ is diagonalizable and
$\s(\A)=\{\underbrace{a_1,\cdots ,a_1}_{m+1~\mbox{times}}
\underbrace{a_2,\cdots ,a_2}_{m~\mbox{times}}\}$.
\end{theorem}

 For $i=1,\cdots ,2m+1$, let $\v_i=(0,\cdots ,0,1,v^{2m+3-i}_i,\cdots ,v^{2m+1}_i)$
be the eigenvector of the matrix $\C$ with the eigenvalue $c_i$.

\begin{lemma}
\label{lemma:odd:eigenbasis_of_C}
\begin{enumerate}
\item
For $1\le i\le m$, we have $v_i^j=0$ for all $2m+2-i<j$.

\item
For $1\le i\le m$, we have $v_{m+i}^j=0$ when $2m+2-i<j\le m+1$ and \\
$v_{m+i}^{m+1+j}= \displaystyle{\frac{(-1)^i\,
p^{21}_{m+i}}{c_{m+1-j}-c_{m+i}}}
\displaystyle{\frac{\prod\limits_{k=m+1\atop k\ne m+i}^{m+j}
q_{m+1-j,k}\, \prod\limits_{k=m+2-i\atop k\ne m+1-j}^m q_{k,m+i}}
{\prod\limits_{k=m+2-j}^m (c_{m+1-j}-c_k)\,
\prod\limits_{k=m+1}^{m-1+i} (c_k-c_{m+i})}}$ for $1\le j\le m$.

\item
For $1\le i\le m$, we have $v_{2m+1}^{1+i}=
\displaystyle{\frac{\prod\limits_{k=1}^i q_{k,2m+1-i}}
{\prod\limits_{k=2m+2-i}^{2m+1} (c_{2m+1-i}-c_k)}}$ and for $1\le
j\le m$, we have

\begin{equation}
\label{equation:odd:v_{2m+1}^{m+1+j}} v_{2m+1}^{m+1+j}=(-1)^m\,
\displaystyle{\frac{p^{21}_{2m+1}}{c_{m+1-j}-c_{2m+1}}}
\displaystyle{\frac{\prod\limits_{k=1\atop k\ne m+1-j}^m
q_{k,2m+1} \prod\limits_{k=m+1}^{m+j} q_{m+1-j,k}}
{\prod\limits_{k=m+1}^{2m}(c_k-c_{2m+1}) \prod\limits_{k=m+2-j}^m
(c_{m+1-j}-c_k)}}.
\end{equation}

\end{enumerate}
\end{lemma}

\noindent Here is an example with $m=3$.

\begin{example}
\label{example:odd:m=2:evecC}
$$
\begin{array}{l}
\v_1 =  \e_7,~~ \v_2 = \e_6,~~ \v_3 = \e_5, \\
 \\
\v_4 = \( 0,0,0,1,-\frac{p^{21}_4}{c_3-c_4},-\frac{p^{21}_4\,
q_{25}}{(c_2-c_3)(c_2-c_4)},
          -\frac{p^{21}_4\, q_{15}\, q_{16}}{(c_1-c_2)(c_1-c_3)(c_1-c_4)} \), \\
 \\
\v_5 = \( 0,0,1,0,\frac{p^{21}_5\,
q_{34}}{(c_3-c_5)(c_4-c_5)},\frac{p^{21}_5\, q_{24}\, q_{35}}
          {(c_2-c_3)(c_2-c_5)(c_4-c_5)}, \frac{p^{21}_5\, q_{14}\, q_{16}\, q_{35}}
          {(c_1-c_2)(c_1-c_3)(c_1-c_5)(c_4-c_5)} \), \\
 \\
\v_6 = \( 0,1,0,0,-\frac{p^{21}_6\, q_{26}\,
q_{34}}{(c_3-c_6)(c_4-c_6)(c_5-c_6)},
          -\frac{p^{21}_6\, q_{24}\, q_{25}\, q_{36}}{(c_2-c_3)(c_2-c_6)(c_4-c_6)(c_5-c_6)},
          -\frac{p^{21}_6\, q_{14}\, q_{15}\, q_{26}\, q_{36}}
          {(c_1-c_2)(c_1-c_3)(c_1-c_6)(c_4-c_6)(c_5-c_6)} \), \\
 \\
\v_7 = \(1, \frac{q_{16}}{c_6-c_7}, \frac{q_{15}\,
q_{25}}{(c_5-c_6)(c_5-c_7)},
            \frac{q_{14}\, q_{24}\, q_{34}}{(c_4-c_5)(c_4-c_6)(c_4-c_7)},
            -\frac{p^{21}_7\, q_{17}\, q_{27}\, q_{34}}{(c_3-c_7)(c_4-c_7)(c_5-c_7)(c_6-c_7)},
       \right. \\
 \\
\phantom{\v_7 = \( \right.} \left. -\frac{p^{21}_7\, q_{17}\,
q_{37}\, q_{24}\, q_{25}}
            {(c_2-c_3)(c_2-c_7)(c_4-c_7)(c_5-c_7)(c_6-c_7)},
            -\frac{p^{21}_7\, q_{27}\, q_{37}\, q_{14}\, q_{15}\, q_{16}}
            {(c_1-c_2)(c_1-c_3)(c_1-c_7)(c_4-c_7)(c_5-c_7)(c_6-c_7)} \).
\end{array}
$$
\end{example}

 We define a scalar product on $V$ by setting

\begin{equation}
\label{equation:odd:form} \<\v_i,\v_j\>=\delta_{ij}\,
\displaystyle{\frac{\prod\limits_{k=1+i}^{2m+1}(c_i-c_k)}
{\prod\limits_{k=1}^{i-1}(c_i-c_k)}}\,\displaystyle{\frac{\prod\limits_
{k=2m+2-i\atop k\ne i}^{2m+1} q_{ik}}{\prod\limits_{k=1\atop k\ne
i}^{2m+1-i}q_{ik}}} \,\times\,\left\{
\begin{array}{lll}
\displaystyle{\frac{p^{31}_i}{p^{21}_i}} & \mbox{,if} & i\le m; \\
                                   &            &         \\
p^{31}_i\, p^{21}_i                      & \mbox{,if} & i> m.
\end{array}
\right.
\end{equation}

\noindent Here is an example with $m=3$.

\begin{example}
\label{example:odd:m=3:norm_of_v_i}
$$
\begin{array}{lccc}
\<\v_1,\v_1\>= &
(c_1-c_2)(c_1-c_3)(c_1-c_4)(c_1-c_5)(c_1-c_6)(c_1-c_7)\times &
\frac{p^{31}_1}{p^{21}_1}\times &
\frac{q_{17}}{q_{12}\, q_{13}\, q_{14}\, q_{15}\, q_{16}} \\
 & & & \\
\<\v_2,\v_2\>= &
\frac{(c_2-c_3)(c_2-c_4)(c_2-c_5)(c_2-c_6)(c_2-c_7)}{c_2-c_1}\times
& \frac{p^{31}_2}{p^{21}_2}\times &
\frac{q_{26}\, q_{27}}{q_{12}\, q_{23}\, q_{24}\, q_{25}} \\
 & & & \\
\<\v_3,\v_3\>= &
\frac{(c_3-c_4)(c_3-c_5)(c_3-c_6)(c_3-c_7)}{(c_3-c_1)(c_3-c_2)}\times
& \frac{p^{31}_3}{p^{21}_3}\times &
\frac{q_{35}\, q_{36}\, q_{37}}{q_{13}\, q_{23}\, q_{34}} \\
 & & & \\
\<\v_4,\v_4\>= &
\frac{(c_4-c_5)(c_4-c_6)(c_4-c_7)}{(c_4-c_1)(c_4-c_2)(c_4-c_3)}\times
&
p^{31}_4\, p^{21}_4 \times & \frac{q_{45}\, q_{46}\, q_{47}}{q_{14}\, q_{24}\, q_{34}} \\
 & & & \\
\<\v_5,\v_5\>= &
\frac{(c_5-c_6)(c_5-c_7)}{(c_5-c_1)(c_5-c_2)(c_5-c_3)(c_5-c_4)}\times
&
p^{31}_5\, p^{21}_5 \times & \frac{q_{35}\, q_{45}\, q_{56}\, q_{57}}{q_{15}\, q_{25}} \\
 & & & \\
\<\v_6,\v_6\>= &
\frac{c_6-c_7}{(c_6-c_1)(c_6-c_2)(c_6-c_3)(c_6-c_4)(c_6-c_5)}\times
&
p^{31}_6\, p^{21}_6 \times & \frac{q_{26}\, q_{36}\, q_{46}\, q_{56}\, q_{67}}{q_{16}} \\
 & & & \\
\<\v_7,\v_7\>= &
\frac{1}{(c_7-c_1)(c_7-c_2)(c_7-c_3)(c_7-c_4)(c_7-c_5)(c_7-c_6)}\times
& p^{31}_7\, p^{21}_7 \times & q_{17}\, q_{27}\, q_{37}\, q_{47}\,
q_{57}\, q_{67}
\end{array}
$$
\end{example}

 The sets $S'((m+1,m),(1,m,m),(1^{2m+1}))$ and $S''((m+1,m),(1,m,m),(1^{2m+1}))$ are constructed
from (\ref{equation:odd:form}) similarly to the hypergeometric
case (see page \pageref{equation:hg:form}
) and have the same properties.

\begin{theorem}
\label{thm:odd:form} The operators $\A$, $\B$, and $\C$ are
self-adjoint with respect to the scalar product
(\ref{equation:odd:form}).
\end{theorem}

 Now suppose that the eigenvalues of the matrices $\A$, $\B$, and $\C$
are real numbers.
Let $b_2>b_3$ and $c_1>c_2>\cdots >c_{2m+1}$.

\begin{theorem}
\label{thm:odd:sign_of_form} Under the condition $b_2>b_3$, the
form $\<*,*\>$ defined by (\ref{equation:odd:form}) is
sign-definite precisely in the following three situations:

$$
\begin{array}{|c|c|c|}
\hline
 & & \\
\begin{array}{l}
b_1>b_2>b_3                   \\
                              \\
\phantom{p^{21}_m>}0>p^{31}_1 \\
p^{21}_m>0>p^{21}_{m+1}       \\
                              \\
q_{1,2m}>0>q_{1,2m+1}         \\
q_{2,2m-1}>0>q_{2,2m}         \\
q_{3,2m-2}>0>q_{3,2m-1}       \\
\vdots                        \\
q_{m,m+1}>0>q_{m,m+2}
\end{array} &

\begin{array}{l}
b_2>b_1>b_3                   \\
                              \\
\phantom{p^{31}_m>}0>p^{31}_1 \\
p^{21}_{m+1}>0>p^{21}_{m+2}   \\
                              \\
q_{1,2m+1}>0>q_{2,2m+1}       \\
q_{2,2m}>0>q_{3,2m}           \\
\vdots                        \\
q_{m-1,m+3}>0>q_{m,m+3}       \\
q_{m,m+2}>0>q_{m+1,m+2}       \\
\end{array} &

\begin{array}{l}
b_2>b_3>b_1                 \\
                            \\
p^{31}_1>0>p^{31}_2         \\
r^{21}_{m+1}>0>r^{21}_{m+2} \\
                            \\
q_{2,2m+1}>0>q_{3,2m+1}     \\
q_{3,2m}>0>q_{4,2m}         \\
\vdots                      \\
q_{m,m+3}>0>q_{m+1,m+3}     \\
q_{m+1,m+2}>0
\end{array} \\
 & & \\
\hline
\end{array}
$$
\medskip

\noindent If the form $\epsilon\, \<*,*\>$ is positive-definite
for $\epsilon=\pm 1$, then $\epsilon=sign((b_1-b_2)(b_1-b_3))$. In
each case, the inequalities between $b_1$ and $b_2$ or $b_3$ are
implied by other inequalities.

\end{theorem}

\subsection{Extra Case}  
\label{subsection:results:extra}

 Let us pick a vector $(a_1,a_1,a_1,a_1,a_2,a_2,b_1,b_1,b_2,b_2,b_3,b_3,c_1,c_2,c_3,c_4,c_5,c_6)$
from $S((4,2),(2,2,2),(1^6))$. Recall that this means $a_1\ne
a_2$, all $b_i$ are distinct, all $c_j$ are distinct, and the trace
condition holds: $4a_1+2a_2=2b_1+2b_2+2b_3+
\sum_{i=1}^6 c_i$. Let us set up the following notation:\\

\begin{equation}
\label{equation:extra_of_S:p,q}
p_{ij}=b_i+c_j-a_1,~~~q_{ijk}=2c_i+2c_j+2c_k-\frac12 \,
\sum_{l=1}^6 c_l.
\end{equation}

 We now define the matrices $\B$ and $\C$ by setting

\begin{center}
  \begin{equation}
  \label{equation:extra_of_S:B,C}
    \begin{array}{l}
      \B= \[ \begin{array}{cc|cc|cc}
        b_1 & 0 & -\frac{p_{16}\, q_{245}}{c_3-c_4} & p_{16} &
        \frac{p_{16}\, q_{245}}{c_1-c_2} & p_{16} \\
         & & & & & \\
        0 & b_1 & \frac{p_{15}\, q_{235}\, q_{246}}{(c_3-c_4)(c_5-c_6)} &
        -\frac{p_{15}\, q_{236}}{c_5-c_6} & -\frac{p_{15}\, q_{236}\, q_{246}}
        {(c_1-c_2)(c_5-c_6)} & -\frac{p_{15}\, q_{235}}{c_5-c_6} \\
         & & & & & \\
        \hline
         & & & & & \\
        0 & 0 & b_2 & 0 & -\frac{p_{24}\, q_{236}}{c_1-c_2} & p_{24} \\
         & & & & & \\
        0 & 0 & 0 & b_2 & -\frac{p_{23}\, q_{245}\, q_{246}}{(c_1-c_2)(c_3-c_4)} &
        \frac{p_{23}\, q_{235}}{c_3-c_4} \\
         & & & & & \\
        \hline
         & & & & & \\
        0 & 0 & 0 & 0 & b_3 & 0 \\
         & & & & & \\
        0 & 0 & 0 & 0 & 0 & b_3
      \end{array} \],
    \\
    \\
    \\
    \\
      \C = \[ \begin{array}{cc|cc|cc}
        c_6 & 0 & 0 & 0 & 0 & 0 \\
         & & & & & \\
        0 & c_5 & 0 & 0 & 0 & 0 \\
         & & & & & \\
        \hline
         & & & & & \\
        -\frac{p_{24}\, q_{236}}{c_5-c_6} & -p_{24} & c_4 & 0 & 0 & 0 \\
         & & & & & \\
        -\frac{p_{23}\, q_{235}\, q_{246}}{(c_3-c_4)(c_5-c_6)} &
        -\frac{p_{23}\, q_{245}}{c_3-c_4} & 0 & c_3 & 0 & 0 \\
         & & & & & \\
        \hline
         & & & & & \\
        -\frac{p_{32}\, q_{235}}{c_5-c_6} & -p_{32} &
        \frac{p_{32}\, q_{235}}{c_3-c_4} & -p_{32} & c_2 & 0 \\
         & & & & & \\
        \frac{p_{31}\, q_{236}\, q_{246}}{(c_1-c_2)(c_5-c_6)} &
        \frac{p_{31}\, q_{245}}{c_1-c_2} &
        \frac{p_{31}\, q_{245}\, q_{246}}{(c_1-c_2)(c_3-c_4)} &
        -\frac{p_{31}\, q_{236}}{c_1-c_2} & 0 & c_1
      \end{array} \].
    \end{array}
  \end{equation}
\end{center}

\noindent It is clear that $\B$ and $\C$ are diagonalizable and
that their spectra are $\s(\B)=\{b_1,b_1,b_2,b_2,b_3,b_3\}$,
$\s(\C)=\{c_1,c_2,c_3,c_4,c_5,c_6\}$.

\begin{theorem}
\label{thm:extra_of_S:A}
For $\B$ and $\C$ as above, let $\A=\B+\C$. Then $\A$ is diagonalizable and \\
$\s(\A)=\{a_1,a_1,a_1,a_1,a_2,a_2\}$.
\end{theorem}

\begin{lemma}
\label{lemma:extra_of_S:eigenbasis_of_C} The following are the
eigenvectors of the matrix $\C$ ($\v_i$ corresponds to the
eigenvalue $c_i$):

$$
\begin{array}{lll}

\v_1 & = & \( 0,0,0,0,0,1 \) \\
 & & \\
\v_2 & = & \( 0,0,0,0,1,0 \) \\
 & & \\
\v_3 & = & \( 0,0,0,1,\frac{p_{32}}{c_2-c_3},\frac{p_{31}\,q_{236}}{(c_1-c_2)(c_1-c_3)} \) \\
 & & \\
\v_4 & = & \( 0,0,1,0,-\frac{p_{32}\,q_{235}}{(c_2-c_4)(c_3-c_4)},
-\frac{p_{31}\,q_{245}\,q_{246}}{(c_1-c_2)(c_1-c_4)(c_3-c_4)} \) \\
 & & \\
\v_5 & = & \(
0,1,\frac{p_{24}}{c_4-c_5},\frac{p_{23}\,q_{245}}{(c_3-c_4)(c_3-c_5)},
-\frac{p_{25}\,p_{32}\,q_{234}}{(c_2-c_5)(c_3-c_5)(c_4-c_5)},
-\frac{p_{25}\,p_{31}\,q_{245}\,q_{256}}{(c_1-c_2)(c_1-c_5)(c_3-c_5)(c_4-c_5)}
\), \\
 & & \\
\v_6 & = & \( 1,0,\frac{p_{24}\,q_{236}}{(c_4-c_6)(c_5-c_6)},
\frac{p_{23}\,q_{235}\,q_{246}}{(c_3-c_4)(c_3-c_6)(c_5-c_6)},
-\frac{p_{26}\,p_{32}\,q_{234}\,q_{235}}{(c_2-c_6)(c_3-c_6)(c_4-c_6)(c_5-c_6)}, \right. \\
 & & \left. -\frac{p_{26}\,p_{31}\,q_{236}\,q_{246}\,q_{256}}{(c_1-c_2)(c_1-c_6)(c_3-c_6)
(c_4-c_6)(c_5-c_6)} \).

\end{array}
$$
\end{lemma}

 Let us define a scalar product on $V$ by setting $\<\v_i,\v_j\>=0$ for $i\ne j$
and setting

\begin{equation}
\label{equation:extra:form}
\begin{array}{lcc}
\<\v_1,\v_1\>= &
-\frac{(c_1-c_2)(c_1-c_3)(c_1-c_4)(c_1-c_5)(c_1-c_6)}
{p_{11}\, p_{21}\, p_{31}} & \\
 & & \\
\<\v_2,\v_2\>= & \frac{(c_2-c_3)(c_2-c_4)(c_2-c_5)(c_2-c_6)}
{(c_2-c_1)\, p_{12}\, p_{22}\, p_{32}}\times & \frac{q_{134}\,
q_{135}\, q_{136}\, q_{145}}
{q_{146}\, q_{156}} \\
 & & \\
\<\v_3,\v_3\>= & \frac{(c_3-c_4)(c_3-c_5)(c_3-c_6)\, p_{33}}
{(c_3-c_1)(c_3-c_2)\, p_{13}\, p_{23}}\times & \frac{q_{124}\,
q_{125}\, q_{126}\, q_{145}}
{q_{146}\, q_{156}} \\
 & & \\
\<\v_4,\v_4\>= & \frac{(c_4-c_5)(c_4-c_6)\, p_{34}}
{(c_4-c_1)(c_4-c_2)(c_4-c_3)\, p_{14}\, p_{24}}\times &
\frac{q_{123}\, q_{125}\, q_{126}\, q_{135}\, q_{136}}{q_{156}} \\
 & & \\
\<\v_5,\v_5\>= & \frac{(c_5-c_6)\, p_{25}\, p_{35}}
{(c_5-c_1)(c_5-c_2)(c_5-c_3)(c_5-c_4)\, p_{15}}\times &
\frac{q_{123}\, q_{124}\, q_{126}\, q_{134}\, q_{136}}{q_{146}} \\
 & & \\
\<\v_6,\v_6\>= & \frac{p_{26}\, p_{36}}
{(c_6-c_1)(c_6-c_2)(c_6-c_3)(c_6-c_4)(c_6-c_5)\, p_{16}}\times &
q_{123}\, q_{124}\, q_{125}\, q_{134}\, q_{135}\, q_{145}
\end{array}
\end{equation}

  The sets $S'((4,2),(2,2,2),(1^6))$ and $S''((4,2),(2,2,2),(1^6))$ are constructed
from (\ref{equation:extra:form}) similarly to the hypergeometric
case (see page \pageref{equation:hg:form}
) and have the same properties.

\begin{theorem}
\label{thm:extra:form} The operators $\A$, $\B$, and $\C$ are
self-adjoint with respect to the scalar product
(\ref{equation:extra:form}).
\end{theorem}

 Now suppose that the eigenvalues of the matrices $\A$, $\B$, and $\C$
are real numbers.
Let $b_1>b_2>b_3$ and $c_1>c_2>\cdots >c_6$.

\begin{theorem}
\label{thm:extra:sign_of_form} The form $\<*,*\>$ defined by
(\ref{equation:extra:form}) is sign-definite precisely in the
following two situations:

$$
\begin{array}{|l|l|}
\hline
 & \\
\begin{array}{l}
b_1+c_4 > a_1 > b_1+c_5           \\
b_2+c_2 > a_1 > b_2+c_3           \\
\phantom{b_2+c_2 \ge} a_1 > b_3+c_1 \\
 \\
c_1+c_4+c_5 > c_2+c_3+c_6 \\
c_1+c_3+c_6 > c_2+c_4+c_5 \\
c_2+c_3+c_5 > c_1+c_4+c_6
\end{array} &
\begin{array}{l}
b_3+c_2 > a_1 > b_3+c_3           \\
b_2+c_4 > a_1 > b_2+c_5           \\
b_1+c_6 > a_1 \\
 \\
c_1+c_4+c_5 >  c_2+c_3+c_6 \\
c_1+c_3+c_6 >  c_2+c_4+c_5 \\
c_2+c_3+c_5 >  c_1+c_4+c_6
\end{array} \\
 & \\
\hline
\end{array}
$$
\medskip

\noindent If the form $\epsilon\, \<*,*\>$ is positive-definite
for $\epsilon=\pm 1$, then $\epsilon=sign(a_1-a_2)$. If the
inequalities of the first column hold, then $a_1>a_2$. If the
inequalities of the second column hold, then $a_1<a_2$.
\end{theorem}

\section{Proofs and More Results}
\label{sec:proofs}

 In this section we prove theorems from Section \ref{sec:main_results}.
In the process, some new results are obtained. \\

 The following simple observation is helpful in this section.
If we replace a triple $(\A=\B+\C,\B,\C)$ by the triple

\begin{equation}
\label{equation:tildeA,B,C}
\begin{array}{l}
\tA=k\,\A+\theta\,\Id, \\
\tB=k\,\B+\phi\,\Id,  \\
\tC=k\,\C+(\theta - \phi)\,\Id \\
k,\theta,\phi\in\mathbb{C}~~\mbox{and}~k\ne 0,
\end{array}
\end{equation}

\noindent  then we still have $\tA = \tB + \tC$. This
transformation changes neither irreducibility nor rigidity of the
triple. If, say, $\v$ is an eigenvector of $\B$ with the
eigenvalue $b$, then $\v$ is an eigenvector of $\tB$ with the
eigenvalue $k\,b+\phi$. If $\A$, $\B$, and $\C$ were self-adjoint
with respect to a symmetric bilinear form , then $\tA$, $\tB$, and
$\tC$ are self-adjoint with respect to the form as well.

\subsection{Hypergeometric Family}
\label{subsec:proofs:hg}

An affine transformation (\ref{equation:tildeA,B,C}) with $k=1$,
$\theta=-a_2$, and $\phi=-a_2/2$ normalizes $\A$ to $\tA$ such
that the eigenvalue of $\tA$ of multiplicity $m-1$ is $0$. So,
without loss of generality, we can assume
that $a_2=0$. Now let us prove Theorem \ref{thm:hg:A}.\\

\noindent {\it Proof of Theorem \ref{thm:hg:A} ---}
\label{proof:hg:A} Consider the matrix $\A$ (with $a_2=0$). Here
is an example with $m=5$.

$$
\A=\left[
\begin{array}{lllll}
b_1+c_5 & b_1+c_5 & b_1+c_5 & b_1+c_5 & b_1+c_5 \\
b_2+c_4 & b_2+c_4 & b_2+c_4 & b_2+c_4 & b_2+c_4 \\
b_3+c_3 & b_3+c_3 & b_3+c_3 & b_3+c_3 & b_3+c_3 \\
b_4+c_2 & b_4+c_2 & b_4+c_2 & b_4+c_2 & b_4+c_2 \\
b_5+c_1 & b_5+c_1 & b_5+c_1 & b_5+c_1 & b_5+c_1
\end{array}
\right]
$$

\noindent Now $\A$ has rank $1$, and its image is the linear span
of the vector ${\bf i}=(b_1+c_m,b_2+c_{m-1},\cdots ,b_m+c_1)$.
$\A\, {\bf i} = a_1 {\bf i} \ne 0$. Thus, $\A$ is diagonalizable,
and $\s(\A)=\{\sum_{i=1}^m (b_i+c_i),0,0,\cdots ,0\}$. Thus,
before the normalizing affine transformation we had
$\s(\A)=\{a_1,\underbrace{a_2,\cdots ,a_2}_{m-1~times}\}$. \endproof

    In our normalized version,
\begin{equation}
\label{equation:hg:normal:B,C}
\begin{array}{ll}
B_{ij}=\left\{
\begin{array}{lll}
0             & ,if & i<j \\
b_i           & ,if & i=j \\
b_i+c_{m+1-i} & ,if & i>j
\end{array},
\right. &

C_{ij}=\left\{
\begin{array}{lll}
0             & ,if & i>j \\
c_{m+1-i}     & ,if & i=j \\
b_i+c_{m+1-j} & ,if & i<j
\end{array}.
\right.
\end{array}
\end{equation}

 We are ready to prove Lemma \ref{lemma:hg:eigenbasis_of_B},
that is to show that for every $i=1,\cdots ,m$ the vector
$\v_i=(v_i^1,\cdots ,v_i^{i-1},1,0,\cdots ,0)$ with

\begin{equation}
\label{equation:eigenvectors_of_B}
v_i^j=\displaystyle{\frac{b_j+c_{m+1-j}}{b_i-b_j}}\prod\limits_{k=1}^{i-j-1}
\displaystyle{\frac{b_i+c_{m+1-j-k}}{b_i-b_{j+k}}}i~~~~(j=1,\cdots
,i-1)
\end{equation}

\noindent is an eigenvector of $\B$ with the eigenvalue $b_i$. \\

\noindent {\it Proof of Lemma \ref{lemma:hg:eigenbasis_of_B} ---}
Remembering the definition (\ref{equation:hg:normal:B,C}) of $\B$,
we need to show the following equality for all $j<i$:
$(b_i-b_j)v_i^j=(b_j+c_{m+1-j})\sum_{k=j+1}^i v_i^k$, or
equivalently

\begin{equation}
\label{equation:hg:sum _of_v_i^k} \sum_{j+1}^i v_i^k =
\prod_{k=1}^{i-j-1} \frac{b_i+c_{m+1-j-k}}{b_i-b_{j+k}}.
\end{equation}

\noindent This identity becomes obvious once we rewrite
(\ref{equation:eigenvectors_of_B}) as

\begin{equation}
\label{equation:hg:telescoping} v_i^j = \prod_{k=1}^{i-j}
\frac{b_i+c_{m+1-j-k}}{b_i-b_{j+k}} - \prod_{k=1}^{i-j-1}
\frac{b_i+c_{m+1-j-k}}{b_i-b_{j+k}},
\end{equation}

\noindent and use telescoping. \endproof

 Now we are ready to prove Theorem \ref{thm:hg:form}, that is to show
that the operators $\A$, $\B$, and $\C$ are self-adjoint with
respect to the scalar product (\ref{equation:hg:form}). \\

\noindent {\it Proof of Theorem \ref{thm:hg:form} ---} The operator
$\B$ is self adjoint with respect to the scalar product by
construction. To show that $\A$ is self-adjoint, we have to show
that $\<\A\, \v_i,\v_j\>=\<\v_i,\A\, \v_j\>$.  As we have seen,
$\A$ has a one-dimensional image spent by the vector ${\bf
i}=(b_1+c_m,b_2+c_{m-1},\cdots ,b_m+c_1)$. Namely, for any vector
${\bf k}=(k_1,k_2,\cdots ,k_m)$, $\A\,{\bf k}=\(\sum_{i=1}^m
k_i\){\bf i}$. In particular, $\A\,\v_i=\(\sum_{j=1}^m v_i^j\){\bf
i}$. In view of (\ref{equation:hg:sum _of_v_i^k}), we have

\begin{equation}
\label{equation:hg:Av} \A\,\v_i = \( \prod_{k=1}^{i-1}
\frac{b_i+c_{m+1-k}}{b_i-b_k} \)\, {\bf i}.
\end{equation}

 It will be convenient to introduce the following notation:

\begin{equation}
\label{equation:s_i,x_i}
s_i=\prod\limits_{k=1}^{i-1}\displaystyle{\frac{b_i+c_{m+1-k}}{b_i-b_k}};~~
x_i=(b_i+c_{m+1-i})\prod\limits_{k=1}^{m-i}\displaystyle{\frac{b_i+c_{m+1-i-k}}
{b_i-b_{i+k}}}.
\end{equation}

\noindent Then $\A\,\v_i=s_i\, {\bf i}$ and
(\ref{equation:hg:form}) can be rewritten as
$\<\v_i,\v_i\>=s_i/x_i$. Now the desired equality
$\<\A\v_i,\v_j\>=s_ix_j\<\v_j,\v_j\>=s_ix_js_j/x_j=s_is_j=\<\v_i,\A\v_j\>$
becomes a consequence of the following lemma. \endproof

\begin{lemma}
\label{lemma:x_i} $\sum\limits_{i=1}^m x_i \v_i={\bf i}.$
\end{lemma}

\proof We have $\v_i=(v_i^1,\cdots ,v_i^{i-1},1,0\cdots ,0)$.
Thus, to prove the lemma we have to prove that the identity
$\sum_{j=i}^m x_j\, v_j^i = b_i+c_{m+1-i}$ holds for $1\le i\le
m$. This is equivalent to  $x_i=b_i+c_{m+1-i}-\sum_{k=1}^{m-i}
x_{i+k} v_{i+k}^i$. This identity after minor simplification
becomes

$$
\sum_{j=i}^m \frac{\prod\limits_{k=1}^{m-i} (b_j+c_k)}
{\prod\limits_{k=i\atop{k\ne j}}^m (b_j-b_k)}=1.
$$

\noindent Let us set $n=m-i+1$; $x_1=b_i$, $x_2=b_{i+1},\cdots
,x_n=b_m$; $y_1=-c_1$, $y_2=-c_2,\cdots ,y_{n-1}=-c_{m-i}$. This
change of variables transforms the last identity into identity
(\ref{equation:hg_id_2}) which we prove in the Appendix.
\endproof

 Now we prove Theorem \ref{thm:irr} for the hypergeometric family. That is,
we show that if the vector $\(\s(\A),\s(\B),\s(\C)\)$ lies in
$S''((1,m-1),(1^m),(1^m))$, then the triple $(\A,\B,\C)$ is
irreducible. \label{proof:hg:irr} \\

\noindent {\it Proof of Theorem \ref{thm:irr} ---}Suppose that the
triple preserves a non-trivial subspace of $V$. Then this subspace
is spanned by some of the eigenvectors $\v_j$ of $\B$. But
$\A\,\v_j=\sum_{i=1}^m s_j x_i \v_i$. If
$\(\s(\A),\s(\B),\s(\C)\)$ lies in $S''((1,m-1),(1^m),(1^m))$,
then all the coefficients $s_jx_i$ are non-zero, so $\A$ does not
preserve any such proper subspace. Thus, the triple $(\A,\B,\C)$
is irreducible. \endproof

 Let us prove Theorem \ref{thm:hg:sign_of_form}, that is determine
the inequalities on the real spectra of $\A$, $\B$, and $\C$ which
make the form (\ref{equation:hg:form}) sign-definite. Recall that
we are working with the normalized version $a_2=0$. Then the trace
identity gives us $a_1=\sum_{i=1}^m (b_i+c_i)$. Also, it is an
assumption of Theorem \ref{thm:hg:sign_of_form} that
$b_1>b_2>\cdots >b_m$
and $c_1>c_2>\cdots >c_m$. \\

\noindent {\it Proof of Theorem \ref{thm:hg:sign_of_form} ---} It
immediately follows from Theorem \ref{thm:hg:form} that

$$
\mbox{sign}(\<\v_i,\v_i\>)=(-1)^{i-1}
\mbox{sign}\(\prod\limits_{j=1}^m (b_i+c_j)\).
$$

\noindent Construct an $m\times m$ matrix T where
$T_{i,j}=b_i+c_j$. Notice that $T_{i,j}>T_{i,j+1}$ and
$T_{i,j}>T_{i+1,j}$ for all $i$ and $j$. Then
$\mbox{sign}(\<\v_i,\v_i\>)=(-1)^{i-1}\times (-1)^{\# \left\{j:
T_{i,j}<0 \right\}}$. Thus, to keep the sign constant, $\#
\left\{j: T_{i,j}<0 \right\}$ must differ from $\# \left\{j:
T_{i+1,j}<0 \right\}$ by an odd number for all $m$ rows of $T$.
This gives us only two possibilities: either
$T_{i,m-i}>0>T_{i,m+1-i}$, or $T_{i,m+1-i}>0>T_{i,m+2-i}$. Here is
a picture which illustrates the two situations for $m=5$. The line
separates $T_{i,j}>0$ from $T_{i,j}<0$.

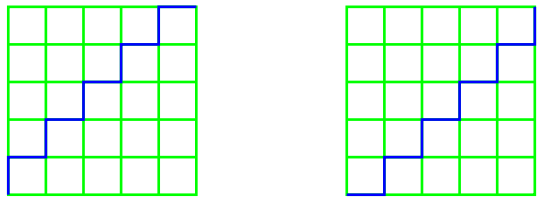
\begin{center}
\pspicture[](0,0)(5,3)
\psset{unit=0.5cm,linewidth=1pt,linecolor=green} \rput(-3,0.5){
  \rput(0,0){
    \pspolygon(0,0)(5,0)(5,5)(0,5)
    \psline(1,0)(1,5)
    \psline(2,0)(2,5)
    \psline(3,0)(3,5)
    \psline(4,0)(4,5)
    \psline(0,1)(5,1)
    \psline(0,2)(5,2)
    \psline(0,3)(5,3)
    \psline(0,4)(5,4)
    \psset{linecolor=blue}
    \psline(0,0)(0,1)(1,1)(1,2)(2,2)(2,3)(3,3)(3,4)(4,4)(4,5)(5,5)
  }
  \rput(9,0){
    \pspolygon(0,0)(5,0)(5,5)(0,5)
    \psline(1,0)(1,5)
    \psline(2,0)(2,5)
    \psline(3,0)(3,5)
    \psline(4,0)(4,5)
    \psline(0,1)(5,1)
    \psline(0,2)(5,2)
    \psline(0,3)(5,3)
    \psline(0,4)(5,4)
    \psset{linecolor=blue}
    \psline(0,0)(1,0)(1,1)(2,1)(2,2)(3,2)(3,3)(4,3)(4,4)(5,4)(5,5)
  }
}
\endpspicture
\end{center}

\noindent The sum of $T_{i,j}$ along the non-main diagonal of $T$
equals $a_1$. Thus, $T_{i,m-i}>0>T_{i,m+1-i}$ forces $a_1<0$ and
$T_{i,m+1-i}>0>T_{i,m+2-i}$ forces $a_1>0$. \endproof

\begin{remark}
\label{remark:hg:pos:standard_form_of_HGM} {\rm Let the
eigenvalues of $\A$, $\B$, and $\C$ be real numbers, and let the
form (\ref{equation:hg:form}) be positive-definite. Then in the
basis ${\tilde \e}_i=\v_i/\sqrt{\<\v_i,\v_i\>}$, the form becomes
standard ($\< {\tilde \e}_i, {\tilde \e}_j \> = \delta_{ij}$). Let
${\tilde \A}$, ${\tilde \B}$, and ${\tilde \C}$, be the matrices
$\A$, $\B$, and $\C$ in the basis ${\tilde \e}_1,{\tilde
\e}_2,\cdots ,{\tilde \e}_m$. Then for $i,j=1,2,\cdots ,m$,
$$
\begin{array}{l}
{\tilde A}_{ij}=\sqrt{x_ix_js_is_j}, \\
{\tilde B}_{ij}=\delta_{ij}b_i, \\
{\tilde C}_{ij}={\tilde A}_{ij}-{\tilde B}_{ij}.
\end{array}
$$
}
\end{remark}

%

\subsection{Even Family}
\label{subsection:proofs:even}

In order to make proofs simpler, let us normalize $\A$, $\B$, and
$\C$ so that they become traceless and $a_1=1$, $a_2=-1$. The
affine transformation (\ref{equation:tildeA,B,C}) with
$k=2/(a_1-a_2)$,
$\theta=-(a_1+a_2)/2$, and $\phi=-(b_1+(m-1)b_2+m\,b_3)/(2m)$ does the job. \\

 Let us prove Theorem \ref{thm:even:A}. In our normalized version, we
have to prove that $\A$ is diagonalizable and that
$\s(\A)=\{\underbrace{1,\cdots ,1,}_{m~times}\underbrace{-1,\cdots
,-1}_{m~times}\}$.

\noindent {\it Proof of Theorem \ref{thm:even:A} ---}
First, let us prove that $\A^2=\Id$. For that, we have to
prove the following eleven identities.

\begin{enumerate}
\item
The identity $\sum_{l=1}^{2m} A_{1l}\, A_{l1}=1$ with the help of
identity (\ref{equation:hg_id_2}) can be reduced to the identity

$$
\begin{array}{l}
((m-1)b_2+mb_3+c_1+c_2+\cdots +c_{2m-1})^2=
\sum\limits_{i=1}^{m-1} \displaystyle{\frac{\prod\limits_{j=1}^m
(b_2+b_3+c_j+c_{m+i})}{\prod\limits_{j=1\atop j\ne i}^{m-1}(c_{m+i}-c_{m+j})}} + \\
\sum\limits_{i=1}^m (b_3+c_i)^2\,
\displaystyle{\frac{\prod\limits_{j=1}^{m-1}
(b_2+b_3+c_i+c_{m+j})}{\prod\limits_{j=1\atop j\ne i}^m
(c_i-c_j)}}.
\end{array}
$$

\noindent For $1\le i\le m-1$, let us set $x_i=c_{m+i}+b_2$. For
$1\le i\le m$, let us set $y_i=c_i+b_3$. Then the last identity
becomes

\begin{equation}
\label{equation:hg_id_4'} \(\sum_{i=1}^{m-1}x_i+\sum_{i=1}^m
y_i\)^2=\sum_{i=1}^{m-1} \displaystyle{\frac{\prod\limits_{j=1}^m
(x_i+y_j)}{\prod\limits_{j=1\atop j\ne i}^{m-1} (x_i-x_j)}} +
\sum_{i=1}^m y_i^2 \displaystyle{\frac{\prod\limits_{j=1}^{m-1}
(x_j+y_i)}{\prod\limits_{j=1\atop j\ne i}^m (y_i-y_j)}}.
\end{equation}

\noindent Introducing $x_m=x_{m+1}=0$, we can rewrite the last
term in (\ref{equation:hg_id_4'}) as

$$
\sum_{i=1}^m y_i^2 \displaystyle{\frac{\prod\limits_{j=1}^{m-1}
(x_j+y_i)}{\prod\limits_{j=1\atop j\ne i}^m (y_i-y_j)}}=
\sum_{i=1}^m \displaystyle{\frac{\prod\limits_{j=1}^{m+1}
(y_i+x_j)}{\prod\limits_{j=1\atop j\ne i}^m (y_i-y_j)}}.
$$

Now we can prove identity (\ref{equation:hg_id_4'}) with the help
of identity (\ref{equation:hg_id_4}) from the Appendix.

\item
For $1\le i\le m-1$, the identity $\sum_{l=1}^{2m} A_{1,l}\,
A_{l,1+i}=0$ after some simplification becomes

$$
\sum_{j=1}^m \((b_3+c_j)^2-1\) \frac{\prod\limits_{k=m+1\atop k\ne
2m-i}^{2m-1} q_{j,k}} {\prod\limits_{k=1\atop k\ne j}^m (c_j-c_k)}
= \sum_{j=1}^m (c_j+b_3)+\sum_{j=m+1\atop j\ne 2m-i}^{2m-1}
(c_j+b_2).
$$

\noindent For $1\le j\le m$, set $x_j=b_3+c_j$. For $1\le j\le
m-1$ and $j\ne m-i$, set $y_j=-(b_2+c_{m+j})$. The identity

$$
\sum_{j=1}^m \frac{\prod\limits_{k=m+1\atop k\ne 2m-i}^{2m-1}
q_{j,k}} {\prod\limits_{k=1\atop k\ne j}^m (c_j-c_k)} = 0
$$

\noindent is equivalent to identity (\ref{equation:hg_id_1}) from
the Appendix. Now the identity to prove becomes

\begin{equation}
\label{equation:eq1} \sum_{j=1}^m (b_3+c_j)^2\,
\frac{\prod\limits_{k=m+1\atop k\ne 2m-i}^{2m-1} q_{j,k}}
{\prod\limits_{k=1\atop k\ne j}^m (c_j-c_k)} = \sum_{j=1}^m
(c_j+b_3)+\sum_{j=m+1\atop j\ne 2m-i}^{2m-1} (c_j+b_2).
\end{equation}

\noindent Introducing $y_{m-i}=y_m=0$, we reduce the last identity
to identity (\ref{equation:hg_id_3}) from the Appendix.

\item
For $1\le j\le m$, the identity $\sum_{l=1}^{2m} A_{1,l}\,
A_{l,m+j}=0$ reduces to identity (\ref{equation:hg_id_3}) from the
Appendix.

\item
For $1\le i\le m-1$, the identity $\sum_{l=1}^{2m} A_{1+i,l}\,
A_{l,1}=0$ reduces to identity (\ref{equation:hg_id_3}) from the
Appendix.

\item
Let $1\le i,j\le m-1$, $i\ne j$. The identity $\sum_{l=1}^{2m}
A_{1+i,l}\, A_{l,1+j}=0$ reduces to identity
(\ref{equation:hg_id_2}) from the Appendix.

\item
For $1\le i\le m-1$, the identity $\sum_{l=1}^{2m} A_{1+i,l}\,
A_{l,1+i}=1$ after some simplification becomes

$$
\sum_{j=1}^m p^{31}_j\,p^{32}_j \frac{\prod\limits_{k=1\atop k\ne
j}^m q_{k,i}\,\prod\limits_{k=m+1\atop k\ne i}^{2m-1}
q_{j,k}}{\prod\limits_{k=1\atop k\ne j}^m (c_j-c_k)\,\prod\limits_
{k=m+1\atop k\ne i}^{2m-1}
(c_i-c_k)}=1-(b_2+c_i)^2+\frac{\prod\limits_{k=1}^m
q_{k,i}}{\prod\limits_{k=m+1\atop k\ne i}^{2m-1} (c_i-c_k)}.
$$

\noindent Recall that in the normalized version
$p_i^{31}=c_i+b_3-1$ and $p_i^{32}=c_i+b_3+1$. For $1\le i\le m$,
let us set $x_i=b_3+c_i$. For $1\le i\le m-1$, let us set
$y_i=-b_2- c_{m+i}$. The above identity splits into two
homogeneous identities: one of degree $0$ and the other of $2$ (in
$x_i$ and $y_j$). The first is equivalent to identity
(\ref{equation:hg_id_5}) from the Appendix. The second is equivalent
to identity (\ref{equation:hg_id_6}) from the Appendix.

\item
For $1\le i\le m-1$ and $1\le j\le m$, the identity
$\sum_{l=1}^{2m} A_{1+i,l}\, A_{l,m+j}=0$ reduces to the trivial
identity $\frac{q_{2m-i,m+1-j}}{q_{2m-i,m+1-j}}-1=0$.

\item
For $1\le i \le m$, the identity $\sum_{l=1}^{2m} A_{m+j,l}\,
A_{l,1}=0$ reduces to the identity

$$
\sum_{j=1}^{m-1} \frac{\prod\limits_{k=1\atop k\ne m+1-i}^m
q_{k,2m-j}} {\prod\limits_{k=m+1\atop k\ne 2m-j}^{2m-1}
(c_{2m-j}-c_k)}=-b_1-b_3-c_{m+1-i}-c_{2m}.
$$

\noindent The latter follows from identity
(\ref{equation:hg_id_3}) of the Appendix and from the fact that the
normalized $\B$ and $\C$ are traceless.

\item
Let $1\le i \le m$ and $1\le j \le m-1$. The identity
$\sum_{l=1}^{2m} A_{m+j,l}\, A_{l,1+i}=0$ reduces to the trivial
identity $\frac{q_{2m-j,m+1-i}}{q_{2m-j,m+1-i}}-1=0$.

\item
Let $1\le i\ne j\le m$. The identity $\sum_{l=1}^{2m} A_{m+j,l}\,
A_{l,m+i}=0$ reduces to identity (\ref{equation:hg_id_2}) from the
Appendix.

\item
Let $1\le i\le m$. To prove that $\sum_{l=1}^{2m} A_{m+i,l}\,
A_{l,m+i}=1$, we set $x_1=b_3+c_1$, $x_2=b_3+c_2,\cdots
,x_m=b_3+c_m$; $y_1=b_2+c_{m+1}$, $y_2=b_2+c_{m+2},\cdots
,y_{m-1}=b_2+c_{2m-1}$. This reduces the identity in question to
identity (\ref{equation:hg_id_7}) of the Appendix.
\end{enumerate}

 Now we are ready to prove that if $(\s(\A),\s(\B),\s(\C))$ is a point of
$S''((m,m),(1,m-1,m),(1^{2m}))$, then
$\s(\A)=\{\underbrace{1,\cdots ,1}_{m~times},\underbrace{-1,\cdots
,-1}_{m~times}\}$. We know that $\A^2=\Id$. Thus, $\A$ is
diagonalizable and the eigenvalues of $\A$ are $1$ and $-1$. For
$1\le i\le m$, let us set

\begin{equation}
\label{equation:evect_of_A}
\begin{array}{lcl}
{\bf a}_i^+    & = & \( \A+\Id \) \e_{m+i}\\
{\bf a}_i^- & = & \( \A-\Id \) \e_{m+i}.
\end{array}
\end{equation}

\noindent Then $\( \A-\Id \) {\bf a}_i^+ = \( \A-\Id \) \( \A+\Id
\) \e_{m+i}=0$ and $\( \A+\Id \) {\bf a}_i^- = \( \A+\Id \) \(
\A-\Id \) \e_{m+i} = 0$. If we take a look at the matrix $\A$, we
see that the condition $\(\s(\A), \s(\B), \s(\C)\) \in
S''((m,m),(1,m-1,m),(1^{2m}))$ guarantees that the vectors $\{
{\bf a}_i^+ \}_{i=1,\cdots ,m}$ are linearly independent as well
as the vectors $\{ {\bf a}_i^- \}_{i=1,\cdots ,m}$. (If this
condition is violated, then ${\bf a}_i^+$ and ${\bf a}_i^-$ are
not necesserily linearly independent. For example, if
$p^{31}_i=0$, then ${\bf a}_i^-=0$.)
\endproof

 Let us prove Lemma \ref{lemma:even:eigenbasis_of_C}, that is compute the
coordinates $v_i^j$ of the eigenvectors $\v_i$ of the matrix $\C$.

\noindent {\it Proof of Lemma \ref{lemma:even:eigenbasis_of_C} ---}
The only non-trivial part of the lemma is formula
(\ref{equation:even:v_{2m}^{m+j}}). A direct computation gives

$$
v_{2m}^{m+i}=\frac{1}{c_{m+1-i}-c_{2m}}\left[ -C_{m+i,1}+
\sum_{j=1}^{m-1}\frac{C_{m+i,j+1}\, C_{j+1,1}}{c_{2m-j}-c_{2m}}
\right],
$$

\noindent where $ C_{j+1,1}$ is given by
(\ref{equation:even:C:IV}), $ C_{m+i,1}$ is given by
(\ref{equation:even:C:VII}), and $ C_{m+i,j+1}$ is given by
(\ref{equation:even:C:VIII}), see page
\pageref{equation:even:C:IV}. Comparing the formula for
$v_{2m}^{m+i}$ given by (\ref{equation:even:v_{2m}^{m+j}}) to the
right hand side of the last formula, we obtain an identity which
reduces to identity (\ref{equation:hg_id_2}) from the Appendix.
\endproof

 The following lemma expresses the vectors $\e_i$ of the standard basis in terms
of the eigenvectors $\v_j$ of the matrix $\C$.

\begin{lemma}
\label{lemma:even:Vinv}
\begin{enumerate}
\item
$\e_1=\v_{2m}-\sum\limits_{i=1}^{m-1}
\displaystyle{\frac{\prod\limits_{k=1}^{m-i} q_{k,m+i}}
{\prod\limits_{k=m+1+i}^{2m} (c_{m+i}-c_k)}}\, \v_{m+i} -
\sum\limits_{i=1}^m \displaystyle{\frac{p^{32}_i
\prod\limits_{k=1\atop k\ne i}^{2m-i}q_{ik}}
{\prod\limits_{k=1+i}^{2m}(c_i-c_k)}}\, \v_i$,

\item
$\e_{1+i} = \v_{2m-i}+\sum\limits_{j=1}^m
(-1)^{m+1-i}\,\displaystyle{\frac{p^{32}_j} {c_j-c_{2m-i}}}\,
\displaystyle{\frac{\prod\limits_{k=m+1\atop k\ne 2m-i}^{2m-j}
q_{j,k} \,\prod\limits_{k=1+i\atop k\ne j}^m
q_{k,2m-i}}{\prod\limits_{k=1+j}^m
(c_j-c_k)\,\prod\limits_{k=m+1}^{2m-1-i} (c_k-c_{2m-i})}} \, \v_j$\\
for $i=1,2,\cdots ,m-1$.

\item
$\e_{m+i}=\v_{m+1-i}$ for $i=1,2,\cdots ,m$.
\end{enumerate}
\end{lemma}

\proof  For $1\le i\le m$, let

\begin{equation}
\label{equation:e_1^i} e_1^i = \sum\limits_{i=1}^m
\displaystyle{\frac{p^{32}_i \prod\limits_{k=1\atop k\ne
i}^{2m-i}q_{ik}}{\prod\limits_{k=1+i}^{2m}(c_i-c_k)}}
\end{equation}

\noindent be the $i$-th coordinate of $\e_1$ in the basis $\v_i$.
The formula (\ref{equation:e_1^i}) is the only nontrivial part of
the lemma. To prove it we have to show that

\begin{equation}
\label{equation:trr}
e_1^i+v_{2m}^{2m+1-i}-\sum_{k=m+1}^{2m-1}v_k^{2m+1-i}\,v_{2m}^{2m+1-k}=0,~
\mbox{where}~ i=1,2,\cdots ,m.
\end{equation}

\noindent Let us use the formulas for $v_i^j$ given in Lemma
\ref{lemma:even:eigenbasis_of_C} and the following change of
variables: for a fixed $1\le i \le m$, set $x_0=c_i+b_3$,
$x_1=c_{m+1}+b_3$, $x_2=c_{m+2}+b_3, \cdots ,x_m=c_{2m}+b_3$,
$y_1=-b_2-c_1$, $y_2=-b_2-c_2,\cdots ,y_{i-1}=-b_2-c_{i-1}$,
$y_i=-b_2-c_{i+1},\cdots ,y_{m-1}=-b_2-c_m$. This reduces
(\ref{equation:trr}) to identity (\ref{equation:hg_id_1}) from
Appendix.
\endproof

 To prove Theorem \ref{thm:even:form}, we have to compute eigenvectors of
the matrix $\B$ (see page \pageref{picture:even:B} for its
description). Let us recall that $\B$ is diagonalizable. Let us
call $\w_1$ the eigenvector of $\B$ corresponding to the
eigenvalue $b_1$, $\w_2,\cdots ,\w_m$ the eigenvectors of $\B$
corresponding to the eigenvalue $b_2$, and $\w_{m+1},\cdots
,\w_{2m}$ the eigenvectors of $\B$ corresponding to the eigenvalue
$b_3$.

\begin{lemma}
\label{lemma:even:evec_of_B}
\begin{enumerate}
\item We have $\w_1=\e_1$.

\item For $1\le i\le m-1$, we have $\w_{1+i}=-\frac{B_{1,1+i}}{b_1-b_2}\e_1+\e_{1+i}$.

\item For $1\le i\le m$, we have $\w_{m+i}=x^{m+i}\, \e_1-\sum\limits_{j=1}^{m-1}\displaystyle
{\frac{B_{j+1,m+i}}{b_2-b_3}}\e_{j+1}+\e_{m+i}$, where

\begin{equation}
\label{equation:even:B:x} x^{m+i}=\frac{(-1)^{m+1-i}
p^{31}_{m+1-i}(b_1+b_2+c_{m+1-i}+c_{2m})}{(b_1-b_3)(b_2-b_3)}
\frac{\prod\limits_{k=m+i}^{2m-1}
q_{m+1-i,k}}{\prod\limits_{k=1}^{m-i} (c_k-c_{m+1-i})}.
\end{equation}

\end{enumerate}
\end{lemma}

Here is an example with $m=3$.

\begin{example}
\label{example:even:evec_of_B:m=3}
$$
\begin{array}{lcl}
\w_1 & = & \(1,0,0,0,0,0\), \\
 & & \\
\w_2 & = & \(\frac{q_{25}\,q_{35}}{(c_4-c_5)(b_1-b_2)},1,0,0,0,0\), \\
 & & \\
\w_3 & = & \(-\frac{q_{34}}{b_1-b_2},0,1,0,0,0\), \\
 & & \\
\w_4 & = & \(-\frac{(b_1+b_2+c_3+c_6) p^{31}_3\,q_{34}\,q_{35}}
{(b_1-b_3)(b_2-b_3)(c_1-c_3)(c_2-c_3)},-\frac{p^{31}_3\,q_{34}\,q_{15}}
{(b_2-b_3)(c_1-c_3)(c_2-c_3)},-\frac{p^{31}_3\,q_{14}\,q_{24}\,q_{35}}
{(b_2-b_3)(c_1-c_3)(c_2-c_3)(c_4-c_5)},1,0,0 \), \\
 & & \\
\w_5 & = & \(\frac{(b_1+b_2+c_2+c_6)
p^{31}_2\,q_{25}}{(b_1-b_3)(b_2-b_3)(c_1-c_2)}, \frac{p^{31}_2\,
q_{15}}{(b_2-b_3)(c_1-c_2)},\frac{p^{31}_2\, q_{14}\, q_{25}}
{(b_2-b_3)(c_1-c_2)(c_4-c_5)},0,1,0 \), \\
 & & \\
\w_6 & = & \(-\frac{(b_1+b_2+c_1+c_6)
p^{31}_1}{(b_1-b_3)(b_2-b_3)},
-\frac{p^{31}_1}{b_2-b_3},-\frac{p^{31}_1\, q_{24}}{(b_2-b_3)(c_4-c_5)},0,0,1 \). \\
\end{array}
$$
\end{example}

\noindent {\it Proof of Lemma \ref{lemma:even:evec_of_B} ---}
All the formulas in Lemma \ref{lemma:even:evec_of_B} are
immediate except for (\ref{equation:even:B:x}). A direct
computation gives

$$
(b_1-b_3)x^{m+i}+\sum_{j=1}^{m-1}
(-1)^{m-j}\frac{\prod\limits_{k=1+j}^{m}q_{k,2m-j}}
{\prod\limits_{k=m+1}^{2m-1-j}(c_k-c_{2m-j})}\,
\frac{B_{1+j,m+i}}{b_2-b_3} + (-1)^{m-i}\, p^{31}_{m+1-i}
\frac{\prod\limits_{k=m+i}^{2m-1} q_{m+1-i,k}}
{\prod\limits_{k=1}^{m-i}(c_k-c_{m+1-i})}=0,
$$

\noindent where $B_{1+j,m+i}$ is given by
(\ref{equation:even:B:VI}). To prove (\ref{equation:even:B:x}), we
have to show that the formula we derive for $x^{m+i}$ from the
equation above equals the formula for $x^{m+i}$ from
(\ref{equation:even:B:x}). This boils down to a proof of identity
(\ref{equation:hg_id_2}) from the Appendix.
\endproof

\begin{lemma}
\label{lemma:even:form_in_good_basis} The following is the matrix
of the scalar product (\ref{equation:even:form}) in the standard
basis $\e_1,\cdots ,\e_{2m}$.\\

\begin{enumerate}
\item
$\<\e_1,\e_1\>= (b_1-b_2)(b_1-b_3)$

\item
$\<\e_1,\e_{1+i}\>= (b_1-b_3)\,
\displaystyle{\frac{\prod\limits_{k=1+i}^m
q_{2m-i,k}}{\prod\limits_{k=m+1}^{2m-1-i} (c_{2m-i}-c_k)}}$~~for~~
$i=1,2,\cdots ,m-1$

\item
$\<\e_1,\e_{m+i}\>=-p^{31}_{m+1-i}\,
\displaystyle{\frac{\prod\limits_{k=m+i}^{2m}
q_{m+1-i,k}}{\prod\limits_{k=1}^{m-i} (c_{m+1-i}-c_k)}}$~~ for~~
$i=1,2,\cdots ,m$

\item
$\<\e_{1+i},\e_{1+j}\>= $\\
$\displaystyle{\frac{\prod\limits_{k=1+i}^m q_{2m-i,k}\,
\prod\limits_{k=1+j}^m q_{2m-j,k}}{\prod\limits_{k=m+1}^{2m-1-i}
(c_{2m-i}-c_k)\,\prod\limits_{k=m+1}^{2m-1-j}
(c_{2m-j}-c_k)}}\left\{1-(b_2-b_3)\, \sum\limits_{r=1}^m c_r\,
\displaystyle{\frac{\prod\limits_{k=m+1\atop{k\ne 2m-i\atop k\ne
2m-j}}^{2m} q_{r,k} \,\prod\limits_{k=m+1\atop{k\ne 2m-i\atop k\ne
2m-j}}^{2m} (c_r-c_k)} {\prod\limits_{k=1\atop k\ne r}^m q_{r,k}\,
\prod\limits_{k=1\atop k\ne r}^m (c_r-c_k)}} \right\}$ for
$i,j=1,2,\cdots ,m-1$, $i\ne j$

\item
$\<\e_{1+i},\e_{1+i}\>= \displaystyle{\frac{\prod\limits_{k=1+i}^m
q_{2m-i,k}}{\prod\limits_{k=m+1}^{2m-1-i} (c_{2m-i}-c_k)}} \times$ \\
$\left\{ \displaystyle{\frac {\prod\limits_{k=1+i}^m
q_{2m-i,k}}{\prod\limits_{k=m+1}^{2m-1-i} (c_{2m-i}-c_k)}}
-(b_2-b_3)\,c_{2m-i}\,\displaystyle{\frac{\prod\limits_{k=2m+1-i}^{2m}
(c_{2m-i}-c_k)\, \prod\limits_{k=m+1\atop k\ne 2m-i}^{2m}
q_{2m-i,k}}{\prod\limits_{k=1}^m (c_{2m-i}-c_k)\,
\prod\limits_{k=1}^i q_{2m-i,k}}}\right.$\\
$\left. -(b_2-b_3)\displaystyle{\frac{\prod\limits_{k=1+i}^m
q_{2m-i,k}}{\prod\limits_{k=m+1}^{2m-1-i} (c_{2m-i}-c_k)}}
\sum\limits_{r=1}^m
\displaystyle{\frac{c_r}{(c_r-c_{2m-i})\,q_{2m-i,r}}}
\displaystyle{\frac{\prod\limits_{k=m+1\atop k\ne 2m-i}^{2m}
q_{r,k}\, \prod\limits_{k=m+1\atop k\ne 2m-i}^{2m}
(c_r-c_k)}{\prod\limits_{k=1\atop k\ne r}^m q_{r,k}\,
\prod\limits_{k=1\atop k\ne r}^m (c_r-c_k)}} \right\}$.

\item
$\<\e_{1+i},\e_{m+j}\>=
\displaystyle{\frac{p^{31}_{m+1-j}}{q_{m+1-j,2m-i}}}\,
\displaystyle{\frac{\prod\limits_{k=1+i}^m
q_{k,2m-i}}{\prod\limits_{k=m+1}^{2m-1-i}
(c_{2m-i}-c_k)}}\,\displaystyle{\frac{\prod\limits_{k=m+j}^{2m}
q_{m+1-j,k}}{\prod \limits_{k=1}^{m-j}
(c_{m+1-j}-c_k)}}\,\displaystyle{\frac{\prod\limits_{k=m+1\atop
k\ne 2m-i}
^{2m} (c_{m+1-j}-c_k)}{\prod\limits_{k=1\atop k\ne m+1-j}^m q_{m+1-j,k}}}$\\
for~~ $i=1,2,\cdots ,m-1$~~ and~~ $j=1,2,\cdots ,m$

\item
$\<\e_{m+j},\e_{m+i}\>=\delta_{ij}\,
\displaystyle{\frac{p^{31}_{m+1-i}}{p^{32}_{m+1-i}}}\,
\displaystyle{\frac{\prod\limits_{k=m+2-i}^{2m}
(c_{m+1-i}-c_k)}{\prod\limits_{k=1}^{m-i}
(c_{m+1-i}-c_k)}}\,\displaystyle{\frac{\prod\limits_{k=m+i}^{2m}
q_{m+1-i,k}}{\prod\limits_{k=1 \atop k\ne m+1-i}^{m-1+i}
q_{m+1-i,k}}}$~~ for~~ $i,j=1,2,\cdots ,m$

\end{enumerate}
\end{lemma}

\proof
\begin{enumerate}

\item
We obtain by direct computation
$$
\<\e_1,\e_1\>= \sum_{i=1}^{2m} p^{31}_i\,p^{32}_i\,
\frac{\prod\limits_{k=1\atop k\ne i}^{2m}
q_{ik}}{\prod\limits_{k=1\atop k\ne i}^{2m} (c_i-c_k)}.
$$

\noindent For $1\le i\le 2m$ let us set $x_i=c_i+(b_2+b_3)/2$.
Identities (\ref{equation:hg_id_8}), (\ref{equation:hg_id_9}), and
(\ref{equation:hg_id_10}) from the Appendix finish the proof.

\item
A direct computation gives
$$
\<\e_1,\e_{1+i}\>=-p^{31}_{2m-i}\, p^{32}_{2m-i}\,
\frac{\prod\limits_{k=1+i\atop k\ne 2m-i} ^{2m}
q_{2m-i,k}}{\prod\limits_{k=1}^{2m-1-i} (c_{2m-i}-c_k)} -
\sum_{r=1}^m \frac{p^{31}_r\,
p^{32}_r}{c_r-c_{2m-i}}\,\frac{\prod\limits_ {k=m+1\atop k\ne
2m-i}^{2m} q_{rk}\, \prod\limits_{k=1+i}^m
q_{k,2m-i}}{\prod\limits_{k=1\atop k\ne r}^m (c_r-c_k)\,
\prod\limits_{k=m+1}^{2m-1-i} (c_{2m-i}-c_k)}.
$$

\noindent After canceling out common multiples, we have to prove
that
$$
\sum_{r=1}^m \frac{p^{31}_r\, p^{32}_r}{c_r-c_{2m-i}}\,
\frac{\prod\limits_{k=m+1\atop k\ne 2m-i}^{2m}
q_{rk}}{\prod\limits_{k=1\atop k\ne r}^m
(c_r-c_k)}=-p^{31}_{2m-i}\, p^{32}_{2m-i}\,
\frac{\prod\limits_{k=m+1\atop k\ne 2m-i}^{2m}
q_{2m-i,k}}{\prod\limits_{k=1}^m (c_{2m-i}-c_k)} -(b_1-b_3).
$$

\noindent Let us set $x_r=c_r+(b_2+b_3)/2$,
$y_r=c_{m+r}+(b_2+b_3)/2$. Now identities
(\ref{equation:hg_id_5}), (\ref{equation:hg_id_11}), and
(\ref{equation:hg_id_12}) from the Appendix finish the proof.

\item
This formula is proved by direct computation.

\item
A direct computation gives
$$
\<\e_{1+i},\e_{1+j}\>= \frac{\prod\limits_{k=1+i}^m q_{2m-i,k}\,
\prod\limits_{k=1+j}^m q_{2m-j,k}}{\prod\limits_{k=m+1}^{2m-1-i}
(c_{2m-i}-c_k)\, \prod\limits_{k=m+1}^{2m-1-j} (c_{2m-j}-c_k)}\,
\sum_{r=1}^m p^{31}_r\,p^{32}_r
\frac{\prod\limits_{k=m+1\atop{k\ne 2m-i\atop k\ne 2m-j}}^{2m}
q_{r,k}\, \prod\limits_{k=m+1\atop{k\ne 2m-i\atop k\ne 2m-j}}^{2m}
(c_r-c_k)} {\prod\limits_{k=1\atop k\ne r}^m q_{r,k}\,
\prod\limits_{k=1\atop k\ne r}^m (c_r-c_k)}.
$$

\noindent Let us set $x_1=c_1+(b_2+b_3)/2$,
$x_2=c_2+(b_2+b_3)/2,\cdots$, $x_m=c_m+(b_2+b_3)/2$;
$y_1=c_{m+1}+(b_2+b_3)/2$, $y_2=c_{m+2} + (b_2+b_3)/2,\cdots$,
$y_m=c_{m+m}+(b_2+b_3)/2$. Identities (\ref{equation:hg_id_1}) and
(\ref{equation:hg_id_2}) from the Appendix finish the proof.

\item
A direct computation gives
$$
\begin{array}{c}
\<\e_{1+i},\e_{1+i}\> = \frac{\prod\limits_{k=1+i}^m q_{2m-i,k}}
{\prod\limits_{k=m+1}^{2m-1-i} (c_{2m-i}-c_k)} \left\{
\frac{\prod\limits_{k=1+i}^m q_{2m-i,k}}
{\prod\limits_{k=m+1}^{2m-1-i} (c_{2m-i}-c_k)}
\sum\limits_{r=1}^m \frac{p^{31}_r\, p^{32}_r}{(c_r-c_{2m-i})\, q_{2m-i,r}}\, \times \right. \\
\left. \frac{\prod\limits_{k=m+1\atop k\ne 2m-i}^{2m} q_{r,k}\,
\prod\limits_{k=m+1\atop k\ne 2m-i}^{2m}
(c_r-c_k)}{\prod\limits_{k=1\atop k\ne r}^m q_{r,k}\,
\prod\limits_{k=1\atop k\ne r}^m (c_r-c_k)} + p^{31}_{2m-i}\,
p^{32}_{2m-i}\, \frac{\prod\limits_{k=2m+1-i}^{2m}
(c_{2m-i}-c_k)}{\prod\limits_{k=1}^m (c_{2m-i}-c_k)}\,
\frac{\prod\limits_{k=m+1\atop k\ne 2m-i}^{2m}
q_{2m-i,k}}{\prod\limits_{k=1}^i q_{2m-i,k}} \right\}.
\end{array}
$$

\noindent Let us set $x_r=c_r+(b_2+b_3)/2$,
$y_r=c_{m+r}+(b_2+b_3)/2$ for $r=1,2,\cdots ,m$. Identities
(\ref{equation:hg_id_2}) and (\ref{equation:hg_id_5}) from Apendix
finish the proof.

\item Proved by direct computation.

\item Proved by direct computation.

\end{enumerate}
\endproof

 Now it is time to prove Theorem \ref{thm:even:form}, that is prove that the matrices
$\A$, $\B$, and $\C$ are self-adjoint with respect to the scalar
product (\ref{equation:even:form}). $\C$ is self-adjoint with
respect to the scalar product by construction. The space $V$
splits into the direct sum $V=V_{b_1}\oplus V_{b_2}\oplus V_{b_3}$
of the spectral subspaces of $\B$. If the subspaces $V_{b_1}$,
$V_{b_2}$, and $V_{b_3}$ are mutually orthogonal with respect to
the scalar product (\ref{equation:even:form}), then $\B$ is
self-adjoint with respect to it as well. Then $\A$ is also
self-adjoint, as $\A=\B+\C$. Proof of the following lemma finishes
the proof of Theorem \ref{thm:even:form}.

\begin{lemma}
\label{lemma:even:B-orthogonality} The subspaces $V_{b_1}$,
$V_{b_2}$, and $V_{b_3}$ are mutually orthogonal with respect to
the scalar product (\ref{equation:even:form}).
\end{lemma}

\proof We use the formulas of Lemma \ref{lemma:even:evec_of_B} to
express the eigenvectors $\w_i$ of the matrix $\B$ in terms of the
standard basis $\{\e_1,\e_2,\cdots ,\e_{2m}\}$. Then we use the
formulas of Lemma \ref{lemma:even:form_in_good_basis} to expand
$\<\e_i,\e_j\>$.

\begin{enumerate}
\item
$$
\begin{array}{c}
\<\w_1,\w_{1+i}\> = \left\<\e_1,\e_{1+i}-\frac{1}{b_1-b_2}\,
\frac{\prod\limits_{k=1+i}^m
q_{k,2m-i}}{\prod\limits_{k=m+1}^{2m-1-i}
(c_{2m-i}-c_k)}\,\e_1 \right\> = \\
\<\e_1,\e_{1+i}\> - \frac{1}{b_1-b_2}\,
\frac{\prod\limits_{k=1+i}^m
q_{k,2m-i}}{\prod\limits_{k=m+1}^{2m-1-i}
(c_{2m-i}-c_k)}\,\<\e_1,\e_1\> = \\
(b_1-b_3)\,\frac{\prod\limits_{k=1+i}^m
q_{k,2m-i}}{\prod\limits_{k=m+1}^{2m-1-i} (c_{2m-i}-c_k)} -
\frac{(b_1-b_2)(b_1-b_3)}{b_1-b_2}\, \frac{\prod\limits_{k=1+i}^m
q_{k,2m-i}} {\prod\limits_{k=m+1}^{2m-1-i} (c_{2m-i}-c_k)}=0
\end{array}
$$

\item
The identity $\<\w_1,\w_{m+i}\>=0$ $(i=1,2,\cdots ,m)$ reduces to
identity (\ref{equation:hg_id_3}) from the Appendix.

\item
Lemma \ref{lemma:even:evec_of_B} gives us
$$
\w_{1+i}=\e_{1+i}-\frac{1}{b_1-b_2}\,\frac{\prod\limits_{k=1+i}^m
q_{k,2m-i}}{\prod\limits_{k=m+1}^{2m-1-i} (c_{2m-i}-c_k)}\,\e_1.
$$
\noindent We know that $\e_1=\w_1$ and that $\<\w_1,\w_{m+i}\>=0$.
Thus, $\<\w_{1+i},\w_{m+j}\>=\<\e_{1+i},\w_{m+j}\>$. The identity
$\<\e_{1+i},\w_{m+j}\>=0$ after expansion and some simplification
becomes the following identity:
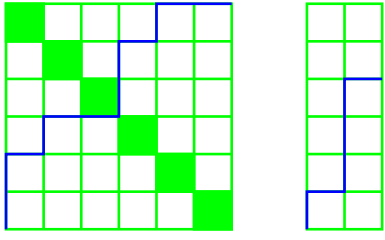
$$
\begin{array}{c}
\sum\limits_{r=1}^{m-1} \frac{\prod\limits_{k=1\atop k\ne m+1-j}^m
q_{2m-r,k}} {\prod\limits_{k=m+1\atop k\ne 2m-r}^{2m-1}
(c_{2m-r}-c_k)}\, \sum\limits_{s=1}^m
\frac{c_s}{q_{2m-i,s}(c_s-c_{2m-i})}\,
\frac{\prod\limits_{k=m+1\atop k\ne 2m-r}^{2m} q_{s,k}\,
\prod\limits_{k=m+1\atop k\ne 2m-r}^{2m} (c_s-c_k)}
{\prod\limits_{k=1\atop k\ne s}^m q_{s,k}\,
\prod\limits_{k=1\atop k\ne s}^m (c_s-c_k)} = \\
1-\frac{c_{2m-i}(c_{2m-i}-c_{2m})}{q_{2m-i,m+1-j}}\,
\frac{\prod\limits_{k=m+1\atop k\ne 2m-i}^{2m} q_{2m-i,k}}
{\prod\limits_{k=1}^m (c_{2m-i}-c_k)} -
\frac{q_{m+1-j,2m}}{q_{m+1-j,2m-i}}\,
\frac{\prod\limits_{k=m+1\atop k\ne 2m-i}^{2m} (c_{m+1-j}-c_k)}
{\prod\limits_{k=1\atop k\ne m+1-j}^m q_{m+1-j,k}}.
\end{array}
$$

\noindent For $i=1,2,\cdots ,m$, let us set $x_i=c_i+(b_2+b_3)/2$
and $y_i=c_{m+i}+(b_2+b_3)/2$. Let us write
$c_s=x_s-y_{2m-i}+y_{2m-i}-(b_2+b_3)/2$ and
$c_{2m-i}=y_{2m-i}-(b_2+b_3)/2$. Now identities
(\ref{equation:hg_id_1}), (\ref{equation:hg_id_5}), and
(\ref{equation:hg_id_13}) of the Appendix finish the proof.
\end{enumerate}
\endproof

 Let us prove Theorem \ref{thm:even:sign_of_form}, that is determine
the inequalities on the real spectra of $\A$, $\B$, and $\C$ which
make the form (\ref{equation:even:form}) sign-definite. It is an
assumption of Theorem \ref{thm:even:sign_of_form} that
$c_1>c_2>\cdots >c_{2m}$. The assumption $a_1>a_2$ of the Theorem
is satisfied automatically, because in our normalized version
$a_1=1$ and $a_2=-1$. \\

\noindent {\it Proof of Theorem \ref{thm:even:sign_of_form} ---}
It is immediately clear from Theorem \ref{thm:even:form} that
$$
\mbox{sign}\(\<\v_i,\v_i\>\)=(-1)^{i-1}\, \mbox{sign}(p^{31}_i\,
p^{32}_i)\, \mbox{sign}\(\prod\limits_{j=1\atop j\ne i}^{2m}
q_{i,j}\).
$$

\noindent Let $Q$ be a $2m\times 2m$ array such that
$Q_{i,j}=q_{i,j}$ for $i\ne j$ and $Q_{i,i}$ are not defined.
Then $Q_{i,j}>Q_{i,j+1}$ and $Q_{i,j}>Q_{i+1,j}$ for all $i$ and
$j$ such that neither of the array elements involved belongs to
the main diagonal. Let $P$ be a $2m\times 2$ matrix such that
$P_{i,1}=p^{32}_i$ and $P_{j,2}=p^{31}_j$. The fact $a_1>a_2$
implies $p^{32}_i>p^{31}_i$. Then $\mbox{sign}(\<\v_i,\v_i\>)=
(-1)^{ i-1+ \# \{j: Q_{i,j}<0 \} + \# \{j: P_{i,j}<0 \}} $.
In order to keep $\mbox{sign}(\<\v_i,\v_i\>)$ constant, the
number of negative elements in the $i$ row of the arrays $Q$ and $P$
must differ from the number of negative elements in the $i+1$ row
by an odd number. This and the fact that $Q_{i,j}=Q_{j,i}$ leaves
room for the following six configurations.
The first is given by the inequalities:

$$
\begin{array}{ll}
p^{31}_{m-1}>0>p^{31}_m     & q_{1,2m-2}>0>q_{1,2m-1} \\
                            & q_{2,2m-3}>0>q_{2,2m-2} \\
p^{32}_{2m-1}>0>p^{32}_{2m} & q_{3,2m-4}>0>q_{3,2m-3} \\
                            & \vdots                  \\
                            & q_{m-1,m}>0>q_{m-1,m+1}
\end{array}.
$$

We have $q_{i,2m-i}<0$ for $i=1,2,\cdots m-1$. Let us sum up these
inequalities with $p^{32}_{2m}<0$ and $p^{31}_m<0$. Recalling that
$\C$ is traceless, we obtain $(m-1)b_2+(m+1)b_3<0$. Recalling that
$\B$ is traceless, we obtain $b_1>b_3$. We have $p^{31}_{m-1}>0$
and $p^{32}_{2m-1}>0$. Thus, $-p^{31}_{m-1}-p^{32}_{2m-1}<0$. We
also have $q_{m-1,2m-1}<q_{m-1,m+1}<0$ for $m>2$ and
$q_{m-1,2m-1}=q_{m-1,m+1}<0$ for $m=2$ because $c_i<c_j$ for
$i>j$. Thus, we have $-p^{31}_{m-1}-p^{32}_{2m-1}+q_{m-1,2m-1}<0$.
This gives us $b_3>b_2$. So, we have  $b_1>b_3>b_2$. Here is a
picture illustrating the case of $m=3$. The line separates
positive elements from negative.

\begin{center}
\pspicture[](0,0)(3,3) \label{pic:em:ineq_tab1}
\psset{unit=0.5cm,linewidth=1pt,linecolor=green} \rput(-3,-1.2){
\pspolygon(1,1)(7,1)(7,7)(1,7) \psline(2,1)(2,7) \psline(3,1)(3,7)
\psline(4,1)(4,7) \psline(5,1)(5,7) \psline(6,1)(6,7)
\psline(1,2)(7,2) \psline(1,3)(7,3) \psline(1,4)(7,4)
\psline(1,5)(7,5) \psline(1,6)(7,6)
\pspolygon[fillstyle=solid,fillcolor=green](1,6)(2,6)(2,7)(1,7)
\pspolygon[fillstyle=solid,fillcolor=green](2,5)(3,5)(3,6)(2,6)
\pspolygon[fillstyle=solid,fillcolor=green](3,4)(4,4)(4,5)(3,5)
\pspolygon[fillstyle=solid,fillcolor=green](4,3)(5,3)(5,4)(4,4)
\pspolygon[fillstyle=solid,fillcolor=green](5,2)(6,2)(6,3)(5,3)
\pspolygon[fillstyle=solid,fillcolor=green](6,1)(7,1)(7,2)(6,2)

\pspolygon(9,1)(11,1)(11,7)(9,7) \psline(10,1)(10,7)
\psline(9,2)(11,2) \psline(9,3)(11,3) \psline(9,4)(11,4)
\psline(9,5)(11,5) \psline(9,6)(11,6)

\psset{linecolor=blue}
\psline(1,1)(1,3)(2,3)(2,4)(4,4)(4,6)(5,6)(5,7)(7,7)
\psline(9,1)(9,2)(10,2)(10,5)(11,5) }
\endpspicture
\end{center}

 The second configuration is given by the following inequalities.
$$
\begin{array}{ll}
p^{31}_1>0>p^{31}_2         & q_{2,2m}>0>q_{3,2m}     \\
                            & q_{3,2m-1}>0>q_{4,2m-1} \\
p^{32}_{m+1}>0>p^{32}_{m+2} & q_{4,2m-2}>0>q_{5,2m-2} \\
                            & \vdots                  \\
                            & q_{m,m+2}>0>q_{m+1,m+2}
\end{array}
$$
We have $p^{31}_2<0$ and $p^{32}_{m+2}<0$. Thus, we have
$-p^{31}_2-p^{32}_{m+2}>0$. We also have $q_{2,m+2}>q_{m,m+2}>0$
for $m>2$ and $q_{2,m+2}=q_{m,m+2}>0$ for $m=2$. Thus,
$q_{2,m+2}-p^{31}_2-p^{32}_{m+2}>0$. This implies $b_2>b_3$. We
have $q_{i,2m+2-i}>0$ for $i=2,3,\cdots ,m$. Summing up these
inequalities with $p^{31}_1>0$ and $p^{32}_{m+1}>0$, we obtain
$(m-1)b_2+(m+1)b_3>0$. Thus, $b_3>b_1$. So, we have
$b_2>b_3>b_1$. Here is the picture illustrating the case of $m=3$.

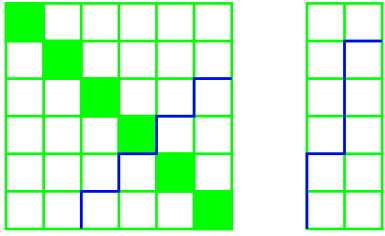
\begin{center}
\pspicture[](0,0)(3,3) \label{pic:em:ineq_tab2}
\psset{unit=0.5cm,linewidth=1pt,linecolor=green} \rput(-3,-1){
\pspolygon(1,1)(7,1)(7,7)(1,7) \psline(2,1)(2,7) \psline(3,1)(3,7)
\psline(4,1)(4,7) \psline(5,1)(5,7) \psline(6,1)(6,7)
\psline(1,2)(7,2) \psline(1,3)(7,3) \psline(1,4)(7,4)
\psline(1,5)(7,5) \psline(1,6)(7,6)
\pspolygon[fillstyle=solid,fillcolor=green](1,6)(2,6)(2,7)(1,7)
\pspolygon[fillstyle=solid,fillcolor=green](2,5)(3,5)(3,6)(2,6)
\pspolygon[fillstyle=solid,fillcolor=green](3,4)(4,4)(4,5)(3,5)
\pspolygon[fillstyle=solid,fillcolor=green](4,3)(5,3)(5,4)(4,4)
\pspolygon[fillstyle=solid,fillcolor=green](5,2)(6,2)(6,3)(5,3)
\pspolygon[fillstyle=solid,fillcolor=green](6,1)(7,1)(7,2)(6,2)

\pspolygon(9,1)(11,1)(11,7)(9,7) \psline(10,1)(10,7)
\psline(9,2)(11,2) \psline(9,3)(11,3) \psline(9,4)(11,4)
\psline(9,5)(11,5) \psline(9,6)(11,6)

\psset{linecolor=blue}
\psline(3,1)(3,2)(4,2)(4,3)(5,3)(5,4)(6,4)(6,5)(7,5)
\psline(9,1)(9,3)(10,3)(10,6)(11,6) }
\endpspicture
\end{center}

 The third configuration is given by the following inequalities.
$$
\begin{array}{ll}
\phantom{p^{32}_m>}0>p^{31}_1 & q_{1,2m-1}>0>q_{1,2m}     \\
                              & q_{2,2m-2}>0>q_{2,2m-1}   \\
p^{32}_m>0>p^{32}_{m+1}       & \vdots                    \\
                              & q_{m-1,m+1}>0>q_{m-1,m+2} \\
                              & \phantom{q_{m-1,m+1}>}0>q_{m,m+1}
\end{array}
$$
The inequalities $p^{32}_{m+1}<0$ and $p^{31}_1<0$ imply the
inequality $2b_3+c_1+c_{m+1}<0$. The inequality $q_{1,2m-1}>0$
implies the inequality $q_{1,m+1}>0$ because $c_{m+1}>c_{2m-1}$
for $m>2$ and $c_{m+1}=c_{2m-1}$ for $m=2$. Now, the inequalities
$b_2+b_3+c_1+c_{m+1}>0$ and $-2b_3-c_1-c_{m+1}>0$ imply the
inequality $b_2-b_3>0$. So, $b_2>b_3$. Let us sum up the
inequalities $q_{i,2m+1-i}<0$ for $i=1,2,\cdots m$. The sum of all
the $c_i$ is equal to zero. Thus we obtain $mb_2+mb_3<0$ which is
equivalent to $b_2-b_1<0$. This gives us $b_1>b_2>b_3$. Here is
the picture illustrating the case of $m=3$.

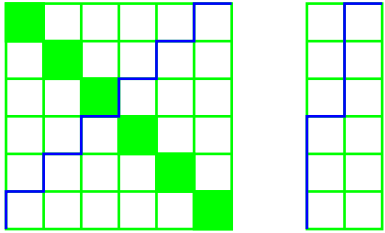
\begin{center}
\pspicture[](0,0)(3,3) \label{pic:em:ineq_tab3}
\psset{unit=0.5cm,linewidth=1pt,linecolor=green} \rput(-3,-1){
\pspolygon(1,1)(7,1)(7,7)(1,7) \psline(2,1)(2,7) \psline(3,1)(3,7)
\psline(4,1)(4,7) \psline(5,1)(5,7) \psline(6,1)(6,7)
\psline(1,2)(7,2) \psline(1,3)(7,3) \psline(1,4)(7,4)
\psline(1,5)(7,5) \psline(1,6)(7,6)
\pspolygon[fillstyle=solid,fillcolor=green](1,6)(2,6)(2,7)(1,7)
\pspolygon[fillstyle=solid,fillcolor=green](2,5)(3,5)(3,6)(2,6)
\pspolygon[fillstyle=solid,fillcolor=green](3,4)(4,4)(4,5)(3,5)
\pspolygon[fillstyle=solid,fillcolor=green](4,3)(5,3)(5,4)(4,4)
\pspolygon[fillstyle=solid,fillcolor=green](5,2)(6,2)(6,3)(5,3)
\pspolygon[fillstyle=solid,fillcolor=green](6,1)(7,1)(7,2)(6,2)

\pspolygon(9,1)(11,1)(11,7)(9,7) \psline(10,1)(10,7)
\psline(9,2)(11,2) \psline(9,3)(11,3) \psline(9,4)(11,4)
\psline(9,5)(11,5) \psline(9,6)(11,6)

\psset{linecolor=blue}
\psline(1,1)(1,2)(2,2)(2,3)(3,3)(3,4)(4,4)(4,5)(5,5)(5,6)(6,6)(6,7)(7,7)
\psline(9,1)(9,4)(10,4)(10,7)(11,7) }
\endpspicture
\end{center}

 The fourth configuration is given by the following inequalities.
$$
\begin{array}{ll}
p^{31}_m>0>p^{31}_{m+1} & q_{1,2m}>0>q_{2,2m}     \\
                        & q_{2,2m-1}>0>q_{3,2m-1} \\
p^{32}_{2m}>0           & \vdots                  \\
                        & q_{m-1,m+2}>0>q_{m,m+2} \\
                        & q_{m,m+1}>0
\end{array}.
$$

We have $q_{m,2m}<q_{m,m+2}<0$ for $m>2$ and
$q_{m,2m}=q_{m,m+2}<0$ for $m=2$. We also have $-p^{32}_{2m}<0$
and $-p^{31}_m<0$. Summing up these inequalities, we obtain
$b_3>b_2$. We also have $q_{i,2m+1-i}>0$ for $i=1,2,\cdots ,m$.
Summing up these inequalities, we obtain $b_2>b_1$. This gives us
$b_3>b_2>b_1$. Here is the picture illustrating the case of $m=3$.

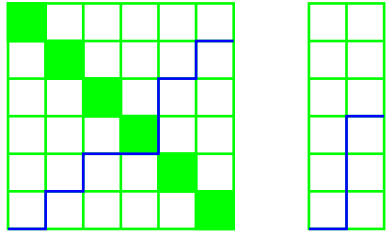
\begin{center}
\pspicture[](0,0)(3,3) \label{pic:em:ineq_tab4}
\psset{unit=0.5cm,linewidth=1pt,linecolor=green} \rput(-3,-1){
\pspolygon(1,1)(7,1)(7,7)(1,7) \psline(2,1)(2,7) \psline(3,1)(3,7)
\psline(4,1)(4,7) \psline(5,1)(5,7) \psline(6,1)(6,7)
\psline(1,2)(7,2) \psline(1,3)(7,3) \psline(1,4)(7,4)
\psline(1,5)(7,5) \psline(1,6)(7,6)
\pspolygon[fillstyle=solid,fillcolor=green](1,6)(2,6)(2,7)(1,7)
\pspolygon[fillstyle=solid,fillcolor=green](2,5)(3,5)(3,6)(2,6)
\pspolygon[fillstyle=solid,fillcolor=green](3,4)(4,4)(4,5)(3,5)
\pspolygon[fillstyle=solid,fillcolor=green](4,3)(5,3)(5,4)(4,4)
\pspolygon[fillstyle=solid,fillcolor=green](5,2)(6,2)(6,3)(5,3)
\pspolygon[fillstyle=solid,fillcolor=green](6,1)(7,1)(7,2)(6,2)

\pspolygon(9,1)(11,1)(11,7)(9,7) \psline(10,1)(10,7)
\psline(9,2)(11,2) \psline(9,3)(11,3) \psline(9,4)(11,4)
\psline(9,5)(11,5) \psline(9,6)(11,6)

\psset{linecolor=blue}
\psline(1,1)(2,1)(2,2)(3,2)(3,3)(5,3)(5,5)(6,5)(6,6)(7,6)
\psline(9,1)(10,1)(10,4)(11,4) }
\endpspicture
\end{center}

 The fifth configuration is given by the following inequalities.

$$
\begin{array}{ll}
\phantom{p^{32}_m>}0>p^{31}_1 & q_{1,2m}>0>q_{2,2m}     \\
                              & q_{2,2m-1}>0>q_{3,2m-1} \\
p^{32}_m>0>p^{32}_{m+1}       & \vdots                  \\
                              & q_{m-1,m+2}>0>q_{m,m+2} \\
                              & q_{m,m+1}>0
\end{array}
$$

We have $q_{i,2m+2-i}<0$ for $i=2,3,\cdots ,m$. Summing up these
inequalities with $p^{31}_1<0$ and $p^{32}_{m+1}<0$, we obtain
$(m-1)b_2+(m+1)b_3<0$. The last is equivalent to $b_1>b_3$. We
also have $q_{i,2m+1-i}>0$ for $i=1,2,\cdots ,m$. Summing up these
inequalities, we obtain $mb_2+mb_3>0$ which is equivalent to
$b_2>b_1$. This gives us $b_2>b_1>b_3$. Here is the picture
illustrating the case of $m=3$.

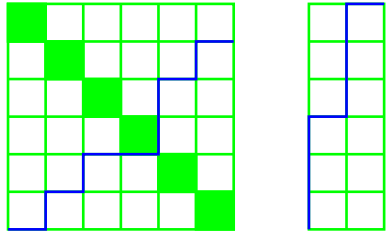
\begin{center}
\pspicture[](0,0)(3,3) \label{pic:em:ineq_tab5}
\psset{unit=0.5cm,linewidth=1pt,linecolor=green} \rput(-3,-1.2){
\pspolygon(1,1)(7,1)(7,7)(1,7) \psline(2,1)(2,7) \psline(3,1)(3,7)
\psline(4,1)(4,7) \psline(5,1)(5,7) \psline(6,1)(6,7)
\psline(1,2)(7,2) \psline(1,3)(7,3) \psline(1,4)(7,4)
\psline(1,5)(7,5) \psline(1,6)(7,6)
\pspolygon[fillstyle=solid,fillcolor=green](1,6)(2,6)(2,7)(1,7)
\pspolygon[fillstyle=solid,fillcolor=green](2,5)(3,5)(3,6)(2,6)
\pspolygon[fillstyle=solid,fillcolor=green](3,4)(4,4)(4,5)(3,5)
\pspolygon[fillstyle=solid,fillcolor=green](4,3)(5,3)(5,4)(4,4)
\pspolygon[fillstyle=solid,fillcolor=green](5,2)(6,2)(6,3)(5,3)
\pspolygon[fillstyle=solid,fillcolor=green](6,1)(7,1)(7,2)(6,2)

\pspolygon(9,1)(11,1)(11,7)(9,7) \psline(10,1)(10,7)
\psline(9,2)(11,2) \psline(9,3)(11,3) \psline(9,4)(11,4)
\psline(9,5)(11,5) \psline(9,6)(11,6)

\psset{linecolor=blue}
\psline(1,1)(2,1)(2,2)(3,2)(3,3)(5,3)(5,5)(6,5)(6,6)(7,6)
\psline(9,1)(9,4)(10,4)(10,7)(11,7) }
\endpspicture
\end{center}

 The last configuration possible is given by the following inequalities.

$$
\begin{array}{ll}
p^{31}_m>0>p^{31}_{m+1} & q_{1,2m-1}>0>q_{1,2m}     \\
                        & q_{2,2m-2}>0>q_{2,2m-1}   \\
p^{32}_{2m}>0           & \vdots                    \\
                        & q_{m-1,m+1}>0>q_{m-1,m+2} \\
                        & \phantom{q_{m-1,m+1}>}0>q_{m,m+1}
\end{array}
$$

We have $q_{i,2m+1-i}<0$ for $i=1,2,\cdots ,m$. Summing up these
inequalities, we obtain $mb_2+mb_3<0$ which is equivalent to
$b_1>b_2$. We also have $q_{i,2m-i}>0$ for $i-1,2,\cdots ,m-1$.
Summing up these inequalities with $p^{32}_{2m}>0$ and
$p^{31}_m>0$, we obtain $(m-1)b_2+(m+1)b_3>0$ which is equivalent
to $b_3>b_1$. This gives us $b_3>b_1>b_2$. Here is the picture
illustrating the case of $m=3$.

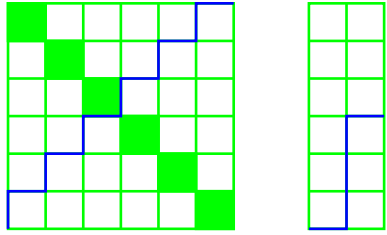
\begin{center}
\pspicture[](0,0)(3,3) \label{pic:em:ineq_tab6}
\psset{unit=0.5cm,linewidth=1pt,linecolor=green} \rput(-3,-1.2){
\pspolygon(1,1)(7,1)(7,7)(1,7) \psline(2,1)(2,7) \psline(3,1)(3,7)
\psline(4,1)(4,7) \psline(5,1)(5,7) \psline(6,1)(6,7)
\psline(1,2)(7,2) \psline(1,3)(7,3) \psline(1,4)(7,4)
\psline(1,5)(7,5) \psline(1,6)(7,6)
\pspolygon[fillstyle=solid,fillcolor=green](1,6)(2,6)(2,7)(1,7)
\pspolygon[fillstyle=solid,fillcolor=green](2,5)(3,5)(3,6)(2,6)
\pspolygon[fillstyle=solid,fillcolor=green](3,4)(4,4)(4,5)(3,5)
\pspolygon[fillstyle=solid,fillcolor=green](4,3)(5,3)(5,4)(4,4)
\pspolygon[fillstyle=solid,fillcolor=green](5,2)(6,2)(6,3)(5,3)
\pspolygon[fillstyle=solid,fillcolor=green](6,1)(7,1)(7,2)(6,2)

\pspolygon(9,1)(11,1)(11,7)(9,7) \psline(10,1)(10,7)
\psline(9,2)(11,2) \psline(9,3)(11,3) \psline(9,4)(11,4)
\psline(9,5)(11,5) \psline(9,6)(11,6)

\psset{linecolor=blue}
\psline(1,1)(1,2)(2,2)(2,3)(3,3)(3,4)(4,4)(4,5)(5,5)(5,6)(6,6)(6,7)(7,7)
\psline(9,1)(10,1)(10,4)(11,4) }
\endpspicture
\end{center}

 So, in these six cases the form $\<*,*\>$ is sign-definite.
Lemma \ref{lemma:even:form_in_good_basis} gives $\<\e_1,\e_1\> =
(b_1-b_2)(b_1-b_3)$. Thus, $\mbox{sign}(\<*,*\>) = \mbox{sign} \(
(b_1-b_2)(b_1-b_3) \)$.
\endproof

\subsection{Odd Family}
\label{subsection:proofs:odd}

For the hypergeometric, odd, and even family, let us call the
objects $\{V;\, \A=\B+\C,\B,\C; \<*,*\>\}$ where $(\A,\B,\C)$ is a
rigid irreducible triple of matrices of the corresponding specrtal
types and $\<*,*\>$ is the non-degenerate scalar product such that
$\A$, $\B$, and $\C$ are self-adjoint with respect to it, $m$-{\it
hypergeometric module}, $m$-{\it even module}, and $m$-{\it odd
module}. Let us denote these objects as $HGM_m$, $EM_m$, and
$OM_m$. The reason for calling these objects modules comes from
the theory of quiver representations and will not be explained
in this paper. \\

 It is possible to prove Theorems \ref{thm:odd:A}, \ref{thm:odd:form}, Lemma
\ref{lemma:odd:eigenbasis_of_C}, etc. in the same fashion as for
the even family. But we choose a different approach. We show that
by means of violating the ``generic eigenvalues'' condition it is
possible to construct $OM_{m-1}$ as a submodule and as a factor
module of $EM_m$. Then all the formulas follow
from the corresponding formulas for the even family. \\

 Let $V$ be the same $2m$-dimensional linear space as in the previous subsection and let
$\e_1,\cdots ,\e_{2m}$ be the standard basis of $V$. Let $\A$,
$\B$, and $\C$ be the matrices from the previous subsection, too.
Fix an integer $i$ such that $1\le i\le m$. Let $V_{\hat{i}}^s$ be
the subspace of $V$ spanned by the vectors $\e_1,\e_2, \cdots
,\e_{2m-i},\widehat{\e_{2m+1-i}},\e_{2m+2-i}, \cdots ,\e_{2m}$. It
follows from the formulas of Lemma
\ref{lemma:even:eigenbasis_of_C} that $V_{\hat{i}}^s$ is spanned
by $\v_1,\v_2,\cdots ,\v_{i-1},\widehat{\v_i},\v_{i+1},\cdots
,\v_{2m}$ (hence the notation). Then the following lemma follows
at once from the formulas for $\A$, $\B$, and $\C$ of Subsection
\ref{subsec:results:even}.

\begin{lemma}
\label{lemma:even:inv_(2m-1)-subsp} If $p^{32}_i=0$, then
$V_{\hat{i}}^s$ is invariant with respect to $\A$, $\B$, and $\C$.
\end{lemma}

 Thus, it makes sense to consider the restrictions of $\A$, $\B$, and $\C$ to
$V_{\hat{i}}^s$ and call them $\A_{\hat{i}}^s$, $\B_{\hat{i}}^s$,
and $\C_{\hat{i}}^s$. We will also call $\<*,*\>_{\hat{i}}^s$ the
form $\<*,*\>$ restricted to $V_{\hat{i}}^s$. Note that
$p^{32}_i=0$ forces $q_{ij}=p^{21}_j$ and $c_i-c_j=-p^{32}_j$.

\begin{theorem}
\label{thm:OM_as_subm_of_EM} If $(\s(\A),\s(\B),\s(\C))$ is a
generic point of the intersection of $S((m,m),(1,m-1,m),(1^{2m}))$
with the hyperplane given by the equation $p^{32}_i=0$ for a fixed
$1\le i\le m$, then $\{V_{\hat{i}}^s;\,\A_{\hat{i}}^s,
\B_{\hat{i}}^s, \C_{\hat{i}}^s; \<*,*\>_{\hat{i}}^s\}$ is
$OM_{m-1}$.
\end{theorem}

 Here is an example of the matrices $\B_{\hat{2}}^s$ and $\C_{\hat{2}}^s$
obtained from the matrices $\B$ and $\C$ of Example
\ref{example:even:m=3:B,C} by setting $p^{32}_2=0$ and restricting
them to $V_{\hat{2}}^s$.

\begin{example}
\label{example:odd_from_even:sub:m=2:B,C}
$$
\B^s_{\hat{2}}=\left[
\begin{array}{c|cc|cc}
b_1 & -\frac{p^{21}_5\, q_{35}}{c_4-c_5} & q_{34} &
-\frac{p^{31}_3\, q_{34}\, q_{35}}{p^{32}_3 (c_1-c_3)} & p^{31}_1 \\
 & & & & \\
\hline
 & & & & \\
0 & b_2 & 0 & -\frac{p^{31}_3\, q_{15}\, q_{34}}{p^{32}_3\, (c_1-c_3)} & p^{31}_1 \\
 & & & & \\
0 & 0 & b_2 & -\frac{p^{31}_3\, p^{21}_4\, q_{14}\,
q_{35}}{p^{32}_3 (c_1-c_3)(c_4-c_5)} &
\frac{p^{31}_1\, p^{21}_4}{c_4-c_5} \\
 & & & & \\
\hline
 & & & & \\
0 & 0 & 0 & b_3 & 0 \\
 & & & & \\
0 & 0 & 0 & 0 & b_3
\end{array}
\right]
$$
\bigskip
$$
\C_{\hat{2}}^s=\left[
\begin{array}{c|cc|cc}
c_6 & 0  & 0 & 0 & 0  \\
 & & & & \\
\hline
 & & & & \\
-q_{15} & c_5 & 0 & 0 & 0 \\
 & & & & \\
-\frac{p^{21}_4\, q_{14}}{c_4-c_5} & 0 & c_4 & 0 & 0 \\
 & & & & \\
\hline
 & & & & \\
-p^{32}_3 & -\frac{p^{32}_3\, p^{21}_5}{c_4-c_5} & p_{33} & c_3 & 0 \\
 & & & & \\
-\frac{q_{14}\, q_{15}}{c_1-c_3} & -\frac{p^{21}_5\, q_{14}\,
q_{35}} {(c_1-c_3)(c_4-c_5)} & \frac{q_{15}\, q_{34}}{c_1-c_3} & 0
& c_1
\end{array}
\right]
$$
\bigskip
\end{example}

 This way  $OM_{m-1}$ is constructed as a submodule of $EM_m$.
It is also possible to construct $OM_{m-1}$ as a factor module of
$EM_m$. Let us fix $i$ such that $1\le i\le m$. The following
lemma follows at once from the formulas for $\A$, $\B$, and $\C$
of Subsection \ref{subsec:results:even}.

\begin{lemma}
\label{lemma:even:inv_1-subsp} If $p^{31}_i=0$, then $\A$, $\B$,
and $\C$ preserve the one-dimensional subspace of $V$ spanned by
the vector $\e_{2m+1-i}$.
\end{lemma}

 Thus, it makes sense to consider $\A$, $\B$, and $\C$ acting on the factor
space $V_{\hat{i}}^f$. Let us call these new operators
$\A_{\hat{i}}^f$, $\B_{\hat{i}}^f$, and $\C_{\hat{i}}^f$. If
$p^{31}_i=0$, then $\<\v_i,\v_i\>=0$. Thus, the form
$\<*,*\>_{\hat{i}}^f$ induced on $V_{\hat{i}}^f$ by the form
$\<*,*\>$ on $V$ is well defined. Note that $p^{31}_i=0$ forces
$q_{ij}=p^{22}_j$ and $c_i-c_j=-p^{31}_j$.

\begin{theorem}
\label{thm:OM_as_factor_of_EM} If $(\s(\A),\s(\B),\s(\C))$ is a
generic point of the intersection of $S((m,m),(1,m-1,m),(1^{2m}))$
with the hyperplane given by the equation $p^{31}_i=0$ for a fixed
$1\le i\le m$, then $\{V_{\hat{i}}^f;\,\A_{\hat{i}}^f,
\B_{\hat{i}}^f, \C_{\hat{i}}^f; \<*,*\>_{\hat{i}}^f\}$ is
$OM_{m-1}$.
\end{theorem}

 Here is an example of the matrices $\B_{\hat{2}}^f$ and $\C_{\hat{2}}^f$
obtained from the matrices $\B$ and $\C$ of Example
\ref{example:even:m=3:B,C} by setting $p^{31}_2=0$ and passing to
the factor space $V_{\hat{2}}^f$.

\begin{example}
\label{example:odd_from_even:m=2:factor:B,C}
$$
\B^f_{\hat{2}}=\left[
\begin{array}{c|cc|cc}
b_1 & -\frac{p^{22}_5\, q_{35}}{c_4-c_5} & q_{34} & -\frac{q_{34}\, q_{35}}{c_1-c_3} & p^{31}_1 \\
 & & & & \\
\hline
 & & & & \\
0 & b_2 & 0 & -\frac{q_{15}\, q_{34}}{c_1-c_3} & p^{31}_1 \\
 & & & & \\
0 & 0 & b_2 & -\frac{p^{22}_4\, q_{14}\,
q_{35}}{(c_1-c_3)(c_4-c_5)} &
\frac{p^{31}_1\, p^{22}_4}{c_4-c_5} \\
 & & & & \\
\hline
 & & & & \\
0 & 0 & 0 & b_3 & 0 \\
 & & & & \\
0 & 0 & 0 & 0 & b_3
\end{array}
\right]
$$
\bigskip
$$
\C_{\hat{2}}^f=\left[
\begin{array}{c|cc|cc}
c_6 & 0  & 0 & 0 & 0  \\
 & & & & \\
\hline
 & & & & \\
-q_{15} & c_5 & 0 & 0 & 0 \\
 & & & & \\
-\frac{p^{22}_4\, q_{14}}{c_4-c_5} & 0 & c_4 & 0 & 0 \\
 & & & & \\
\hline
 & & & & \\
-p^{32}_3 & -\frac{p^{32}_3\, p^{22}_5}{c_4-c_5} & p^{32}_3 & c_3 & 0 \\
 & & & & \\
-\frac{p^{32}_1\, q_{14}\, q_{15}}{p^{31}_1\, (c_1-c_3)} &
-\frac{p^{32}_1\, p^{22}_5\, q_{14}\, q_{35}} {p^{31}_1
(c_1-c_3)(c_4-c_5)} & \frac{p^{32}_1\, q_{15}\, q_{34}}{p^{31}_1
(c_1-c_3)} & 0 & c_1
\end{array}
\right]
$$
\bigskip
\end{example}

 We first prove Theorem \ref{thm:OM_as_subm_of_EM}
and then we derive all the proofs for the odd family from what we
already know about the
even family. Let us prove Theorem \ref{thm:OM_as_subm_of_EM}. \\

\noindent {\it Prove of Theorem \ref{thm:OM_as_subm_of_EM} ---~}
\label{proof1} It is clear that $\B_{\hat{i}}^s$ is
diagonalizable and that $\s(\B_{\hat{i}}^s)=\{b_1,\underbrace{b_2,
\cdots ,b_2}_{m-1~times}, \underbrace{b_3,\cdots
,b_3}_{m-1~times}\}$. It is clear that $\C_{\hat{i}}^s$ is
diagonalizable and that $\s(\C_{\hat{i}}^s)=\{c_1,c_2,\cdots
,c_{i-1},\widehat{c_i},c_{i+1},\cdots ,c_{2m}\}$. In view of Lemma
\ref{lemma:even:inv_(2m-1)-subsp} and Theorem \ref{thm:even:A}, it
is clear that $A_{\hat{i}}^s$ is diagonalizable as well. In the
notations of Subsection \ref{subsection:proofs:even}, $\A$ has
eigenvectors ${\bf a}_j^+$ corresponding to the eigenvalue $a_1$
(normalized to $1$). Vectors ${\bf a}_j^-$ are eigenvectors of
$\A$ corresponding to the eigenvalue $a_2$ (normalized to $-1$).
Once we set $p^{32}_i=0$, all the eigenvectors of $\A$ belong to
$V^s_{\hat i}$ except for ${\bf a}_{m+1-i}^-$ and the proof
follows immediately.
\endproof

 To finish the rest of the proofs for the odd family, we just have to say that all
the formulas for $OM_m$ in this paper were obtained from the
formulas for $EM_{m+1}$ by setting $p^{32}_{m+1}=0$ and
renumbering the remaining $c_1,c_2,\cdots ,c_m,c_{m+2},\cdots
,c_{2m+2}$ as $c_1,c_2,\cdots ,c_{2m+1}$.

\begin{remark}
\label{remark:EM_from_OM} {\rm  In exactly the same fashion, we
can construct $EM_m$ as a factor module of $OM_m$ by setting
either $p^{31}_i=0$ for a fixed $1\le i\le m$ or $p^{21}_i=0$ for
a fixed $m+1\le i\le 2m$. Also similarly, one can show that
setting $b_i+c_{m+1-i}-a_2=0$ for a fixed $1\le i\le m$ creates
$HGM_{m-1}$ as a submodule of $HGM_m$. }
\end{remark}

\subsection{$HGM_m$ from $EM_m$}
\label{subsection:HGM_from_EM}

In this subsection we show that there exist $HGM_m$ and
$HGM_{m+1}$ which naturally
``live inside'' $EM_m$ with roles of $\A$ and $\B$ interchanged.\\

 Let us consider $EM_m$ constructed in Subsection~\ref{subsec:results:even}.
Let $V_1$ be the subspace of $V$ spanned by $\e_1$. Let $V_2$ be
the subspace of $V$ spanned by the vectors $\e_2$, $\e_3,\cdots
,\e_m$. Finally, let $V_3$ be the subspace of $V$ spanned by the
vectors $\e_{m+1}$, $\e_{m+2},\cdots ,\e_{2m}$. Let ${\bf
I}':V_1\oplus V_2\hookrightarrow V$ be an inclusion. Let ${\bf
P}': V \rightarrow V_1\oplus V_2$ be the projection along $V_3$.
Set $\A'={\bf P}'\, (-\B)\, {\bf I}'$, $\B'={\bf P}'\, (-\A)\,
{\bf I}'$, and $\C'={\bf P}'\, \C\, {\bf I}'$.

\begin{theorem}
\label{thm:HGM_m_in_EM_m} If $(\s(\A),\s(\B),\s(\C))$ is a generic
point of $S''((m,m),(1,m-1,m),(1^{2m}))$, then $\{V_1\oplus V_2;\,
\A',\B',\C'\}$ is $HGM_m$.
\end{theorem}

 Let ${\bf I}'':V_1\oplus V_3\hookrightarrow V$ be an inclusion. Let
${\bf P}'': V \rightarrow V_1\oplus V_3$ be the projection along
$V_2$. Set $\A''={\bf P}''\, (-\B)\, {\bf I}''$, $\B''={\bf P}''\,
(-\A)\, {\bf I}''$, and $\C''={\bf P}''\, \C\, {\bf I}''$.

\begin{theorem}
\label{thm:HGM_m+1_in_EM_m} If $(\s(\A),\s(\B),\s(\C))$ is a
generic point of $S''((m,m),(1,m-1,m),(1^{2m}))$, then
$\{V_1\oplus V_3;\, \A'',\B'',\C''\}$ is $HGM_{m+1}$.
\end{theorem}

 The following lemma is a simple exercise in linear algebra.
\begin{lemma}
\label{lemma:det}
$$
\left|
\begin{array}{ccccc}
\alpha_1     & \beta_1  & \beta_2  & \ldots & \beta_{m-1}\\
\gamma_1     & \alpha_2 & 0        & \ldots & 0          \\
\gamma_2     & 0        & \alpha_3 & \ldots & 0          \\
\vdots       & \vdots   & \vdots   & \ddots & \vdots     \\
\gamma_{m-1} & 0        & 0        & \ldots & \alpha_m
\end{array}
\right|=\prod\limits_{i=1}^m \alpha_i-\sum\limits_{j=1}^{m-1}
\beta_j \gamma_j \prod\limits_{k=2\atop k\ne 1+i}^m \alpha_k.
$$
\end{lemma}
\bigskip

\noindent {\it Proof of Theorem \ref{thm:HGM_m_in_EM_m} ---}
It is clear that $\A'=\B'+\C'$. It is clear that $\A'$ is
diagonalizable and that $\s(\A')=(-b_1,\underbrace{-b_2,\cdots
,-b_2}_ {m-1~times})$. It is also clear that $\C'$ is
diagonalizable and that $\s(\C_m)=(c_{m+1},c_{m+2},\cdots
,c_{2m})$. Let us prove that $\B'$ is diagonalizable and
$\s(\B')=(b_3+c_1-a_1-a_2,b_3+c_2-a_1-a_2,\cdots
,b_3+c_m-a_1-a_2)$. Let us consider the matrix $\B'$. Application
of Lemma \ref{lemma:det} immediately yields the following formula
for the characteristic polynomial of $\B'$:

\begin{equation}
\label{equation:char_pol_of_B'}
\begin{array}{l}
|\B'-x\,\Id|=(-1)^m\,\[(b_1+c_{2m}+x)\,\prod\limits_{k=1}^{m-1} (b_2+c_{m+k}+x) +\right. \\
\left. \sum\limits_{i=1}^{m-1}
\displaystyle{\frac{\prod\limits_{k=1}^m q_{k,2m-i}}
{\prod\limits_{k=m+1\atop k\ne 2m-i}^{2m-1} (c_{2m-i}-c_k)}}\,
\prod\limits_{k=2\atop k\ne 1+i}^m (b_2+c_{2m+1-k}+x)\].
\end{array}
\end{equation}

\noindent where $p_i^{jk}$ and $q_{i,j}$ are as in (\ref{eq:p,q}). \\

 Let us prove that $x_j=b_3+c_j-a_1-a_2$ where $j=1,2,\cdots ,m$ are the roots of
the polynomial. For that, we have to prove the following identity:

$$
\sum_{i=1}^{m-1} \frac{\prod\limits_{k=1\atop k\ne j}^m
q_{k,2m-i}} {\prod\limits_{k=m+1\atop k\ne 2m-i}^{2m-1}
(c_{2m-i}-c_k)}=a_1+a_2-b_1-b_3-c_j-c_{2m}.
$$

\noindent Identity (\ref{equation:hg_id_3}) from the Appendix reduces
this identity to the trace identity. The matrix $\B'$ is
diagonalizable because it has all eigenvalues distinct. If
$(\s(\A),\s(\B),\s(\C))$ is a generic point of
$S''((m,m),(1,m-1,m),(1^{2m}))$, then $(\s(\A'),\s(\B'),\s(\C'))$
is a generic point of $S''((1,m-1),(1^m),(1^m))$. Then it follows
from Theorems \ref{thm:form} and \ref{thm:irr} and Corollary
\ref{corollary:uniqueness_of_form} that there exists the unique
non-trivial non-degenerate symmetric bilinear form $\<*,*\>'$ such
that $\A'$, $\B'$, and $\C'$ are self-adjoint with respect to it.
\endproof

%
%
%

\noindent {\it Prove of Theorem \ref{thm:HGM_m+1_in_EM_m} ---}
It is clear that $\A''$ is diagonalizable and that
$\s(\A'')=(-b_1,\underbrace{-b_3\cdots ,-b_3}_{m~times})$. It is
also clear that $\C''$ is diagonalizable and that
$\s(\C'')=(c_1,c_2,\cdots ,c_m,c_{2m})$. Let us prove that $\B''$
is diagonalizable and that
$\s(B'')=(-a_1,-a_2,b_2+c_{m+1}-a_1-a_2,b_2+c_{m+2}-a_1-a_2,\cdots
,b_2+c_{2m-1}-a_1-a_2)$. Let us consider the matrix $\B''$.
Application of Lemma \ref{lemma:det} immediately yields the
following formula for the characteristic polynomial of $\B''$:

\begin{equation}
\label{equation:char_pol_of_B''}
\begin{array}{l}
|\B''-x\,\Id|=(-1)^{m+1}\, \[(b_1+c_{2m}+x)\,\prod\limits_{k=1}^m (b_3+c_{m+1-k}+x) + \right. \\
\left. \sum\limits_{i=1}^m p^{31}_{m+1-i}\, p^{32}_{m+1-i}\,
\displaystyle{\frac{\prod\limits_{k=m+1}^{2m-1}
q_{k,m+1-i}}{\prod\limits_{k=1\atop k\ne m+1-i}^m
(c_{m+1-i}-c_k)}}\, \prod\limits_{k=1\atop k\ne i}^m (b_3+c_{m+1-k}+x) \] = \\
(-1)^{m+1}\, \[ (b_1+c_{2m}+x)\,\prod\limits_{k=1}^m (b_3+c_k+x) +
\sum\limits_{i=1}^m p^{31}_i \,p^{32}_i\,
\displaystyle{\frac{\prod\limits_{k=m+1}^{2m-1}
q_{k,i}}{\prod\limits_{k=1\atop k\ne i}^m (c_i-c_k)}}\,
\prod\limits_{k=1\atop k\ne i}^m (b_3+c_k+x) \].
\end{array}
\end{equation}

 Let us prove that $-a_1$ is a root of this polynomial. For that, we have to show
$$
\sum_{i=1}^m p^{32}_i\, \frac{\prod\limits_{k=m+1}^{2m-1} q_{i,k}}
{\prod\limits_{k=1\atop k\ne i}^m (c_i-c_k)}=a_1-b_1-c_{2m}.
$$

\noindent Let us set $x_1=p^{32}_1$, $x_2=p^{32}_2,\cdots
,x_m=p^{32}_m$; $y_1=-p^{31}_{m+1}$, $y_2=-p^{31}_{m+2},\cdots
,y_{m-1}=-p^{31}_{2m-1}$. Then $q_{ik}=x_i-y_k$,
$c_i-c_k=x_i-x_k$, and  the formula to prove becomes equivalent to
identity (\ref{equation:hg_id_3}) from the Appendix. The proof that
$-a_2$ is a root
of the polynomial is exactly the same. \\

 Let us fix $1\le j\le m-1$. The same identity (\ref{equation:hg_id_3}) proves that
$b_2+c_{m+j}-a_1-a_2$ is a root of the polynomial. The rest of the
proof of this theorem is the same as the corresponding part of the
proof of Theorem \ref{thm:HGM_m_in_EM_m}.
\endproof

\begin{remark}
\label{remark:HGM_inside_OM} {\rm Similarly, one can construct
$HGM_{m+1}$ ``inside'' $OM_m$. }
\end{remark}

\subsection{Extra case of Simpson}
\label{subsection:E_8_hat}

Consider the following vectors.

\begin{equation}
\label{equation:E_8_hat:evB}
\begin{array}{lll}

\w_1 & = & \( 1,0,0,0,0,0 \), \\
 & & \\
\w_2 & = & \( 0,1,0,0,0,0 \), \\
 & & \\
\w_3 & = & \( \frac{p_{16}\, q_{245}}{(b_1-b_2)(c_3-c_4)},
-\frac{p_{15}\, q_{235}\, q_{246}}{(b_1-b_2)(c_3-c_4)(c_5-c_6)},1,0,0,0 \), \\
 & & \\
\w_4 & = & \( -\frac{p_{16}}{b_1-b_2},\frac{p_{15}\,
q_{236}}{(b_1-b_2)(c_5-c_6)},
0,1,0,0 \), \\
 & & \\
\w_5 & = & \( -\frac{p_{16}\, q_{245}\,
(q_{126}-p_{31})}{(b_1-b_3)(b_2-b_3)(c_1-c_2)}, \frac{p_{15}\,
q_{236}\, q_{246}\,
(q_{125}-p_{31})}{(b_1-b_3)(b_2-b_3)(c_1-c_2)(c_5-c_6)},
\frac{p_{24}\, q_{236}}{(b_2-b_3)(c_1-c_2)}, \right. \\
 & & \left. \frac{p_{23}\, q_{245}\, q_{246}}{(b_2-b_3)(c_1-c_2)(c_3-c_4)},1,0 \), \\
 & & \\
\w_6 & = & \( -\frac{p_{16}\,
(q_{126}-p_{32})}{(b_1-b_3)(b_2-b_3)}, \frac{p_{15}\, q_{235}\,
(q_{125}-p_{32})}{(b_1-b_3)(b_2-b_3)(c_5-c_6)},
-\frac{p_{24}}{b_2-b_3},-\frac{p_{23}\,
q_{235}}{(b_2-b_3)(c_3-c_4)},0,1 \)
\end{array}
\end{equation}

Theorem \ref{thm:extra_of_S:A} and Lemma
\ref{lemma:extra_of_S:eigenbasis_of_C} are proved by direct
computation as well as the following two lemmas.

\begin{lemma}
\label{lemma:extra:evect:B} Let $\B$ be as in
(\ref{equation:extra_of_S:B,C}). Then $\w_1$ and $\w_2$ are
eigenvectors of $\B$ with the eigenvalue $b_1$, $\w_3$ and $\w_4$
are eigenvectors of $\B$ with the eigenvalue $b_2$, and $\w_5$ and
$\w_6$ are eigenvectors of $\B$ with the eigenvalue $b_3$.
\end{lemma}

\begin{lemma}
\label{lemma:extra:orthogonality} Let $\<*,*\>$ be defined by
(\ref{equation:extra:form}). Let $V_{b_1}$ be the subspace of $V$
spanned by $\w_1$ and $\w_2$. Let $V_{b_2}$ be the subspace of $V$
spanned by $\w_3$ and $\w_4$. Let $V_{b_3}$ be the subspace of $V$
spanned by $\w_5$ and $\w_5$. Then $V_{b_1}$, $V_{b_2}$, and
$V_{b_3}$ are mutually orthogonal with respect to $\<*,*\>$.
\end{lemma}

\noindent  Theorem \ref{thm:extra:form} follows from Lemmas
\ref{lemma:extra_of_S:eigenbasis_of_C} and
\ref{lemma:extra:orthogonality}. Finally, Theorem
\ref{thm:extra:sign_of_form} can be proved similarly to Theorems
\ref{thm:hg:sign_of_form} and \ref{thm:even:sign_of_form}.

\section{Indecomposable Triple Flag Varieties with Finitely Many Orbits}
\label{sec:triple_flag_var}

 Let $i\in\{1,2,3\}$. For a triple of flags
$\emptyset =V^i_0\subset V^i_1\subset V^i_2\subset\cdots \subset
V^i_{k_i-1}\subset V^i_{k_i}=V$, we call the {\it dimension vector
in the jump coordinates} the vector
$\((\dim (V^1/ V^1_0),\dim (V^1_2/ V^1_1),\cdots ,\right.$ \\
$\left. \dim (V^1_{k_1}/ V^1_{k_1-1})),(\dim (V^2_1/ V^2_0),\dim
(V^2_2/ V^2_1), \cdots ,\dim (V^2_{k_2}/ V^2_{k_2-1})),(\dim
(V^3_1/ V^3_0),
\dim (V^3_2/ V^3_1),\cdots ,\right.$ \\
$\left. \dim (V^3_{k_3}/ V^3_{k_3-1}))\)$. We say that this triple
of flags is in a {\it standard form}, if $V$ is given a basis
$\z_1,\cdots ,\z_n$ with the following property: for the flag
$\emptyset =V^2_0\subset V^2_1\subset V^2_2\subset\cdots \subset
V^2_{k_2-1}\subset V^2_{k_2}=V$, the subspace $V^2_i$ of dimension
$d^2_i$ is spanned by the first $d^2_i$ basis vectors $\z_1$,
$\z_2,\cdots$; for the flag $\emptyset =V^3_0\subset V^3_1\subset
V^3_2\subset\cdots \subset V^3_{k_3-1}\subset V^3_{k_3}=V$, the
subspace $V^3_j$ of dimension $d^3_j$ is
spanned by the last $d^3_j$ basis vectors $\z_n$, $\z_{n-1},\cdots$. \\

 P.~Magyar, J.~Weyman, and A.~Zelevinsky classify in \cite{MWZ} all indecomposable
triple partial flag varieties with finitely many orbits of the
diagonal action of the general linear group. They work over an
algebraically closed field. ${\mathbb C}$ is enough for our
purposes. Among other results, they give the dimension vectors in
the jump coordinates as well as explicit representatives of the
open orbit in the standard form. For every element of their list,
the first flag turns to be just a subspace $V_1^1\subset V$. It
also turns out that this subspace is spanned by vectors ${\bf
a}_i$ such that all their coordinates in the standard basis
$\z_1,\cdots ,\z_n$ are equal to either $0$ or $1$. Their list is
given on page \pageref{equation:MWZ_list}.

\begin{remark}{\rm
Our definition of a standard form for a triple of flags is weaker
than that of Magyar, Weyman, and Zelevinsky (includes more triple
flags).}
\end{remark}

 Recall that we proved Theorem \ref{thm:irr} only for the hypergeometric family
so far. Now we use the results of Magyar, Weyman and Zelevinsky to
prove the counterparts of this result for all other families of
Simpson. Let us begin with the even family. Recall that we work
with the normalized matrices $\A$, $\B$, and $\C$. This means that
they are traceless and the eigenvalues of $\A$ are $1$ and $-1$.
Let ${\bf Z}$ be the following matrix

\begin{center}
\pspicture[](0,0)(9,7) \label{picture:even:Z}
\psset{unit=1cm,linewidth=1pt,linecolor=blue} \rput(0,-0.5){
\pspolygon(1,1)(7,1)(7,7)(1,7) \psline(2,1)(2,7) \psline(4,1)(4,7)
\psline(1,4)(7,4) \psline(1,6)(7,6) \rput(0,4.1){${\bf Z}=$}
\rput(7.8,4.1){,where} \rput(1.5,7.3){$1$} \rput(3.1,7.3){$m-1$}
\rput(5.4,7.2){$m$} \rput(0.7,2.7){$m$} \rput(0.4,5){$m-1$}
\rput(0.8,6.6){$1$} \rput(1.5,6.5){$1$} \rput(3,6.5){$0$}
\rput(5.4,6.5){$0$} \rput(1.5,5){$0$} \rput(3.1,5){$Z_{1+i,1+j}$}
\rput(5.4,5){$0$} \rput(1.5,2.7){$0$} \rput(3,2.7){$0$}
\rput(5.4,2.7){$Z_{m+i,m+j}$} }
\endpspicture
\end{center}

\begin{equation}
\label{equation:Z}
\begin{array}{lll}
Z_{1+i,1+j} & = & \left\{
\begin{array}{ll}
0 & if~~i<j \\
1 & if~~i=j \\
\frac{\prod\limits_{k=1+j}^i
q_{k,2m-i}}{\prod\limits_{k=2m+1-i}^{2m-j} (c_{2m-i}-c_k)} &
if~~i>j
\end{array},
\right. \\
 & & \\
Z_{m+i,m+j} & = & \left\{
\begin{array}{ll}
0 & if~~i<j \\
1 & if~~i=j \\
\frac{\prod\limits_{k=m+j}^{m-1+i}
q_{m+1-i,k}}{\prod\limits_{k=m+2-i}^{m+1-j} (c_{m+1-i}-c_k)} &
if~~i>j
\end{array}.
\right.
\end{array}
\end{equation}

\noindent Note that ${\bf Z}$ is lower-triangular with all the diagonal elements equal to $1$. \\

 For $1\le i\le 2m$, let $\z_i={\bf Z}\, \e_i$. The matrix ${\bf Z}$ is non-degenerate, so
$\z_i$ is a basis of $V$. Consider the following flags:
$V^2_1\subset V^2_2\subset V$ and $V^3_1\subset V^3_2\subset
\cdots \subset V^3_{2m-1}\subset V$ where $V^2_1$ is spanned by
$\z_1$, $V^2_2$ is spanned by $\z_1,\z_2,\cdots ,\z_m$, and
$V^3_i$ is spanned by $\z_{2m},\z_{2m-1},\cdots ,\z_{2m+1-i}$.
They are the second and the third flags of the even family with
the dimension vector $((m,m),(1,m-1,m),(1^{2m}))$ in
(\ref{equation:MWZ_list}).

\begin{equation}
\label{equation:MWZ_list}
\begin{tabular}{|l|l|}
\hline
\multicolumn{2}{|c|}{ hypergeometric family } \\
\hline
 & \\
$(m-1,1),(1^m),(1^m)$ & ${\bf a}_k=\z_1+\z_{k+1}~(1\le k\le m-1)$ \\
$(1,m-1),(1^m),(1^m)$ & ${\bf a}_1=\z_1+\z_2+\cdots +\z_m$ \\
 & \\
\hline
\multicolumn{2}{|c|}{ even family } \\
\hline
 & \\
$(m,m),(1,m-1,m),(1^{2m})$ & ${\bf a}_k=\z_1+\z_{k+1}+\z_{2m+1-k}~(1\le k\le m-1)$, \\
 & ${\bf a}_m=\z_1+\z_{m+1}$ \\
$(m,m),(1,m,m-1),(1^{2m})$ & ${\bf a}_1=\z_1+\z_2$, \\
 & ${\bf a}_k=\z_1+\z_{k+1}+\z_{2m+2-k}~(2\le k\le m)$ \\
$(m,m),(m-1,m,1),(1^{2m})$ & ${\bf a}_k=\z_k+\z_{2m-k}+\z_{2m}~(1\le k\le m-1)$, \\
 & ${\bf a}_m=\z_m+\z_{2m}$ \\
$(m,m),(m-1,1,m),(1^{2m})$ & ${\bf a}_k=\z_k+\z_m+\z_{2m+1-k}~(1\le k\le m-1)$, \\
 & ${\bf a}_m=\z_m+\z_{m+1}$ \\
$(m,m),(m,m-1,1),(1^{2m})$ & ${\bf a}_1=\z_1+\z_{2m}$, \\
 & ${\bf a}_k=\z_k+\z_{2m+1-k}+\z_{2m}~(2\le k\le m)$ \\
$(m,m),(m,1,m-1),(1^{2m})$ & ${\bf a}_1=\z_1+\z_{m+1}$, \\
 & ${\bf a}_k=\z_k+\z_{m+1}+\z_{2m+2-k}~(2\le k\le m)$ \\
 & \\
\hline
\multicolumn{2}{|c|}{ odd family } \\
\hline
 & \\
$(m,m+1),(1,m,m),(1^{2m+1})$ & ${\bf a}_k=\z_1+\z+{k+1}+\z_{2m+2-k}~(1\le k\le m)$\\
$(m+1,m),(1,m,m),(1^{2m+1})$ & ${\bf a}_1=\z_1+\z_2$, ${\bf a}_{m+1}=\z_1+\z_{m+2}$, \\
 & ${\bf a}_k=\z_1+\z_{k+1}+\z_{2m+3-k}~(2\le k\le m)$ \\
$(m,m+1),(m,1,m),(1^{2m+1})$ & ${\bf a}_k=\z_k+\z_{m+1}+\z_{2m+2-k}~(1\le k\le m)$ \\
$(m+1,m),(m,1,m),(1^{2m+1})$ & ${\bf a}_1=\z_1+\z_{m+1}$, ${\bf a}_{m+1}=\z_{m+1}+\z_{m+2}$, \\
 & ${\bf a}_k=\z_k+\z_{m+1}+\z_{2m+3-k}~(2\le k\le m)$ \\
$(m,m+1),(m,m,1),(1^{2m+1})$ & ${\bf a}_k=\z_k+\z_{2m+1-k}+\z_{2m+1}~(1\le k\le m)$ \\
$(m+1,m),(m,m,1),(1^{2m+1})$ & ${\bf a}_1=\z_1+\z_{2m+1}$, ${\bf a}_{m+1}=\z_{m+1}+\z_{2m+1}$, \\
 & ${\bf a}_k=\z_k+\z_{2m+2-k}+\z_{2m+1}~(2\le k\le m)$ \\
 & \\
\hline
\multicolumn{2}{|c|}{ ${\hat E}_8$ family } \\
\hline
 & \\
$(2,4),(2,2,2),(1,1,1,1,1,1)$ & ${\bf a}_1=\z_1+\z_2+\z_3+\z_6$, ${\bf a}_2=\z_1+\z_4+\z_5$ \\
$(4,2),(2,2,2),(1,1,1,1,1,1)$ & ${\bf a}_1=\z_1+\z_5$, ${\bf a}_2=\z_2+\z_3$, \\
 & ${\bf a}_3=\z_2+\z_5+\z_6$, ${\bf a}_4=\z_4+\z_5$ \\
 & \\
\hline
\multicolumn{2}{|c|}{ $E_8$ family } \\
\hline
 & \\
$(3,3),(2,2,2),(2,1,1,1,1)$ & ${\bf a}_1=\z_1+\z_2+\z_3$, ${\bf
a}_2=\z_1+\z_6$,
${\bf a}_3=\z_2+\z_4+\z_5$ \\
$(3,3),(2,2,2),(1,2,1,1,1)$ & ${\bf a}_1=\z_1+\z_2+\z_4+\z_6$, \\
 & ${\bf a}_2=\z_1+\z_3$, ${\bf a}_3=\z_1+\z_5$ \\
$(3,3),(2,2,2),(1,1,2,1,1)$ & ${\bf a}_1=\z_1+\z_5+\z_6$, ${\bf
a}_2=\z_2+\z_3+\z_6$,
${\bf a}_3=\z_4+\z_5$ \\
$(3,3),(2,2,2),(1,1,1,2,1)$ & ${\bf a}_1=\z_1+\z_2+\z_4+\z_6$,
${\bf a}_2=\z_1+\z_3$, ${\bf a}_3=\z_1+\z_5$ \\
$(3,3),(2,2,2),(1,1,1,1,2)$ & ${\bf a}_1=\z_1+\z_4+\z_6$, ${\bf
a}_2=\z_2+\z_4+\z_5$,
${\bf a}_3=\z_2+\z_3$ \\
 & \\
\hline
\end{tabular}
\end{equation}
\newpage

 Here is an example of the matrix ${\bf Z}$ with $m=3$.

\begin{example}
\label{example:even:Z}
$$
{\bf Z}=\[
\begin{array}{c|cc|ccc}
1 & 0 & 0 & 0 & 0 & 0 \\
 &  &  &  &  &  \\
\hline
 &  &  &  &  &  \\
0 & 1 & 0 & 0 & 0 & 0 \\
 &  &  &  &  &  \\
0 & \frac{q_{24}}{c_4-c_5} & 1 & 0 & 0 & 0 \\
 &  &  &  &  &  \\
\hline
 &  &  &  &  &  \\
0 & 0 & 0 & 1 & 0 & 0 \\
 &  &  &  &  &  \\
0 & 0 & 0 & \frac{q_{24}}{c_2-c_3} & 1 & 0 \\
 &  &  &  &  &  \\
0 & 0 & 0 & \frac{q_{14}\, q_{15}}{(c_1-c_2)(c_1-c_3)} &
\frac{q_{15}}{c_1-c_2} & 1
\end{array} \]
$$
\end{example}

\begin{definition}
\label{definition:spec_flag} Let $\A$ be a diagonalizable complex
linear operator with the spectrum $(\lambda_1,\lambda_2,\cdots
,\lambda_k)$. Let $V_{\lambda_i}$ be the eigenspace of $\A$
corresponding to the eigenvalue $\lambda_i$. We will call the flag
$V_{\lambda_1}\subset V_{\lambda_1}\oplus V_{\lambda_2} \subset
\cdots \subset V$ the spectral flag of $\A$ corresponding to the
ordering $(\lambda_1,\lambda_2,\cdots ,\lambda_k)$ of its
spectrum.
\end{definition}

 If $(m_1,m_2,\cdots ,m_k)$ are the multiplicities of the spectrum of $\A$ from the
above definition, then the dimension vector in the jump
coordinates of its spectral
flag is also $(m_1,m_2,\cdots ,m_k)$. \\

 If we take another look at the eigenvectors of $\B$ (Lemma \ref{lemma:even:evec_of_B})
and at the eigenvectors of $\C$ (Lemma
\ref{lemma:even:eigenbasis_of_C}), we see that the spectral flags
of these matrices are exactly the second and the third flags of
the Magyar, Weyman, Zelevinsky triple
$((m,m),(1,m-1,m),(1^{2m}))$.

\begin{lemma}
\label{lemma:even:spec_flag_of_A} The subspace $V^1_1$ spanned by
the vectors ${\bf a}_1,\cdots ,{\bf a}_m$ (from the
$((m,m),(1,m-1,m),(1^{2m}))$ entry in (\ref{equation:MWZ_list}))
is the spectral subspace of the matrix $\A$ corresponding to the
eigenvalue $-1$.
\end{lemma}

\proof In order to prove $\( \A + \Id \)\, {\bf a}_i = 0$ for
$1\le i\le m$, we have to prove the following identities.
\begin{enumerate}
\item The first identity says that the first component of $(\A+\Id)\, {\bf a}_i$ is zero.

\begin{equation}
\label{equation:z1} b_1+c_{2m}+1 + \sum_{j=i}^{m-1} Z_{1+j,1+i}\,
B_{1,1+j}  + \sum_{m+1-i}^m Z_{m+j,2m+1-i}\, B_{1,m+j} = 0
\end{equation}

\item The second identity says that the components $2$ through $m$ of
$(\A+\Id)\, {\bf a}_i$ are zero.

\begin{equation}
\label{equation:z2}
C_{1+j,1}+(c_{2m-j}+b_2+1)Z_{1+j,1+i}+\sum_{k=m+1-i}^m
Z_{m+k,2m+1-i}\, B_{1+j,m+k}=0
\end{equation}

\item The third identity says that the components $m+1$ through $2m$ of
$(\A+\Id)\, {\bf a}_i$ are zero.

\begin{equation}
\label{equation:z3}
C_{m+j,1}+(c_{m+1-j}+b_3+1)Z_{m+j,2m+1-i}+\sum_{k=i}^{m-1}
Z_{1+k,1+i}\, C_{m+j,1+k} = 0
\end{equation}
\end{enumerate}

 Recall that the matrix elements of the matrices $\B$ and $\C$ are given by
formulas (\ref{equation:even:B:II}), (\ref{equation:even:B:III}),
(\ref{equation:even:B:VI}), (\ref{equation:even:C:IV}),
(\ref{equation:even:C:VII}), (\ref{equation:even:C:VIII}) on pages
\pageref{equation:even:B:II} -- \pageref{equation:even:C:VIII}.
Then the first identity becomes
\begin{equation}
\label{equation:z4} b_1+c_{2m}+1 + \sum_{j=i}^{m-1}
\frac{\prod\limits_{k=1+i}^m q_{k,2m-j}}
{\prod\limits_{k=m+1\atop{k\ne 2m-j}}^{2m-i} (c_{2m-j}-c_k)} +
\sum_{j=m+1-i}^m p^{31}_{m+1-j}\,
\frac{\prod\limits_{k=2m+1-i}^{2m-1} q_{m+1-j,k}}
{\prod\limits_{k=1\atop{k\ne m+1-j}}^i (c_{m+1-j}-c_k)} = 0.
\end{equation}
In our normalized version, $p^{31}_{m+1-j}=c_{m+1-j}+b_3-1$. Thus,
(\ref{equation:z4}) splits into two identities of homogeneous
degrees $0$ and $1$. The part of degree $0$ is
$$
\sum_{j=m+1-i}^m \frac{\prod\limits_{k=2m+1-i}^{2m-1} q_{m+1-j,k}}
{\prod\limits_{k=1\atop{k\ne m+1-j}}^i (c_{m+1-j}-c_k)} = 1.
$$
This identity is equivalent to (\ref{equation:hg_id_2}) from
Appendix. The part of degree $1$ is
$$
b_1+c_{2m} + \sum_{j=i}^{m-1} \frac{\prod\limits_{k=1+i}^m
q_{k,2m-j}} {\prod\limits_{k=m+1\atop{k\ne 2m-j}}^{2m-i}
(c_{2m-j}-c_k)} + \sum_{j=m+1-i}^m (b_3+c_{m+1-j})\,
\frac{\prod\limits_{k=2m+1-i}^{2m-1} q_{m+1-j,k}}
{\prod\limits_{k=1\atop{k\ne m+1-j}}^i (c_{m+1-j}-c_k)} = 0.
$$
This one is proved similarly to identity (\ref{equation:eq1}) on
page \pageref{equation:eq1} with the help of identity
(\ref{equation:hg_id_3}) from the Appendix applied separately to each
of the sums. \\

 To prove the second identity (\ref{equation:z2}), let us recall that
$$
Z_{1+j,1+i} = \left\{
\begin{array}{ll}
0 & if~~j<i \\
1 & if~~j=i \\
\frac{\prod\limits_{k=1+i}^j
q_{k,2m-j}}{\prod\limits_{k=2m+1-j}^{2m-i} (c_{2m-j}-c_k)} &
if~~j>i
\end{array}.
\right.
$$
Thus, the second identity (\ref{equation:z2}) splits into three
different identities. In case $j<i$, we have
$$
C_{1+j,1}+\sum_{k=m+1-i}^m Z_{m+k,2m+1-i}\, B_{1+j,m+k}=0.
$$
This identity, after cancelling out common multiples, splits into
a sum of two identities: one of degree $-1$ and the other of
degree $0$. The first reduces to identity (\ref{equation:hg_id_1})
and the second reduces to
identity (\ref{equation:hg_id_2}) from the Appendix. \\

 In case $i=j$, we have $C_{1+j,1}+(c_{2m-j}+b_2+1)+\sum_{k=m+1-i}^m Z_{m+k,2m+1-i}\,
B_{1+j,m+k}=0$. This identity splits into two parts of degree $0$
and $1$. The part of degree $0$ reduces to identity
(\ref{equation:hg_id_5}) from the Appendix. The part of degree $1$ is
a sum of two identities: one is equivalent to
(\ref{equation:hg_id_2}) and the other is equivalent to
(\ref{equation:hg_id_5}) from the Appendix. The case $j>i$,
after cancelling out common multiples, becomes equivalent to the case $j=i$. \\

 To prove the third identity (\ref{equation:z3}), recall that
$$
Z_{m+j,2m+1-i}=\left\{
\begin{array}{ll}
0 & if~~i+j<m+1 \\
1 & if~~i+j=1 \\
\frac{\prod\limits_{k=2m+1-i}^{m-1+j}
q_{m+1-j,k}}{\prod\limits_{k=m+2-j}^i (c_{m+1-j}-c_k)} &
if~~i+j>m+1
\end{array}.
\right.
$$
Thus, the third identity splits into three different identities.
In case $i+j<m+1$, after cancellig out common multiples,
(\ref{equation:z3}) becomes equivalent to identity
(\ref{equation:hg_id_2}) from the Appendix. In cases $i+j=m+1$ and
$i+j>m+1$, identity
(\ref{equation:hg_id_7}) from the Appendix does the job. \\

 Finally, it is clear from (\ref{equation:MWZ_list}) and (\ref{equation:Z}) that for
$1\le i\le m$, the vectors ${\bf a}_i$ are linearly independent.
\endproof

 Now let us get back to the proof of Theorem \ref{thm:irr}.
\label{proof:irr}

\noindent {\it Proof of Theorem \ref{thm:irr} ---}
For a triple of partitions $S''(\alpha,\beta,\gamma)$ from
Simpson's list (\ref{equation:S_list}), we want to prove that if
$(\s(\A),\s(\B),\s(\C))$ is a point of $S''(\alpha,\beta,\gamma)$,
then the triple $(\A,\B,\C)$ is irreducible. We have proved it
already for the hypergeometric family in subsection
\ref{subsec:proofs:hg}. To prove it for the even family, let us
recall that if $(\s(\A),\s(\B),\s(\C))\in
S''(\alpha,\beta,\gamma)$, then the scalar product
(\ref{equation:even:form}) is well-defined and non-degenerate.
Assume that the triple of matrices $(\A,\B,\C)$ is reducible. Then
they all preserve a non-trivial subspace $V'$. This subspace is
spanned by some eigenvectors of $\C$. Let $V''$ be the subspace of
$V$ spanned by the complementary eigenvectors of $\C$. But all
eigenvectors of $\C$ are orthogonal to each other with respect to
(\ref{equation:even:form}). Thus, the space $V$ splits into the
orthogonal direct sum $V'\oplus V''$. Thus, the matrices $\A$ and
$\B$ preserve the subspace $V''$ as well. So, no matter how we
introduce linear orders on the spectra of $\A$, $\B$, and $\C$,
the corresponding triple of flags will decompose. However, if
$\(\s(\A), \s(\B), \s(\C)\)$ is a point of
$S''((m,m),(1,m-1,m),(1^{2m}))$, then as follows from Lemma
\ref{lemma:even:spec_flag_of_A} and the preceeding discussion, the
spectral flags of the matrices $\A$, $\B$, and $\C$ give the
Magyar, Weyman, Zelevinsky representative of the open orbit of the
corresponding triple flag variety. According to Magyar, Weyman,
and Zelevinsky, this triple of flags is indecomposable. Thus, the
assumption that the triple
$(\A,\B,\C)$ is reducible cannot be true. This proves Theorem \ref{thm:irr} in the even case. \\

 The $Z$-matrix for the odd family triple $((m+1,m),(1,m,m),(1^{2m+1}))$ can be obtained from
the $Z$-matrix for the even family triple
$((m+1,m+1),(1,m,m+1),(1^{2m+2}))$ by restricting the latter to
$V^s_{m+1}$ as in the proof of Theorem \ref{thm:OM_as_subm_of_EM}
on page \pageref{proof1}. The rest of the argument is the same.
Finally, let us give the $Z$-matrix for the extra case of Simpson
(or, more precisely, for the triple of compositions
$((4,2),(2,2,2),(1^6))$ of the ${\hat E}_8$-familly from
(\ref{equation:MWZ_list})).

\begin{equation}
\label{equation:^E_8-Z} Z=\[
    \begin{array}{cc|cc|cc}
       1 & 0 & 0 & 0 & 0 & 0 \\
        &  &  &  &  &  \\
       \frac{q_{145}}{c_5-c_6} & 1 & 0 & 0 & 0 & 0 \\
        &  &  &  &  &  \\
       \hline
        &  &  &  &  &  \\
       0 & 0 & 1 & 0 & 0 & 0 \\
        &  &  &  &  &  \\
       0 & 0 & \frac{q_{245}}{c_3-c_4} & 1 & 0 & 0 \\
        &  &  &  &  &  \\
       \hline
        &  &  &  &  &  \\
       0 & 0 & 0 & 0 & 1 & 0 \\
       0 & 0 & 0 & 0 & \frac{q_{236}}{c_1-c_2} & 1
    \end{array}
  \].
\end{equation}
\endproof

 The only family of Magyar, Weyman, Zelevinsky which does not appear in the
list of Simpson, is the $E_8$-family. We now construct the
matrices $\A=\B+\C, \B, \C$ such that their spectral flags form
the Magyar, Weyman, Zelevinsky representative for the open orbit
of the triple flag variety of dimension
$((3,3),(2,2,2),(1,1,1,1,2))$. This time the standard basis $\e_i$
and the $z$-basis of Magyar, Weyman, and Zelevinsky
$\z_i$ coinside ($\z_i=\e_i$). \\

 Let $(a_1,a_1,a_1,a_2,a_2,a_2,b_1,b_1,b_2,b_2,b_3,b_3,c_1,c_2,c_3,c_4,c_5,c_5)\in
S((3,3),(2,2,2),(1,1,1,1,2))$.

\begin{equation}
\label{equation:nonsph:B} \mbox{Let~~} \B = \[
\begin{array}{cc|cc|cc}
b_1 & 0 & 0 & \begin{array}{c} a_1+a_2 \\ -b_1-b_3 \\ -c_1-c_5
\end{array} &
\begin{array}{c} -a_1-a_2+ \\ b_1+b_3+ \\ c_1+c_5 \end{array} & -a_2+b_3+c_1 \\
 & & & & & \\
0 & b_1 & a_1-b_1-c_5 & \begin{array}{c} -a_1-2a_2+ \\
b_1+b_2+b_3+ \\ c_1+c_3+c_5 \end{array} & \begin{array}{c}
-a_1-a_2+ \\ b_2+b_3+ \\c_2+c_4 \end{array} &
\begin{array}{c} a_1+2a_2 \\ -b_1-b_2-b_3 \\-c_1-c_3-c_5 \end{array} \\
 & & & & & \\
\hline
 & & & & & \\
0 & 0 & b_2 & 0 & -a_1+b_2+c_4 & \begin{array} {c} a_1+2a_2 \\ -b_1-b_2-b_3 \\
-c_1-c_3-c_5 \end{array} \\
 & & & & & \\
0 & 0 & 0 & b_2 & \begin{array}{c} 2a_1+a_2 \\ -b_1-2b_2 \\
-c_1-c_3-c_5 \end{array} &
-a_2+b_3+c_1 \\
 & & & & & \\
\hline
 & & & & & \\
0 & 0 & 0 & 0 & b_3 & 0 \\
 & & & & & \\
0 & 0 & 0 & 0 & 0 & b_3 \\
\end{array} \],
\end{equation}

\begin{equation}
\label{equation:nonsph:C} \C = \[
\begin{array}{cc|cc|cc}
c_5 & 0 & 0 & 0 & 0 & 0 \\
 & & & & & \\
0 & c_5 & 0 & 0 & 0 & 0 \\
 & & & & & \\
\hline
 & & & & & \\
\begin{array}{c} -a_1-2a_2+ \\ b_1+b_2+b_3+ \\ c_1+c_3+c_5 \end{array} & a_1-b_2-c_4 &
c_4 & 0 & 0 & 0 \\
 & & & & & \\
\begin{array}{c} a_1+a_2 \\ -b_2-b_3 \\ -c_1-c_3 \end{array} & \begin{array}{c} -a_1-a_2+ \\
b_1+b_2+ \\ c_4+c_5 \end{array} & \begin{array}{c} a_1+a_2 \\
-b_1-b_2 \\ -c_4-c_5 \end{array} &
c_3 & 0 & 0 \\
 & & & & & \\
\hline
 & & & & & \\
\begin{array}{c} -a_1-2a_2+ \\ b_1+b_2+b_3+ \\ c_1+c_3+c_5 \end{array} & \begin{array}{c}
-a_1-a_2+ \\ b_1+b_2+ \\ c_4+c_5 \end{array} & \begin{array}{c} a_1+a_2 \\ -b_1-b_2 \\
-c_4-c_5 \end{array} & \begin{array}{c} -a_1-2a_2+ \\ b_1+b_2+b_3+
\\ c_1+c_3+c_5 \end{array} &
c_2 & 0 \\
 & & & & & \\
-a_2+b_1+c_5 & 0 & 0 & \begin{array}{c} a_1+a_2 \\ -b_1-b_3 \\
-c_1-c_5 \end{array} &
\begin{array}{c} -a_1-a_2+ \\ b_1+b_3+ \\ c_1+c_5 \end{array} & c_1
\end{array}
\].
\end{equation}

 It is clear that $\B$ and $\C$ are diagonalizable and that $\s(\B)=\{b_1,b_1,b_2,b_2,b_3,b_3\}$,
$\s(\C)=\{c_1,c_2,c_3,c_4,c_5,c_5\}$. The following is proved by
direct computation.

\begin{theorem}
\label{thm:nonsph:A} For $\B$ and $\C$ as in
(\ref{equation:nonsph:B}) and (\ref{equation:nonsph:C}), let
$\A=\B+\C$. Then $\A$ is diagonalizable and
$\s(\A)=\{a_1,a_1,a_1,a_2,a_2,a_2\}$.
\end{theorem}

 The existence of a scalar product on $V$ such that $\A$, $\B$, and $\C$ are self-adjoint
with respect to it can be proved using methods of the theory of
quiver representations. It follows from Schur's lemma that if
$(\s(\A),\s(\B),\s(\C))$ is a generic point in
$S((3,3),(2,2,2),(1,1,1,1,2))$, then the form is unique up to a
constant multiple. It would be interesting to find a basis in
which the form is ``nice'' (for instance, having matrix entries as
ratios of products of linear forms in the eigenvalues of $\A$,
$\B$, and $\C$).

\section{Connections with the Littlewood-Richardson rule}
\label{sec:LR}

An irreducible rational representation of $GL(n,\mathbb{C})$ is
uniquely determined
by its {\it highest weight} \\
$\lambda=(\lambda_1,\lambda_2,\cdots,\lambda_n)$ where $\lambda_i$
are integers such that $\lambda_1\ge\lambda_2\ge
\cdots\ge\lambda_n$. We can decompose tensor products of
irreducible representations into sums of irreducibles:

\begin{equation}
\label{iequation:LR-coeff} V_{\lambda}\otimes
V_{\mu}=\sum_{\nu}c_{\lambda\mu}^{\nu}V_{\nu}.
\end{equation}

The number $c_{\lambda\mu}^{\nu}$ of copies of $V_{\nu}$ in
$V_{\lambda}\otimes V_{\mu}$ is called the {\it
Littlewood-Richardson coefficient}. There exists a famous
combinatorial algorithm to compute the Littlewood-Richardson
coefficient called the {\it Littlewood-Richardson rule} (see
\cite{F2} for more information). It follows from the results of
A.~Klyachko \cite{Kl} combined with a refinement by A.~Knutson and
T.~Tao \cite{KT}, that the lattice points of the Klyachko cone are
exactly the triples of highest weights with non-zero
Littlewood-Richardson coefficients (see also \cite{F1} for a nice
survey). The question whether all the lattice points of the
Klyachko cone were such triples was raised in \cite{Z} under the
name of the {\it saturation conjecture}. The conjecture was proved
by Knutson and Tao in \cite{KT}. Some of the Klyachko inequalities
describing the Klyachko cone are redundant. Knutson, Tao, and
Woodward give in \cite{KTW} the set of necessary inequalities for
the Klyachko cone. H.~Derksen and J.~Weyman in \cite{DW1} give a
proof of the saturation conjecture different from that of Knutson
and Tao. They use methods of the theory of quiver representations,
developing further ideas of A.~Schofield from \cite{Sch1}.
Moreover, in \cite{DW2}, Derksen and Weyman give description of
{\it all} the faces of the Klyachko cone of arbitrary dimension.
However, all these results involve recursive computations. \\

    The inequalities of Theorem \ref{thm:hg:sign_of_form}, Theorem
\ref{thm:even:sign_of_form}, Theorem \ref{thm:odd:sign_of_form},
and Theorem \ref{thm:extra:sign_of_form} give nonrecursive
description of some faces of the Klyachko cone. Thus, integral
solutions to these inequalities give non-recursive description of
some triples of highest weights with $c_{\lambda\mu}^{\nu}\ne 0$. \\

    Let us show a different way to derive these inequalities and also
show that the corresponding $c_{\lambda\mu}^{\nu}=1$. For that, we
use the {\it Berenstein-Zelevinsky triangle}. It was invented in
\cite{BZ} as a geometric version of the Littlewood-Richardson
rule. A variation of the BZ--triangle was used in \cite{KT} under
the name of a {\it honeycomb tinkertoy}. A different variation of
the triangle was used in \cite{GP} under the name of a {\it
web-function} to examine relations of the Littlewood-Richardson
coefficients with a Quantum-Yang-Baxter-Type equation. \\

 Consider the following graph.

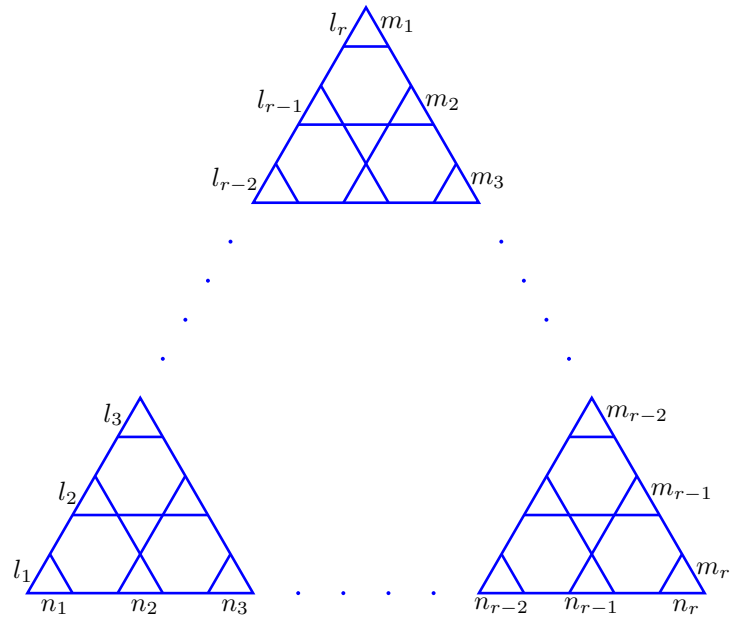
\begin{center}
\pspicture[](0,0)(7,9) \label{BZtr}
\psset{unit=0.6cm,linewidth=1pt,linecolor=blue} \rput{0}(-2,1.5){
\pspolygon(0,0)(5,0)(2.5,4.33) \psline(1,0)(0.5,0.866)
\psline(3,0)(1.5,2.598) \psline(2,0)(3.5,2.598)
\psline(4,0)(4.5,0.866) \psline(1,1.732)(4,1.732)
\psline(2,3.464)(3,3.464) \rput{0}(10,0){
\pspolygon(0,0)(5,0)(2.5,4.33) \psline(1,0)(0.5,0.866)
\psline(3,0)(1.5,2.598) \psline(2,0)(3.5,2.598)
\psline(4,0)(4.5,0.866) \psline(1,1.732)(4,1.732)
\psline(2,3.464)(3,3.464) } \rput{0}(5,8.66){
\pspolygon(0,0)(5,0)(2.5,4.33) \psline(1,0)(0.5,0.866)
\psline(3,0)(1.5,2.598) \psline(2,0)(3.5,2.598)
\psline(4,0)(4.5,0.866) \psline(1,1.732)(4,1.732)
\psline(2,3.464)(3,3.464) } \psdots[dotsize=1.5pt
0](6,0)(7,0)(8,0)(9,0)(3,5.196)(3.5,6.062)(4,6.928)(4.5,7.794)
\psdots[dotsize=1.5pt
0](12,5.196)(11.5,6.062)(11,6.928)(10.5,7.794)
\rput{0}(0.6,-0.3){$n_1$} \rput{0}(2.6,-0.3){$n_2$}
\rput{0}(4.6,-0.3){$n_3$} \rput{0}(10.5,-0.3){$n_{r-2}$}
\rput{0}(12.5,-0.3){$n_{r-1}$} \rput{0}(14.6,-0.3){$n_r$}
\rput{0}(-0.1,0.5){$l_1$} \rput{0}(0.9,2.232){$l_2$}
\rput{0}(1.9,3.964){$l_3$} \rput{0}(4.6,9.16){$l_{r-2}$}
\rput{0}(5.6,10.892){$l_{r-1}$} \rput{0}(6.9,12.624){$l_r$}
\rput{0}(15.2,0.5){$m_r$} \rput{0}(14.5,2.232){$m_{r-1}$}
\rput{0}(13.5,3.964){$m_{r-2}$} \rput{0}(10.2,9.16){$m_3$}
\rput{0}(9.2,10.892){$m_2$} \rput{0}(8.2,12.624){$m_1$} }
\endpspicture
\end{center}

This is the Berenstein-Zelevinsky triangle for $sl_r$. In order to
define it formally, it is convenient to use the baricentric
coordinates in~$\R^2$. Namely, we represent a point in~$\R^2$ by a
triple $(\alpha,\beta,\gamma)$ such that $\alpha+\beta+\gamma=0$.
The {\it $r$-th Berenstein-Zelevinsky triangle} $\BZ_r$ is the set
of points in~$\R^2$ with baricentric coordinates $(\alpha, \beta,
\gamma)$, such that

\begin{enumerate}
\item
$0 < \beta < -\alpha < r+1$;
\item
the numbers $2\alpha$, $2\beta$, and $2\gamma$ are integers;
\item
at least one $\alpha$, $\beta$, or $\gamma$ is not integer.
\end{enumerate}

Every integer point $(a,b,c)$, $a+b+c=0$, with $0 < b < -a < r+1$
has six neighbors in $\BZ_r$ that form vertices of the following
hexagon: \psset{xunit=.3cm,yunit=.3cm}
\begin{center}
\pspicture[.1](-3,-4)(9,16) \psset{linecolor=blue}
\psline[linewidth=.5pt,showpoints=false]{-}(0,0)(6,0)(9,6)(6,12)(0,12)
(-3,6)(0,0)
  \qdisk(0,0){1.5pt}
  \qdisk(6,0){1.5pt}
  \qdisk(-3,6){1.5pt}
  \qdisk(9,6){1.5pt}
  \qdisk(0,12){1.5pt}
  \qdisk(6,12){1.5pt}

 \rput[u](3,6){$(a,b,c)$}
 \rput[u](-5,-2){$E=(a+{1\over 2},\,b,\,c-{1\over 2})$}
 \rput[u](11,-2){$D=(a,\,b+{1\over 2},\,c-{1\over 2})$}
 \rput[u](-5,14){$A=(a,\,b-{1\over 2},\,c+{1\over 2})$}
 \rput[u](11,14){$B=(a-{1\over 2},\,b,\,c+{1\over 2})$}
 \rput[u](-10,6){$F=(a+{1\over 2},\,b-{1\over 2},\,c)$}
 \rput[u](16,6){$C=(a-{1\over 2},\,b+{1\over 2},\,c)$}
\endpspicture
\end{center}

\begin{definition}
\label{definition:BZ_filling} {\rm A function
$f:\BZ_r\to\{0,1,2,\dots\}$ is called a {\it BZ-filling} if for
any hexagon as above we have $f(A)+f(B)=f(D)+f(E)$,
$f(B)+f(C)=f(E)+f(F)$, and $f(C)+f(D)=f(F)+f(A)$ (the last
condition follows from the first two). We call this the hexagon
condition. }
\end{definition}

Let $\lambda=\sum_{i=1}^{r}l_i\,\omega_i$,
$\mu=\sum_{i=1}^{r}m_j\,\omega_i$, and
$\nu=\sum_{i=1}^{r}n_i\,\omega_i$, where the $\omega_i$ are the
{\it fundamental weights} of $sl(r+1,{\mathbb C})$. Let us assign
$l_i$, $m_i$, and $n_i$ to the boundary segments of the $\BZ_r$ as
shown on the picture on page \pageref{BZtr}. Note that
$l_i=\lambda_i-\lambda_{i+1}$, $m_i=\mu_i-\mu_{i+1}$,
$n_i=\nu_i-\nu_{i+1}$.

\begin{definition}
\label{definition:BZ_boundary_condition} A filling $f$ of $\BZ_r$
{\it satisfies boundary conditions} if for any boundary segment
with vertices $A$, $B$, and a nonnegative integer value $v$
assigned to the segment, $f(A)+f(B)=v$.
\end{definition}

\begin{theorem}[Berenstein,~Zelevinsky]
\label{thm:BZ} Let $\lambda=\sum_{i=1}^{r}l_i\omega_i$,
$\mu=\sum_{i=1}^{r}m_j\omega_i$, and
$\nu=\sum_{i=1}^{r}n_i\omega_i$ be dominant weights of
$sl(r+1,{\mathbb C})$. Then $c_{\lambda\mu}^{\nu}=\#\{${\it of
fillings of $\BZ_r$ satisfying the boundary conditions.}$\}$
\end{theorem}

 Let us use the BZ-triangle for a different proof of Theorem \ref{thm:hg:sign_of_form},
and also to show that the corresponding Littlewood-Richardson coefficient is equal to one. \\

\noindent {\it Proof of Theorem \ref{thm:hg:sign_of_form} ---}
Let us assume that $a_1<a_2$. Consider $\BZ_r$ for the
hypergeometric case ($r=m-1$). For that, we have to switch from
$gl(n,{\mathbb C})$ to $sl(n,{\mathbb C})$. Let us set
$\tA=\A-\frac{1}{r+1}\,\trace(\A)\,\Id$,
$\tB=\B-\frac{1}{r+1}\,\trace(\B)\,\Id$, and
$\tC=\C-\frac{1}{r+1}\,\trace(\C)\,\Id$. Then
$\s(\tA)=\{\ta,\ta,\cdots,\ta,-r\,\ta\}$. Let us call the
eigenvalues of $\tB$ and $\tC$ $\tb_i$ and $\tc_j$. We have
$\sum_{j=1}^{r+1} \tb_j=\sum_{j=1}^{r+1} \tc_j=0$. Recall that
$l_i=\tb_i-\tb_{i+1}$; $m_i=\tc_i-\tc_{i+1}$; $n_1=n_2=\cdots
=n_{r-1}=0$, and $n_r=(r+1)a$.

\begin{center}
\pspicture[](0,0)(16,44) \label{HGtr}
\psset{unit=1cm,linewidth=1pt,linecolor=blue} \rput{0}(-5,0){
\pspolygon(0,0)(5,0)(2.5,4.33) \psline(1,0)(0.5,0.866)
\psline(3,0)(1.5,2.598) \psline(2,0)(3.5,2.598)
\psline(4,0)(4.5,0.866) \psline(1,1.732)(4,1.732)
\psline(2,3.464)(3,3.464)

\rput{0}(-0.25,0){$0$}

\rput{0}(1.15,0.25){$0$} \rput{0}(1.85,0.25){$0$}
\rput{0}(1.15,1.55){$0$} \rput{0}(1.85,1.55){$0$} \rput{0}(2,0){
\rput{0}(1.15,0.25){$0$} \rput{0}(1.85,0.25){$0$}
\rput{0}(1.15,1.55){$0$} \rput{0}(1.85,1.55){$0$} }
\rput{0}(1,1.73){ \rput{0}(1.15,1.55){$0$}
\rput{0}(1.85,1.55){$0$} } \rput{0}(0.8,0.9){$l_1$}
\rput{0}(2.8,0.9){$l_1$} \rput{0}(4.8,0.9){$l_1$}
\rput{0}(1.8,2.63){$l_2$} \rput{0}(3.8,2.63){$l_2$}
\rput{0}(2.8,4.36){$l_3$}

\rput{0}(10,0){ \pspolygon(0,0)(5,0)(2.5,4.33)
\psline(1,0)(0.5,0.866) \psline(3,0)(1.5,2.598)
\psline(2,0)(3.5,2.598) \psline(4,0)(4.5,0.866)
\psline(1,1.732)(4,1.732) \psline(2,3.464)(3,3.464)

\rput{0}(-0.25,0){$0$}

\rput{0}(1.15,0.25){$0$} \rput{0}(1.85,0.25){$0$}
\rput{0}(1.15,1.55){$0$} \rput{0}(1.85,1.55){$0$} \rput{0}(2,0){
\rput{0}(1.15,0.25){$0$} } \rput{0}(0.8,0.9){$l_1$}
\rput{0}(2.8,0.9){$l_1$} \rput{0}(1.8,2.63){$l_2$} }
\rput{0}(5,8.66){ \pspolygon(0,0)(5,0)(2.5,4.33)
\psline(1,0)(0.5,0.866) \psline(3,0)(1.5,2.598)
\psline(2,0)(3.5,2.598) \psline(4,0)(4.5,0.866)
\psline(1,1.732)(4,1.732) \psline(2,3.464)(3,3.464)

\rput{0}(-0.25,0){$0$}

\rput{0}(1.15,0.25){$0$} \rput{0}(1.85,0.25){$0$}
\rput{0}(1.15,1.55){$0$} \rput{0}(1.85,1.55){$0$} \rput{0}(2,0){
\rput{0}(1.15,0.25){$0$} } \rput{0}(0.95,0.9){$l_{r-2}$}
\rput{0}(2.95,0.9){$l_{r-2}$} \rput{0}(1.95,2.63){$l_{r-1}$} }
\psdots[dotsize=1.5pt
0](6,0)(7,0)(8,0)(9,0)(3,5.196)(3.5,6.062)(4,6.928)(4.5,7.794)
\psdots[dotsize=1.5pt
0](12,5.196)(11.5,6.062)(11,6.928)(10.5,7.794)
\rput{0}(0.5,-0.2){$0$} \rput{0}(2.5,-0.2){$0$}
\rput{0}(4.5,-0.2){$0$} \rput{0}(10.5,-0.2){$0$}
\rput{0}(12.5,-0.2){$0$} \rput{0}(14.6,-0.2){$n_r$}
\rput{0}(-0.1,0.5){$l_1$} \rput{0}(0.9,2.232){$l_2$}
\rput{0}(1.9,3.964){$l_3$} \rput{0}(4.7,9.16){$l_{r-2}$}
\rput{0}(5.7,10.892){$l_{r-1}$} \rput{0}(6.9,12.624){$l_r$}
\rput{0}(15.1,0.5){$m_r$} \rput{0}(14.3,2.232){$m_{r-1}$}
\rput{0}(13.3,3.964){$m_{r-2}$} \rput{0}(10.1,9.16){$m_3$}
\rput{0}(9.1,10.892){$m_2$} \rput{0}(8.1,12.624){$m_1$} }
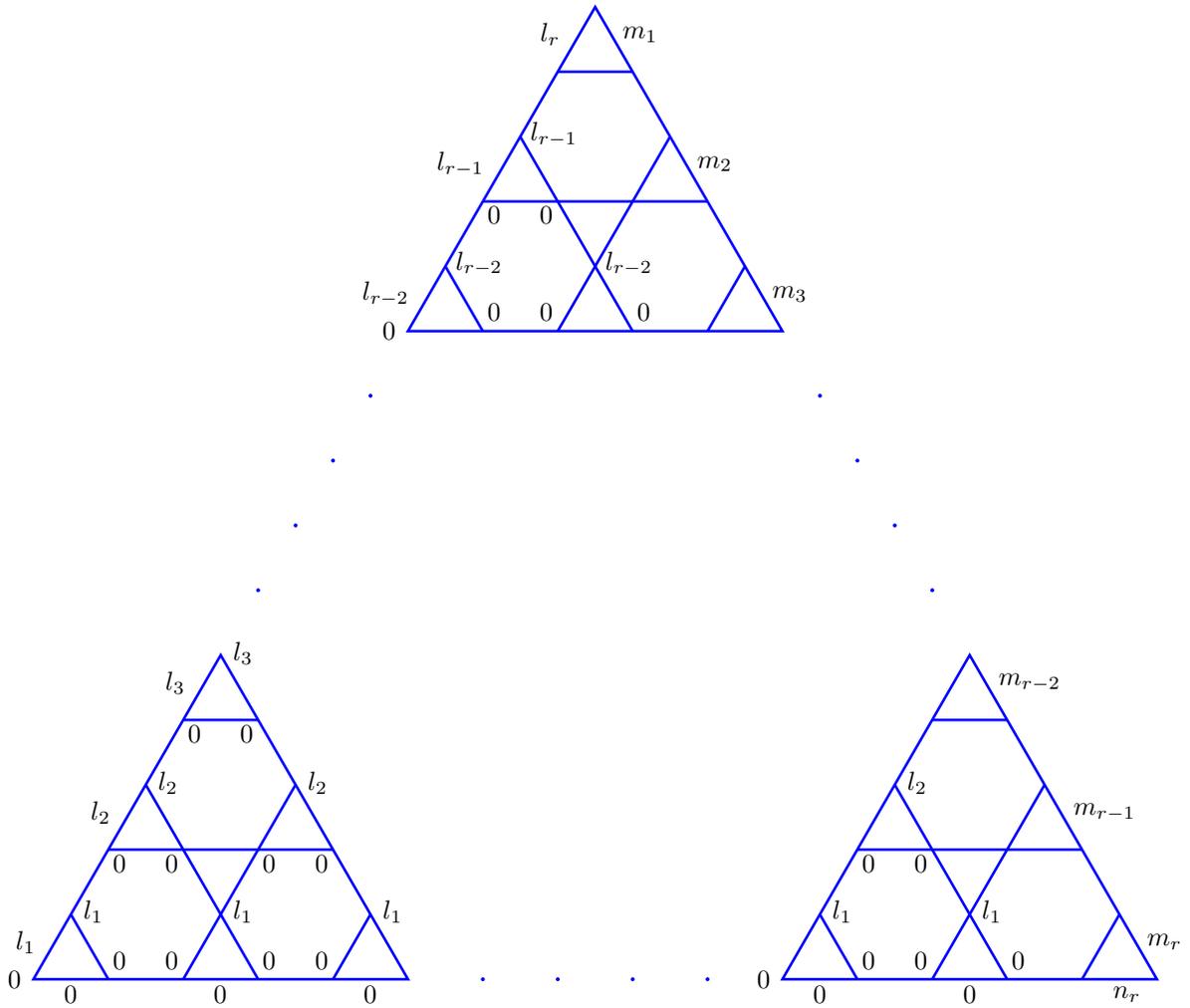
\endpspicture
\end{center}

 If a hexagon has two zeros on a side, then the nonnegativity of BZ-fillings
and the hexagon condition force two zeros on the opposite side. In
the hypergeometric case, this mechanism reduces the BZ-triangle to
a strip. Let us put a variable $x$ in an unfilled vertex of the
strip. Then the filling is expressed in terms of $x$, and the
boundary conditions. Also the $l_r$-boundary condition gives a
linear equation on $x$.

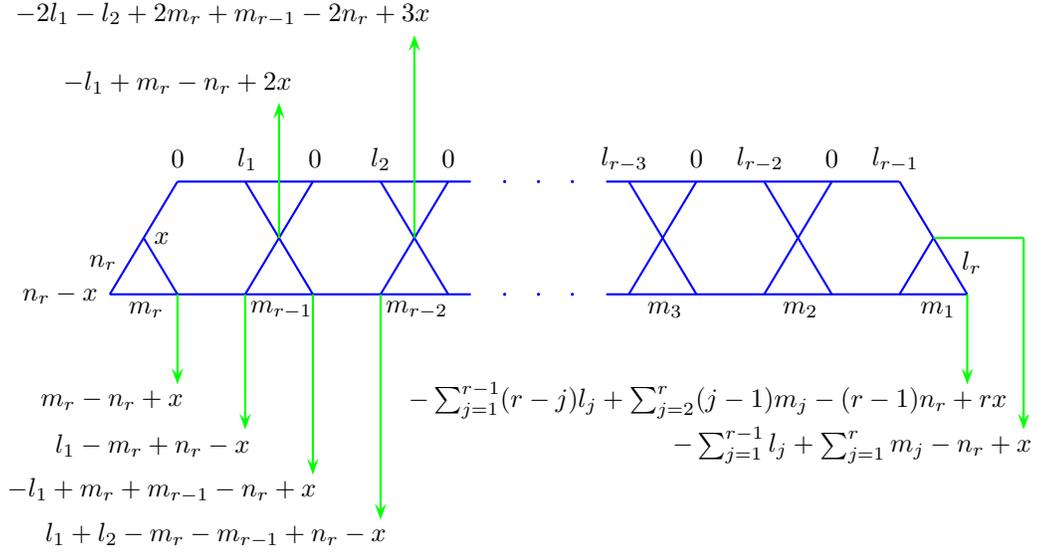
\pspicture(30,30) \psset{unit=0.03cm, linecolor=blue}
\rput(30,135){
 \psline(20,0)(180,0)
 \psline(240,0)(400,0)
 \psline(20,0)(50,50)
 \psline(50,50)(180,50)
 \psline(240,50)(370,50)
 \psline(80,0)(110,50)
 \psline(140,0)(170,50)
 \psline(250,0)(280,50)
 \psline(310,0)(340,50)
 \psline(370,0)(385,25)
 \psline(400,0)(370,50)
 \psline(340,0)(310,50)
 \psline(280,0)(250,50)
 \psline(50,0)(35,25)
 \psline(110,0)(80,50)
 \psline(170,0)(140,50)

 \pscircle*(195,0){0.7}
 \pscircle*(210,0){0.7}
 \pscircle*(225,0){0.7}
 \pscircle*(195,50){0.7}
 \pscircle*(210,50){0.7}
 \pscircle*(225,50){0.7}

 \rput(-2,0){$n_r-x$}
 \rput(43,25){$x$}
 \rput(50,60){$0$}
 \rput(81,60){$l_1$}
 \psline[linecolor=green]{->}(95,25)(95,85)
 \rput(50,94){$-l_1+m_r-n_r+2x$}
 \rput(111,60){$0$}
 \rput(140,60){$l_2$}
 \psline[linecolor=green]{->}(155,25)(155,115)
 \rput(70,124){$-2l_1-l_2+2m_r+m_{r-1}-2n_r+3x$}
 \rput(170,60){$0$}
 \psline[linecolor=green]{->}(50,0)(50,-40)
 \rput(21,-47){$m_r-n_r+x$}
 \psline[linecolor=green]{->}(80,0)(80,-60)
 \rput(39,-67){$l_1-m_r+n_r-x$}
 \psline[linecolor=green]{->}(110,0)(110,-80)
 \rput(43,-87){$-l_1+m_r+m_{r-1}-n_r+x$}
 \psline[linecolor=green]{->}(140,0)(140,-100)
 \rput(67,-107){$l_1+l_2-m_r-m_{r-1}+n_r-x$}
 \rput(36,-7){$m_r$}
 \rput(96,-7){$m_{r-1}$}
 \rput(156,-7){$m_{r-2}$}
 \rput(266,-7){$m_3$}
 \rput(326,-7){$m_2$}
 \rput(387,-7){$m_1$}
 \psline[linecolor=green]{->}(400,0)(400,-40)
 \rput(285,-47){$-\sum_{j=1}^{r-1}(r-j)l_j+\sum_{j=2}^r (j-1)m_j -(r-1)n_r+rx$}
 \psline[linecolor=green](385,25)(425,25)
 \psline[linecolor=green]{->}(425,25)(425,-60)
 \rput(349,-67){$-\sum_{j=1}^{r-1} l_j +\sum_{j=1}^r m_j -n_r +x$}
 \rput(402,14){$l_r$}
 \rput(368,60){$l_{r-1}$}
 \rput(340,60){$0$}
 \rput(308,60){$l_{r-2}$}
 \rput(17,14){$n_r$}
 \rput(280,60){$0$}
 \rput(248,60){$l_{r-3}$}
}
\endpspicture

The $l_r$-boundary condition gives us the equation
$-\sum_{j=1}^{r-1}(r+1-j)l_j+\sum_{j=1}^{r}jm_j-rn_r+(r+1)x=l_r$.
Thus,
$$
x=\displaystyle{\frac{\sum_{j=1}^{r}(r+1-j)l_j-\sum_{j=1}^{r}jm_j+rn_r}{r+1}}
$$
and the filling is defined uniquely. Let us list the Klyachko
inequalities. First, $x={\tilde b}_1+{\tilde c}_n+r{\tilde a}>0$.
However, this inequality is not a generating one. If we have
another look at the strip above, we see that $x$ has a neighboring
$0$-vertex. So, $x$ is the sum of the numbers at the opposite edge
($-{\tilde b}_2-{\tilde c}_{n-1}+{\tilde a}$, and ${\tilde
b}_1+{\tilde b}_2+ {\tilde c}_{n-1}+{\tilde c}_n+(n-2){\tilde
a}$). All the numbers in the middle part of the strip except for
the atmost right one ($-{\tilde b}_n-{\tilde c}_1+{\tilde a}$) do
not produce generating inequalities for the same reason. The
numbers on the lower part of the strip together with the last
middle number produce the following generating inequalities:

\begin{equation}
\label{equation:hg:Kl_ineq:a_2>a_1}
\begin{array}{lcccl}
\tb_1+\tc_{m-1} & > &     & > & \tb_1+\tc_m     \\
\tb_2+\tc_{m-2} & > &     & > & \tb_2+\tc_{m-1} \\
\cdots          &   & \ta &   & \cdots          \\
\tb_{m-1}+\tc_1 & > &     & > & \tb_{m-1}+\tc_2 \\
                &   &     & > & \tb_m+\tc_1.
\end{array}
\end{equation}

Switching back to $gl(n,{\mathbb C})$ proves the theorem in this
case. The case $a_1>a_2$ is obtained from the case $a_1<a_2$ in
the following way. Let us multiply $\A$, $\B$, and $\C$ by $-1$.
Then let us renumber $b_i$ and $c_j$ so that $b_1>b_2>\cdots >b_m$
and $c_1>c_2>\cdots >c_m$ again.
\endproof

 One can similarly prove Theorems \ref{thm:even:sign_of_form},
\ref{thm:odd:sign_of_form}, \ref{thm:extra:sign_of_form} and also
show that the
corresponding Littlewood-Richardson coefficients are equal to one. \\

 In the $E_8$ case of the Magyar, Weyman, Zelevinsky  list (\ref{equation:MWZ_list}),
we do not have an explicit criterion for positivity of the
corresponding Hermitian form. However, the BZ-triangle enables us
to compute the inequalities on the eigenvalues of $\A$, $\B$, and
$\C$ which make the form sign-definite. In the notations of
(\ref{equation:nonsph:B}), (\ref{equation:nonsph:C}), let
$a_1>a_2$ and $b_1>b_2>b_3$. Then the form is sign-definite
precisely in the following situations.

\begin{equation}
\label{equation:E_8_pos_ineq}
\begin{array}{|lcccl|}
\hline
 & & & & \\
-a_1-a_2+b_1+b_3+c_1+c_5 & > & 0 & > & -a_1-a_2+b_1+b_3+c_2+c_5 \\
-a_1-a_2+b_1+b_2+c_3+c_5 & > & 0 & > & -a_1-a_2+b_1+b_2+c_4+c_5 \\
a_1-a_2-c_1+c_2-c_3+c_4  & > & 0 & > & a_1-a_2-c_1-c_2+c_3+c_4  \\
                         &   & 0 & > & -a_2+b_2+c_5             \\
 & & & & \\
\hline
 & & & & \\
-a_1-a_2+b_1+b_3+c_2+c_5 & > & 0 & > & -a_1-a_2+b_1+b_3+c_3+c_5 \\
-a_2+b_2+c_5             & > & 0 & > & -a_2+b_2+c_4             \\
-a_1+a_2+c_1-c_2+c_3-c_4 & > & 0 &   &                          \\
                         &   & 0 & > & -a_1-a_2+b_2+b_3+c_1+c_5 \\
                         &   & 0 & > & -a_1+b_1+c_4             \\
 & & & & \\
\hline
 & & & & \\
-a_1-a_2+b_2+b_3+c_1+c_5 & > & 0 & > & -a_1-a_2+b_2+b_3+c_2+c_5 \\
-a_1-a_2+b_1+b_2+c_3+c_5 & > & 0 & > & -a_1-a_2+b_1+b_2+c_4+c_5 \\
a_1-a_2-c_1+c_2-c_3+c_4  & > & 0 &   &                          \\
-a_1+b_1+c_5             & > & 0 &   &                          \\
                         &   & 0 & > & -a_2+b_3+c_5             \\
 & & & & \\
\hline
 & & & & \\
-a_1-a_2+b_1+b_3+c_2+c_5 & > & 0 & > & -a_1-a_2+b_1+b_3+c_3+c_5 \\
-a_1+b_2+c_1             & > & 0 & > & -a_1+b_2+c_5             \\
-a_1-a_2+b_1+b_2+c_4+c_5 & > & 0 &   &                          \\
-a_2+b_3+c_1             & > & 0 &   &                          \\
                         &   & 0 & > & a_1-a_2-c_1+c_2-c_3+c_4  \\
 & & & & \\
\hline
 & & & & \\
-a_1-a_2+b_1+b_3+c_3+c_5 & > & 0 & > & -a_1-a_2+b_1+b_3+c_4+c_5 \\
-a_1-a_2+b_2+b_3+c_1+c_5 & > & 0 & > & -a_1-a_2+b_2+b_3+c_2+c_5 \\
-a_1+a_2+c_1+c_2-c_3-c_4 & > & 0 & > & -a_1+a_2+c_1-c_2+c_3-c_4 \\
-a_1+b_2+c_5             & > & 0 & > &                          \\
 & & & & \\
\hline
\end{array}
\end{equation}
\bigskip

 The first set of inequalities forces $c_1>c_2>c_3>c_4>c_5$ realizing the
dimension vector $(3,3),(2,2,2),(1,1,1,1,2)$. The second set of
inequalities forces $c_1>c_2>c_3>c_5>c_4$ realizing the dimension
vector $(3,3),(2,2,2),(1,1,1,2,1)$. The third set of inequalities
forces $c_1>c_2>c_5>c_3>c_4$ realizing the dimension vector
$(3,3),(2,2,2),(1,1,2,1,1)$. The forth set of inequalities forces
$c_1>c_5>c_2>c_3>c_4$ realizing the dimension vector
$(3,3),(2,2,2),(1,2,1,1,1)$. The last set of inequalities forces
$c_5>c_1>c_2>c_3>c_4$ realizing the dimension vector
$(3,3),(2,2,2),(2,1,1,1,1)$. Thus, all the members of the
$E_8$-family from (\ref{equation:MWZ_list}) can be constructed
this way with the help of the corresponding eigenvectors of $\A$,
$\B$, and $\C$.

\section{Fuchsian systems, Fuchsian equations, Okubo normal forms, and the list
of Haraoka-Yokoyama.}
\label{sec:applications_to_Fuchsian_systems}

    Let us consider a system of linear differential equations on a
${\mathbb C}^n$-valued function $f$ on $\CP^1$.

\begin{equation}
\label{def:lin-reg_sys}
df=\omega f,
\end{equation}

\noindent where $\omega$ is a ($n \times n$) matrix-valued
1-differential form on $\mathbb{CP}^1$ . Let the form be
holomorphic everywhere on $\CP^1$ except for a finite set of
points ${\cal D}=\{z_1, z_2, \cdots z_k \}$. Let us consider a
solution of (\ref{def:lin-reg_sys}) restricted to a sectorial
neighborhood centered at any $z_i \in {\cal D}$. If any such
solution has polynomial growth when it approaches $z_i$ within any
such sector, then the system (\ref{def:lin-reg_sys}) is called
{\it linear regular}. If $\omega$ has only first order poles at
$\cal D$, then the system is called {\it Fuchsian}. Any Fuchsian
system is linear regular, but there exist linear regular systems
which are not Fuchsian (for more detailed treatment, see
\cite{Bo} or \cite{V1}). \\

    An $n$ order Fuchsian equation is a linear differential equation

\begin{equation}
\label{def:Fuc_eq}
f^{n}(z)+q_1(z)f^{n-1}(z)+\cdots +q_n(z)f(z)=0
\end{equation}

\noindent such that its coefficients $q_j(z)$ have a finite set of
poles ${\cal D}=\{z_1, z_2, \cdots z_k \}$ and in a small
neighborhood of a pole $z_i$ the coefficients of
(\ref{def:Fuc_eq}) have the form

\begin{equation}
\label{coef_of_Fuc_eq}
q_j(z)=\displaystyle{\frac{r_j(z)}{(z-z_i)^j}},~~~j=1,\cdots ,n,
\end{equation}

\noindent where the $r_j(z)$ are holomorphic functions. Solutions
of Fuchsian equations have polynomial growth when continued
analytically towards a pole. This distinguishes Fuchsian
differential equations from all other linear differential
equations on $\CP^1$. Thus, for linear differential equations the
notions ``Fuchsian'' and ``linear regular'' coincide. \\

%
%
%
%

    The matrix $R_i=\res_{z=z_i} \omega (z)$ is called the {\it
residue} of a linear regular system at $z_i$. By Cauchy residue
theorem, $\sum_{i=1}^k R_i=0$.

\begin{theorem}[see \cite{Bo}]
\label{thm:standard_ form_of_Fuchsian_system} Any Fuchsian system
has the standard form

\begin{equation}
\label{def:stand_form_of_Fuch_sys} \frac{df}{dz}=\sum_{i=1}^k
\frac{R_i}{z-z_i} f(z).
\end{equation}
\end{theorem}

\begin{theorem}[see \cite{Bo}]
\label{thm:Fuchsian_syst_from_eq} For any Fuchsian equation on the
Riemann sphere, it is possible to construct a Fuchsian system with
the same singular points and the same monodromy. Converse is not
true.
\end{theorem}

\begin{remark}
Thus, the notion of a residue matrix makes sense for a Fuchsian
equation as well as for a Fuchsian system.
\end{remark}

    To study Fuchsian differential equations, K.~Okubo invented
what became later known as the {\it Okubo normal form} of a
Fuchsian equation. In \cite{O}, he proves that any Fuchsian
equation can be written in the following form:
\begin{equation}
\label{equation:okubo_nf}
(t\,Id-B)\displaystyle{\frac{dx}{dt}}=A\, x,
\end{equation}
where t is a complex variable, $x\in {\mathbb C}^n$ is an unknown
vector, $Id$ is the identity matrix of order $n$, $B$ is a
constant diagonal $n\times n$ matrix, and $A$ is a constant
$n\times n$ matrix. Let
\begin{equation}
\label{equation:O:B}
B=\mbox{diag}(\underbrace{z_1,\cdots,z_1}_{n_1},
\underbrace{z_2,\cdots,\z_2}_{n_2},\cdots
,\underbrace{z_k,\cdots,z_k}_{n_k}),
\end{equation}
such that $z_i\ne z_j$ for $i\ne j$, $n_1+n_2+\cdots +n_k=n$, and
$n_1\ge n_2\ge\cdots \ge n_k$. The partition $(n_1,n_2,\cdots
,n_k)$ of $n$ endows $A$ with the block decomposition
$A=(A_{ij})_{1\le i,j\le d}$. Let us call $\Lambda_i$ the set of
eigenvalues of $A_{ii}$ and let us call $\Lambda_{\infty}$ the set
of eigenvalues of $A$. Then $z_1$, $z_2,\cdots ,z_k$ and $\infty$
are the singular points of (\ref{equation:okubo_nf}). At the point
$z_i$, (\ref{equation:okubo_nf}) has $n_i$ non-holomorphic
solutions with local exponents $\lambda_j\in \Lambda_i$. At
$\infty$, (\ref{equation:okubo_nf}) has local exponents
$\lambda_j\in \Lambda_{\infty}$. \\

    In \cite{Y1}, T.~Yokoyama used Okubo theory to classify the
spectral types of rigid irreducible Fuchsian equations. For such,
all $A_{ii}$ are diagonalizable as well as $A$ itself. Quoting the
result of Yokoyama, we will not give the spectral types of
$A_{ii}$ and $A$ the way he does. Instead, we will list spectral
types of the residue matrices (which are diagonalizable, too).

\begin{equation}
\label{equation:yokoyama:list}
\begin{tabular}{|l|l|l|}
\hline
 & & \\
$I$ & $(m-1,1),(1^m),(1^m)$ & $m\ge 2$ \\
 & & \\
\hline
 & & \\
$II$ & $(m,1^m),(m,1^m),(m,m-1,1)$ & $m\ge 2$ \\
 & & \\
\hline
 & & \\
$III$ & $(m,1^{m+1}),(m+1,1^m),(m,m,1)$ & $m\ge 2$ \\
 & & \\
\hline
 & & \\
$IV$ & $(4,1,1),(2,1,1,1,1),(2,2,2)$ & \\
 & & \\
\hline
 & & \\
$I^{*}$ & $\underbrace{(m-1,1),\cdots ,(m-1,1)}_{m~times}$ & $m\ge 2$ \\
 & & \\
\hline
 & & \\
$II^{*}$ & $(m,1^m),(m+1,1^{m-1}),(2m-1,1),(m,m)$ & $m\ge 2$ \\
 & & \\
\hline
 & & \\
$III^{*}$ & $(m+1,1^m),(m+1,1^m),(2m,1),(m+1,m)$ & $m\ge 2$ \\
 & & \\
\hline
 & & \\
$IV^{*}$ & $(4,1,1),(4,1,1),(4,1,1),(4,2)$ & \\
 & & \\
\hline
\end{tabular}
\end{equation}
\medskip

    Y.~Haraoka explicitly constructed the equations of the above
spectral types in \cite{H1}. In \cite{H2}, he explored the
solutions of these equations: computed their monodromies, found
monodromy invariant forms in their spaces of solutions, etc. It
turns out that the solutions of these equations are important
hypergeometric functions. It also turns out that the Fuchsian
systems constructed in our paper are closely related to
Yokoyama-Haraoka equations: sometimes the $A$ matrices are just
the same! We think it is interesting to understand the nature of
this relation, find solutions to our systems, and their
monodromies. We plan to do it in a subsequent publication.

\section{Appendix}
\label{sec:appendix} In this section we collect the identities
needed for the proofs in the previous sections.

 For $k<n-1$,
\begin{equation}
\label{equation:hg_id_1} \sum\limits_{i=1}^n
\displaystyle{\frac{\prod\limits_{j=1}^k (x_i-y_j)}
{\prod\limits_{j=1\atop j\ne i}^n (x_i-x_j)}} = 0
\end{equation}

\begin{equation}
\label{equation:hg_id_2} \sum\limits_{i=1}^n
\displaystyle{\frac{\prod\limits_{j=1}^{n-1} (x_i-y_j)}
{\prod\limits_{j=1\atop j\ne i}^n (x_i-x_j)}} = 1
\end{equation}

\begin{equation}
\label{equation:hg_id_3} \sum\limits_{i=1}^n
\displaystyle{\frac{\prod\limits_{j=1}^n (x_i-y_j)}
{\prod\limits_{j=1\atop j\ne i}^n (x_i-x_j)}} =
\sum\limits_{i=1}^n (x_i-y_i)
\end{equation}

\begin{equation}
\label{equation:hg_id_4} \sum\limits_{i=1}^n
\displaystyle{\frac{\prod\limits_{j=1}^{n+1} (x_i+y_j)}
{\prod\limits_{j=1\atop j\ne i}^n (x_i-x_j)}}=\sum\limits_{i=1}^n
x^2_i + \sum\limits_{1\le i<j\le n} x_i x_j + \sum\limits_{1\le
i<j\le n+1} y_i y_j + \(\sum\limits_{i=1}^n
x_i\)\(\sum\limits_{i=1}^{n+1} y_i\)
\end{equation}

 For $1\le i\le m-1$,
\begin{equation}
\label{equation:hg_id_5} \sum\limits_{j=1}^m
\displaystyle{\frac{\prod\limits_{k=1\atop k\ne i}^m
(x_j-y_k)}{\prod\limits_{k=1\atop k\ne j}^m
(x_j-x_k)}}\,\displaystyle{ \frac{\prod\limits_{k=1\atop k\ne j}^m
(y_i-x_k)}{\prod\limits_{k=1\atop k\ne i}^m (y_i-y_k)}}=1
\end{equation}

\begin{equation}
\label{equation:hg_id_6} y_i^2
\displaystyle{\frac{\prod\limits_{k=1\atop k\ne i}^{m-1}
(y_i-y_k)}{\prod\limits_{k=1}^m (y_i-x_k)}}+\sum\limits_{j=1}^m
\displaystyle{\frac{x_j^2}{(y_i-x_j)^2}}\,\displaystyle{\frac
{\prod\limits_{k=1}^{m-1} (x_j-y_k)}{\prod\limits_{k=1\atop k\ne
j}^m (x_j-x_k)}}=1
\end{equation}

For $1\le i\le m$,
\begin{equation}
\label{equation:hg_id_7}
\displaystyle{\frac{\prod\limits_{k=1\atop k\ne i}^m
(x_i-x_k)}{\prod\limits_{k=1}^{m-1} (x_i+y_k)}}\, +
\,\sum\limits_{j=1}^{m-1} \displaystyle{\frac{1}{x_i+y_j}}\,
\displaystyle{ \frac{\prod\limits_{k=1\atop k\ne i}^m
(y_j+x_k)}{\prod\limits_{k=1\atop k\ne j}^{m-1} (y_j-y_k)}}=1
\end{equation}

\begin{equation}
\label{equation:hg_id_8} \sum\limits_{i=1}^n
\displaystyle{\frac{(x_i+x_1)(x_i+x_2)\cdots\widehat{(x_i+x_i)}\cdots
(x_i+x_n)}{(x_i-x_1)(x_i-x_2)\cdots\widehat{(x_i-x_i)}\cdots(x_i-x_n)}}=\left\{
\begin{array}{l}
0,~~\mbox{if~~}n~~\mbox{is even}\\
1,~~\mbox{if~~}n~~\mbox{is odd}
\end{array}
\right.
\end{equation}

\begin{equation}
\label{equation:hg_id_9} \sum\limits_{i=1}^n
x_i\,\displaystyle{\frac{(x_i+x_1)(x_i+x_2)\cdots\widehat{(x_i+x_i)}\cdots
(x_i+x_n)}{(x_i-x_1)(x_i-x_2)\cdots\widehat{(x_i-x_i)}\cdots(x_i-x_n)}}=x_1+x_2+\cdots
+x_n
\end{equation}

\begin{equation}
\label{equation:hg_id_10} \sum\limits_{i=1}^n
x_i^2\,\displaystyle{\frac{(x_i+x_1)(x_i+x_2)\cdots\widehat{(x_i+x_i)}\cdots
(x_i+x_n)}{(x_i-x_1)(x_i-x_2)\cdots\widehat{(x_i-x_i)}\cdots(x_i-x_n)}}=(x_1+x_2+\cdots
+x_n)^2
\end{equation}

\begin{equation}
\label{equation:hg_id_11}
y_j\,\displaystyle{\frac{\prod\limits_{k=1\atop k\ne j}^m
(y_j+y_k)}{\prod\limits_{k=1}^m (y_j-x_k)}}+\sum\limits_{r=1}^m
\displaystyle{\frac{x_r}{x_r-y_j}}\,\displaystyle{\frac{\prod
\limits_{k=1\atop k\ne j}^m (x_r+y_k)}{\prod\limits_{k=1\atop k\ne
r}^m (x_r-x_k)}}=1~~~ \mbox{for~~}j=1,2,\cdots ,m-1
\end{equation}

\begin{equation}
\label{equation:hg_id_12}
y_j^2\,\displaystyle{\frac{\prod\limits_{k=1\atop k\ne j}^m
(y_j+y_k)}{\prod\limits_{k=1}^m (y_j-x_k)}}+\sum\limits_{r=1}^m
\displaystyle{\frac{x_r^2}{x_r-y_j}}\,\displaystyle{\frac{\prod
\limits_{k=1\atop k\ne j}^m (x_r+y_k)}{\prod\limits_{k=1\atop k\ne
r}^m (x_r-x_k)}}=\sum \limits_{r=1}^m
(x_r+y_r)~~~\mbox{for~~}j=1,2,\cdots ,m-1
\end{equation}

For $i=1,2,\cdots m-1$ and $j=1,2,\cdots ,m$,
\begin{equation}
\label{equation:hg_id_13} \sum\limits_{r=1}^{m-1}
\displaystyle{\frac{\prod\limits_{k=1\atop k\ne j}^m (x_k+y_r)}
{\prod\limits_{k=1\atop k\ne r}^{m-1}
(y_r-y_k)}}\,\sum\limits_{s=1}^m
\displaystyle{\frac{1}{x_s+y_i}}\,\displaystyle{\frac{\prod\limits_{k=1\atop
k\ne r}^m (x_s^2-y_k^2)}{\prod\limits_{k=1\atop k\ne s}^m
(x_s^2-x_k^2)}} +
\displaystyle{\frac{x_j+y_m}{x_j+y_i}}\,\displaystyle{\frac{\prod\limits_{k=1\atop
k\ne i}^m (x_j-y_k)}{\prod\limits_{k=1\atop k\ne j}^m (x_j+x_k)}}
= 1
\end{equation}

 All these identities have the following features: the left hand side $L(x,y)$ is a
rational homogeneous function in $x_i$ and $y_j$. All the
denominators of $L(x,y)$ are products of linear forms $\alpha$ of
the form $(x_i \pm x_j)$, $(y_i \pm y_j)$, or $(x_i \pm y_j)$. The
power of every such form in any denominator is $1$. The right hand
sides $R(x,y)$ are constants or
homogeneous polynomials in $x_i$ and $y_j$ of degree $1$ or $2$. \\

 The first step to prove such an identity is to prove that $L(x,y)$ is in fact a
polynomial. For that, it is enough to prove that $\alpha\,
L(x,y)|_{\alpha=0}=0$ for every form $\alpha$ from any denominator
of the identity. For all the identities except for
(\ref{equation:hg_id_13}), the restriction of $\alpha\, L(x,y)$ to
the hyperplane $\alpha=0$ turns to be a sum of just two terms with
equal absolute values and opposite signs. For example, consider
identity (\ref{equation:hg_id_2}). Let us fix $p$ and $q$ such
that $1\le p<q\le n$. Consider

$$
(x_p-x_q)\, \sum_{i=1}^n \frac{\prod\limits_{j=1}^{n-1}
(x_i-y_j)}{\prod\limits_{j=1\atop j\ne i}^n (x_i-x_j)}
$$

\noindent restricted to the hyperplane $x_p=x_q$. The restriction
equals
$$
\frac{\prod\limits_{j=1}^{n-1} (x_p-y_j)}{\prod\limits_{j=1\atop
{j\ne p \atop j\ne q}}^n (x_p-x_j)} -
\frac{\prod\limits_{j=1}^{n-1} (x_q-y_j)} {\prod\limits_{j=1\atop
{j\ne q\atop j\ne p}}^n (x_q-x_j)}=0.
$$

 For identity (\ref{equation:hg_id_13}), the same technique works for
all the forms in the denominators except for $\alpha=x_p+x_m$
where $1\le p\le m-1$. If $\alpha=x_p+x_m$, then the restriction
of $\alpha\, L(x,y)$ to the hyperplane $\alpha=0$ is
$$
-\sum_{r=1}^{m-1} \frac{\prod\limits_{k=1}^{m-1} (x_k+y_r)}
{\prod\limits_{k=1\atop k\ne r}^{m-1} (y_r-y_k)}\,
\frac{\prod\limits_{k=1\atop k\ne r}^m
(x_p^2-y_k^2)}{\prod\limits_{k=1\atop k\ne p}^{m-1}
(x_p^2-x_k^2)}\, \frac{1}{y_1^2-x_p^2}+ \frac{y_m-x_p}{y_1-x_p}\,
\frac{\prod\limits_{k=2}^m (x_p+y_k)}{\prod\limits_{k=1\atop k\ne
p}^{m-1} (x_p-x_k)}.
$$

\noindent The fact that this restriction equals zero is equivalent
to the identity

\begin{equation}
\label{equation:hg_id_14} \sum\limits_{r=1}^m
\displaystyle{\frac{\prod\limits_{k=1\atop k\ne p}^m
(y_r+x_k)}{\prod\limits_{k=1\atop k\ne r}^m
(y_r-y_k)}}\,\displaystyle{ \frac{\prod\limits_{k=1\atop k\ne r}^m
(x_p-y_k)}{\prod\limits_{k=1\atop k\ne p}^m (x_p+x_k)}}=1
\end{equation}

\noindent which is similar to identity (\ref{equation:hg_id_5}),
but different from it. In order to prove
(\ref{equation:hg_id_14}), it is convenient to rewrite it as
$$
\sum_{r=1}^m \frac{\prod\limits_{k=1\atop k\ne p}^m (y_r+x_k)}
{\prod\limits_{k=1\atop k\ne r}^m (y_r-y_k)}\,
\frac{1}{x_p-y_r}-\frac{\prod\limits_{k=1\atop k\ne p}^m
(x_p+x_k)} {\prod\limits_{k=1}^m (x_p-y_k)}=0
$$

\noindent and use the same technique over again. \\

 The second step in the proofs is to show that a  polynomial $L(x,y)$ equals the
corresponding polynomial $R(x,y)$. Let us, for example, consider
(\ref{equation:hg_id_10}). In this case, $L(x)$ and $R(x)$ are
symmetric homogeneous polynomials in $x$ of degree $2$. The space
of such polynomials is two-dimensional. It is spanned by
$s_2=x_1^2+\cdots +x_n^2$ and $s_1^2$, where $s_1=x_1+\cdots
+x_n$. To prove that $L(x)\equiv R(x)$, we have to find two
linearly independent functionals $f_1$ and $f_2$ on this space
such that $f_i(L)=f_i(R)$ for $i=1,2$. We will treat the cases
$n=2k$ and
$n=2k+1$ separately. \\

 Let $n=2k$. Let $p_1=(-k,-k+1,\cdots ,-1,1,2,\cdots ,k)$ and
$p_2=(-k+1,-k+2,\cdots ,-1,1,2,\cdots ,$\\
$k+1)$. For a symmetric homogeneous polynomial $s$ of degree $2$,
let $f_i(s)=s(p_i)$, where $i=1,2$. Then

$$
\left| \begin{array}{cc}
f_1(s_2) & f_1(s_1^2) \\
 & \\
f_2(s_2) & f_2(s_1^2)
\end{array} \right| =
\left| \begin{array}{cc}
\frac{(2k+1)(k+1)k}{3} & 0 \\
 & \\
\frac{(2k-1)k(k-1)}{3}+k^2+(k+1)^2 & (2k+1)^2
\end{array} \right| = \frac{(2k+1)^3(k+1)k}{3}\ne 0.
$$

\noindent Thus, $f_1$ and $f_2$ are linearly independent. We have
$f_1(L)=L(p_1)=0=R(p_1)=f_1(R)$ and
$f_2(L)=L(p_2)=(2k+1)^2=R(p_2)=f_2(R)$. This finishes the proof
for $n=2k$. For $n=2k+1$, take $p_1=(-k,-k+1,\cdots ,k)$ and
$p_2=(-k+1,-k+2,\cdots ,k+1)$. The rest of the proof is analogous
to the case $n=2k$. Proofs of other identities are finished
similarly.


\medskip


\begin{thebibliography}{20}

\bibitem{BZ} A.~Berenstein, A.~Zelevinsky, {\it ``Triple Multiplicities for sl(r+1)
and the Spectrum of the Exterior Algebra of the Ajoint
Representation'',} J. of Algebraic Combinatorics 1 (1992), 7-22

\bibitem{BH} F.~Beukers, G.~Heckman, {\it ``Monodromy for the Hypergeometric Function
$\phantom{|}_nF_{n-1}$'',} Inventiones Mathematicae 95 (1989),
325-354

\bibitem{Bo} A.~Bolibrukh, {\it ``The 21st Hilbert Problem for Linear Fuchsian Systems'',}
proc. of the Steklov Institute of Mathematics, vol.206, AMS,
Providence, RI, 1995

\bibitem{CB1} W.~Crawley-Boevey, {\it ``On Matrices in Prescribed Conjugacy Classes with
No Common Invariant Subspace and Sum Zero'',} preprint
arXiv:math.LA/0103101, 15 March 2001

\bibitem{DR} M.~Dettweiler, S.~Reiter, {\it ``An Algorithm of Katz and its Applications
to the Inverse Galois Problem'',} J. Symbolic Computations (2000)
30, 761-798

\bibitem{DW1} H.~Derksen, J.~Weyman, {\it ``Semi-Invariants of Quivers and Saturation
Conjecture for Littlewood-Richardson Coefficients'',} JAMS,
vol.13, n.3 , 467-479

\bibitem{DW2} H.~Derksen, J.~Weyman, {\it ``On the $\sigma$-Stable Decomposition of
Quiver Representations'',} preprint (available on Derksen's
homepage)

\bibitem{F1} W.~Fulton, {\it ``Eigenvalues of Sums of Hermitian Matrices (after
A.~Klyachko)'',} S\'eminaire Bourbaki (1998)

\bibitem{F2} W.~Fulton, {\it ``Young Tableaux with Applications to Representation
Theory and Geometry'',} Cambridge University Press, 1997

\bibitem{GP} O.~Gleizer, A.~Postnikov, {\it ``Littlewood-Richardson Coefficients via
Yang-Baxter Equation'',} IMRN 14 (2000), 740-774

\bibitem{H1} Y.~Haraoka, {\it ``Canonical Forms of Differential Equations Free from
Accessory Parameters'',} SIAM J. of Mathematical Analysis, vol.25,
n.4 (1994), 1203-1226

\bibitem{H2} Y.~Haraoka, {\it ``Monodromy Representations of Systems of Differential Equations
Free from Accessory Parameters'',} SIAM J. of Mathematical
Analysis, vol.25, n.6 (1994), 1595-1621


\bibitem{VKac1} V.~Kac, {\it ``Infinite Root Systems, Representations of Graphs and
Invariant Theory'',} Inventiones Mathematicae 56 (1980), 57-92

\bibitem{VKac2} V.~Kac, {\it ``Infinite Root Systems, Representations of Graphs and
Invariant Theory II'',} Journal of Algebra 78 (1982), 141-162

\bibitem{NKatz} N.~Katz, {\it ``Rigid Local Systems''}, Annals of Mathematics
Studies, n.139, Princeton University Press (1996).

\bibitem{Kl} A.~Klyachko, {\it ``Stable Bundles, Representation Theory, and Hermitian Operators''},
Selecta Math. 4 (1998), 419-445

\bibitem{KT} A.~Knutson, T.~Tao, {\it ``The Honeycomb Model of $\Gl_n(\mathbb{C})$ Tensor
Products I: Proof of the Saturation Conjecture'',} JAMS, vol.12,
n.4 (1999), 1055-1090

\bibitem{KTW} A.~Knutson, T.~Tao, C.~Woodward, {\it ``The Honeycomb Model of $\Gl_n(\mathbb{C})$
Tensor Products II: Facets of the Littlewood-Richardson cone'',}
preprint

\bibitem{Ko1} V.~Kostov, {\it ``Monodromy Groups of Regular Systems on Riemann's
Sphere''}, prepublication N 401 de l'Universite de Nice, 1994; to
appear in Encyclopedia of Mathematical Sciences, Springer

\bibitem{Ko2} V.~Kostov, {\it ``On the Existence of Monodromy
Groups of Fuchsian systems on Riemann's Sphere with Unipotent
Generators'',} Journal of Dynamical and Control Systems, v.2, n.10
(1996) 125-155.

\bibitem{Ko3} V.~Kostov, {\it ``On the Deligne-Simpson problem''}, manuscrit, 47 p.
Electronic preprint math.AG/0011013. A para\^{\i}tre dans
Proceedings of the Steklov Institute, v. 238 (2002) , ``Monodromy
in problems of algebraic geometry and differential equations''.

\bibitem{Ko4} V.~Kostov, {\it ``On some aspects of the Deligne-Simpson problem''},
manuscrit, 48 p., \`a para\^{\i}tre dans un volume de ``Trudy
Seminara Arnol'da''. Electronic preprint math.AG/0005016.


\bibitem{Ko5} V.~Kostov, {\it ``On the Deligne-Simpson problem''}, C.R.A.S. Paris, t.
329, S\'er. I, pp. 657 -- 662, 1999

\bibitem{Ko6} V.~Kostov, {\it ``Some examples of rigid representations''}, Serdica
Mathematical Journal 26 (2000), p. 253 -- 276. Electronic preprint
math.AG/0006021.

\bibitem{Ko7} V.~Kostov, {\it ``Some examples related to the Deligne-Simpson problem''},
Proceedings of Second International Conference on Geometry,
Integrability and Quantization, St. Constantine and Elena,
Bulgaria, June 7-15, 2000, Edited by Ivailo M. Mladenov and
Gregory L. Naber, ISBN 954-90618-2-5, pp. 208 -- 227. Coral Press,
Sofia (2001) Electronic preprint math.AG/0011015.

\bibitem{Ko8} V.~Kostov, {\it ``The Deligne-Simpson problem for zero index of
rigidity''}, PERSPECTIVES IN COMPLEX ANALYSIS, DIFFERENTIAL
GEOMETRY AND MATHEMATICAL PHYSICS Editors: Stancho Dimiev and
Kouei Sekigawa (Proceedings of the 5th International Workshop on
Complex Structures and Vector Fields, St Constantin (Bulgaria),
3-9 September 2000), World Scientific, 2001, p. 1 -- 35.
Electronic preprint math.AG/0011107.

\bibitem{Ko9} V.~Kostov, {\it ``Examples illustrating some aspects of the
Deligne-Simpson problem''}, Serdica Mathematical Journal 27, No. 2
(2001), p. 143-158, Electronic preprint math.AG0101141.

\bibitem{MWZ} P.~Magyar, J.~Weyman, A.~Zelevinsky, {\it ``Multiple Flag Varieties of
Finite Type'',} Advances in Mathematics 141, 97-118 (1999)

\bibitem{O} K.~Okubo, {\it ``On the Group of Fuchsian
Equations'',} Seminar reports of Tokyo Metropolitan University,
1987.

\bibitem{R} C.~Ringel, {\it ``Exceptional Modules are Tree Modules'',} proc. of the
6-th conference of the International Linear Algebra Soc.
(Chermnitz, 1966), Linear Algebra Applications 275/276 (1998),
471-493

\bibitem{Sch1} A.~Schofield, {\it General Representations of Quivers,} Poc. London
Math. Soc. (3) 65 (1992), 46-64.

\bibitem{S} C.~Simpson, {\it ``Products of Matrices'',} AMS Proceedings 1 (1992),
157-185

\bibitem{V1} V.~Varadarjan, {\it ``Meromorphic Differential
Equations'',} Expositiones Mathematicae, vol.9, n.2 (1991),
97-188.

\bibitem{WW} E.~Whittaker, G.~Watson, {\it ``A Course of Modern Analysis'',}
Cambridge University Press, 1927

\bibitem{Y1} T.~Yokoyama, {\it ``On an Irreducibility Condition for Hypergeometric Systems'',}
Funk. Ekvac. 38, (1995), 11-19

\bibitem{Z} A.~Zelevinsky, {\it ``Littlewood-Richardson Semigroups'',} MSRI Preprint
n.1997-044.

\end{thebibliography}
\end{document}